\documentclass{article}

\usepackage{epic} 
\usepackage{eepic} 
\usepackage{latexsym}
\usepackage{amsthm}
\usepackage{amssymb}
\usepackage[all]{xy}
\usepackage{graphicx}
\usepackage{makeidx}
\makeindex

\providecommand{\noglossaryignore}[1]{}
\newcommand{\globalglossaryentry}[3]{\makebox[1.5in][l]{\tt $\backslash${#1}} 
\makebox[1.1in][l]{{$#2$}} \makebox[2.5in][l]{{#3}}\newline} 
\newcommand{\newcommandabbreviation}[3]{\newcommand{#1}{#2}%
\noglossaryignore{\globalglossaryentry{#3}{#2}{}}}
\newcommand{\renewcommandabbreviation}[3]{\renewcommand{#1}{#2}%
\noglossaryignore{\globalglossaryentry{#3}{#2}{}}}
\newcommand{\newcommandmacro}[4]{\newcommand{#1}{#2}%
\noglossaryignore{\globalglossaryentry{#3}{#2}{#4}}}
\newcommand{\gge}[3]{\noglossaryignore{\globalglossaryentry{#1}{#2}{#3}}}
\newcommand{\myaddress}%
{\parbox{3in}{\footnotesize \begin{center} 
Mathematics Department, City University, \\  
Northampton Square, London EC1V 0HB, UK.\end{center}}}
    
\newcounter{minidef}[section]
\renewcommand{\theminidef}{\thesection.\arabic{minidef}}
\newcommand{\mdef}{\refstepcounter{minidef} 
\medskip \noindent ({\bf \theminidef}) }
\newcounter{minicapt}

\newcounter{perchapter}
\renewcommand{\theperchapter}{\arabic{perchapter}}

\newtheoremstyle{puuu}
{7pt}%
{7pt}%
{\it}
{}
{}
{.}
{ }
{\thmnumber{({\bf #2}) }\thmname{\textbf{#1}}\thmnote{#3}}
\newtheoremstyle{puuuu}
{3pt}%
{3pt}%
{\rm}
{}
{}
{.}
{ }
{\thmnumber{({\bf #2}) }\thmname{\textsc{#1}}\thmnote{#3}}

\theoremstyle{puuu}
\newtheorem{pr}[minidef]{Proposition}
\newtheorem{lem}[minidef]{Lemma}
\newtheorem{theo}[minidef]{Theorem}

\newenvironment{cor}{\refstepcounter{perchapter}
\par \noindent
{\bf Corollary \theperchapter.}\it}

\newenvironment{de}{\refstepcounter{perchapter}
\par \noindent
{\bf Definition \theperchapter.}}

\newenvironment{rem}{\refstepcounter{perchapter}
\par \noindent
{\bf Remark \theperchapter.}}

\newenvironment{exa}{\refstepcounter{perchapter}
\par \noindent
{\bf Example \theperchapter.}}

\noglossaryignore{GREEK ETC.\newline}
\newcommandabbreviation{\e}{\epsilon}{e}        
\newcommandabbreviation{\lam}{\lambda}{lam}  
\newcommandabbreviation{\la}{\langle}{la}        
\newcommandabbreviation{\ran}{\rangle}{ran}
\newcommandabbreviation{\ha}{\#}{ha}             
\newcommandabbreviation{\rmap}{\rightarrow}{rmap}
\newcommandabbreviation{\aaa}{\alpha}{aaa}        
\newcommandabbreviation{\ab}{\alpha,\beta}{ab}
\newcommandabbreviation{\aab}{a(\ab )}{aab}       
\noglossaryignore{\newline RINGS\newline}
\newcommandabbreviation{\HH}{H \!\!\! I}{HH}               
\newcommandabbreviation{\C}{\mathbb C}{C}
\newcommandabbreviation{\N}{\mathbb N}{N}   
\newcommandabbreviation{\Z}{\mathbb Z}{Z}      
\renewcommandabbreviation{\Re}{\mathbb R}{Re}
\newcommandabbreviation{\R}{{\mathbb R}}{R}
\newcommandabbreviation{\Q}{\mathbb Q }{Q}
\renewcommandabbreviation{\H}{\mathbb H }{H}
\noglossaryignore{\newline SYMMETRIC GROUP\newline}
\def\Sym(#1){\Sigma(#1)}                   
\gge{Sym(-)}{\Sym(-)}{symmetric group on - objects}
\def\Sy(#1){\Sigma_{#1}}                   
\gge{Sy(-)}{\Sy(-)}{symmetric group irreducible -}
\def\sym(#1){\mbox{\LARGE s}(#1)}        
\gge{sym(-)}{\sym(-)}{symmetric group on - objects (variant)}
\def\sy(#1){\mbox{\LARGE s}({#1})}        
\gge{sy(-)}{\sy(-)}{symmetric group irreducible - (variant)}
\newcommandmacro{\cs}{\C \, \sy(n)}{cs}{symmetric group algebra over $\C$}
\noglossaryignore{\newline PARTITIONS/SETS\newline}
\newcommand{\Nset}[1]{\underline{#1}}
\gge{Nset\{-\}}{\Nset{-}}{set of natural numbers to -} 
\def\nset(#1){ \{ #1 \}_{ \underline{n} }} 
\gge{nset(-)}{\nset(-)}{a set $-\times\Nset$}
\def\ul(#1){_{\underline{#1}}}             
\gge{ul(-)}{{}\ul(-)}{subscript underline -}
\def\Ee(#1){{\bf E}_{#1}}                  
\gge{Ee(-)}{\Ee(-)}{set of equivalence relations on set -}
\def\Eee(#1){{\bf E}_{\{ #1 \}_{\underline{n}}}}   
\gge{Eee(-)}{\Eee(-)}{ditto for nset}
\def\Een(#1,#2){{\bf E}_{\{ #1 \}_{\underline{#2}}}}   
\def\Ssn(#1,#2){{\bf S}_{\{ #1 \}_{\underline{#2}}}}   
\def\Ss(#1){{\bf S}_{#1}}                  
\def\Sss(#1){{\bf S}_{\{ #1 \}_{\underline{n}}}}   
\def\bbc(#1){((\beta_1)(\beta_2)...(\beta_{#1}))}      
\newcommandmacro{\Ln}{{\Gamma}^{n}}{Ln}{large index set}
\newcommandmacro{\LnQ}{{\Gamma}^{n}_Q}{LnQ}{index set}
\newcommandmacro{\Zz}{\zeta}{Zz}{`shape' function}
\noglossaryignore{\newline PARTITION ALGEBRA\newline}
\def\ka(#1){\kappa_{#1}}                   
\def\Sm(#1){\Sigma_{#1}}                   
\newcommandmacro{\com}{\bullet}{com}{bullet composition}
\newcommandmacro{\enm}{\; e^n(\! m\! ) \;}{enm}{product of idempotents}
\def\Ai(#1){ A^{ #1 \cdot } }              
\def\Aij(#1,#2){ A^{ #1  #2 } }            
\newcommandmacro{\One}{\mbox{\bf $1 \!\!\! 1$}}{One}{algebra unit 1}
\newcommandmacro{\Bp}{B_p}{Bp}{partition basis}
\def\Bb(#1){B_p[#1]}                       
\def\Pp(#1){P_n[#1]}                       
\def\Ps(#1){P_n[#1] \! /}                  
\newcommandmacro{\Ph}{\hat{P}}{Ph}{P hat  algebra}
\def\Is(#1){\sim^{#1}}                     
\noglossaryignore{\newline MODULES\newline}
\def\Wm(#1){{\cal S}_{#1}}                 
\gge{Wm(-)}{\Wm(-)}{Weyl module with index -}
\def\wm(#1,#2){{}_{#1}{\cal S}_{#2}}       
\gge{wm(-1,-)}{\wm(-1,-)}{Weyl module with index -}
\def\Ind(#1,#2,#3){\mbox{Ind}_{#1}^{#2}#3} 
\gge{Ind(-1,-2,-)}{\Ind(-1,-2,-)}{induction}
\def\Res(#1,#2,#3){\mbox{Res}_{#1}^{#2}#3} 
\gge{Res(-1,-2,-)}{\Res(-1,-2,-)}{restriction}
\newcommandabbreviation{\weyl}{standard}{weyl}
\newcommandabbreviation{\head}{\mbox{head }}{head}
\newcommandabbreviation{\Weyl}{Weyl}{Weyl}
\def\SS(#1){{\cal S}_{#1}}                 
\gge{SS(-)}{\SS(-)}{Specht/Weyl module index -}
\def\LL(#1){{\cal L}_{#1}}                 
\gge{LL(-)}{\LL(-)}{simple module index -}
\newcommandabbreviation{\mod}{\mbox{mod}}{mod}
\newcommand{\modl}[1]{\mbox{$#1$-}\mod}     
\gge{modl\{-\}}{\modl{-}}{set of left modules for ring -}
\newcommand{\modr}[1]{\mod\mbox{-$#1$}}     
\gge{modr\{-\}}{\modr{-}}{set of right modules for ring -}
\noglossaryignore{\newline FUNCTORS/MAPS\newline}
\newcommandmacro{\Gg}{{\cal G}}{Gg}{G Functor}
\newcommandmacro{\Fg}{{\cal F}}{Fg}{F Functor}
\newcommandmacro{\ra}{\rightarrow}{ra}{}
\def\ses(#1,#2,#3){0\ra #1 \ra #2 \ra #3 \ra 0}   
\gge{ses(1,2,-)}{\ses(1,2,-)}{\hspace{.5in} short exact sequence}
\def\starr(#1){ \stackrel{ #1 }{\longrightarrow} }
\gge{starr(-)}{\starr(-)}{}
\newcommandmacro{\doublerightarrow}{\; -\!\!\! -\!\!\!\!\!\! \gg \;}
{doublerightarrow}{}
\noglossaryignore{\newline PARTITION ALGEBRA MAPS\newline}
\newcommandmacro{\smap}{s}{smap}{`inclusion' map}
\newcommandmacro{\tmap}{t}{tmap}{$ P_n -> S_n$}
\newcommandmacro{\pmap}{\psi}{pmap}{$ S_n -> P_n $}
\noglossaryignore{\newline MISC.\newline}
\def\Amap(#1){{\cal A}_{#1}}               
\gge{Amap(-)}{\Amap(-)}{}
\def\Rr(#1){R_{#1}}                        
\gge{Rr(-)}{\Rr(-)}{restriction of E}
\def\Cr(#1){C_{#1}}                        
\gge{Cr(-)}{\Cr(-)}{restriction of E to N}
\newcommandmacro{\Tm}{{\cal T}}{Tm}{Transfer Matrix}
\def\On(#1){{\cal I}_{#1}}
\gge{On(-)}{\On(-)}{}
\newcommandmacro{\UU}{\underline{\sqcup}}{UU}{}  
\newcommandmacro{\UUU}{\sqcup}{UUU}{}  
\newcommandmacro{\Vq}{V_Q^{\otimes n}}{Vq}{Potts config. space}
\def\bs(#1,#2){\mbox{{\Large $\ast$}}^{#1}_{#2}}  
\gge{bs(-,-)}{\bs(-,-)}{general plumbing multiplier}
\newcommand{\ignore}[1]{}
\gge{ignore\{-\}}{\ignore{-}}{ignore argument!}
\def\choo(#1,#2){ \left( \begin{array}{c} #1 \\ #2 \end{array} \right) } 
\gge{choo(-1,-)}{\choo(-1,-)}{choose}
\newcommand{\Qed}{$\Box$}
\gge{Qed}{\mbox{\Qed}}{QED}
\def\staq(#1){\stackrel{#1}{=}}            
\gge{staq(-)}{\staq(-)}{}
\def\stam(#1){\stackrel{#1}{\rightarrow}}  
\gge{stam(-)}{\stam(-)}{}
\def\mat{ \left( \begin{array} }    
\def\tam{ \end{array}  \right) }
\gge{mat/tam}{...}{matrix delimiters}
\newcommand{\beq}{\begin{equation} }
\def\eql(#1){ \begin{equation} \label{#1} 
}
\newcommand{\eq}{\end{equation} }
\def\eqal(#1){\begin{eqnarray} \label{#1} }
\def\eqa{\end{eqnarray} }
\def\lab(#1){\label{#1}
}
\def\prl(#1){ \begin{pr} \label{#1} 
}

\def\theol(#1){ \begin{theo} \label{#1} 
}

\def\leml(#1){ \begin{lem} \label{#1} 
}

\def\corl(#1){ \begin{cor} \label{#1} 
}

\def\del(#1){ \begin{de} \label{#1} 
}

\def\reml(#1){ \begin{rem} \label{#1} 
}

\def\exal(#1){ \begin{exa} \label{#1} 
}

\gge{smeq\{-\}}{...}{small equation}
\gge{fneq\{-\}}{...}{very small equation}
\newcommand{\beqa}{\begin{eqnarray}}%
\newcommand{\eeqa}{\end{eqnarray}}%

\noglossaryignore{\newline HECKE/BLOB\newline}
\newcommandmacro{\Hnq}{H_n(q)}{Hnq}{ * freestanding symbol}
\newcommandmacro{\Hn}{H_n}{Hn}{      *-mod etc.}
\newcommandmacro{\A}{{\cal A}}{A}{}
\newcommandmacro{\Cwts}{C}{Cwts}{}
\newcommandmacro{\CA}{{\cal A}}{CA}{}

\newcommandmacro{\calA}{{\cal A}}{calA}{}
\newcommandmacro{\modi}{\mbox{Mod} }{modi}{was mod not modi!}
\newcommandmacro{\Wgen}{{\Bbb S}}{Wgen}{}
\def\ol(#1){\overline{#1}}
\newcommandmacro{\st}{\mbox{St}}{st}{}
\def\CMult(#1,#2){(#1:#2)}
\def\CM(#1,#2){( #1 : #2 )}
\def\FMult#1,#2{(#1:#2)}
\def\CF#1,#2{(#1:#2)}

\newcommandmacro{\Top}{\mbox{Top}}{Top}{}
\newcommandmacro{\Soc}{\mbox{Soc}}{Soc}{}
\newcommandmacro{\Head}{\mbox{Head}}{Head}{}
\newcommandmacro{\Filt}{{\cal F}}{Filt}{}
\newcommandmacro{\Mod}{\mbox{mod}}{Mod}{}
\newcommandmacro{\Resi}{\mbox{Res }}{Resi}{was without i!}
\newcommandmacro{\Indi}{\mbox{Ind }}{Indi}{was without i!}
\def\RR(#1,#2){R^{#1}_{#2}}   
\def\TT(#1,#2){T^{#1}_{#2}}   
\def\implies{\Rightarrow}

\def\Hom{\mbox{Hom}}
\def\bigplus{\mbox{\LARGE $+$}}

\newcommandmacro{\Ann}{\mbox{Ann}}{Ann}{}
\newcommandmacro{\Cen}{\mbox{Cen}}{Cen}{}
\newcommandmacro{\End}{\mbox{End}}{End}{}
\newcommandabbreviation{\semisimple}{semisimple}{semisimple}
\newcommandabbreviation{\Bratteli}{Bratteli}{Bratteli}
\newcommandabbreviation{\JBC}{Jones Basic Construction}{JBC}
\newcommandabbreviation{\pa}{partition algebra}{pa}
\newcommandabbreviation{\TM}{transfer matrix}{TM}
\newcommandabbreviation{\PM}{Potts model}{PM}
\newcommandabbreviation{\QSC}{quantum spin chain}{QSC}
\newcommandabbreviation{\Hamiltonian}{Hamiltonian}{Hamiltonian}
\newcommandabbreviation{\YS}{Young symmetrizer}{YS}

\newcommand{\mystufffont}{\textsc} 
  
\newtheoremstyle{pu}
{7pt}%
{7pt}%
{\it}
{}
{}
{.}
{ }
{\thmnumber{({\bf #2}) }\thmname{\textbf{#1}}\thmnote{#3}}
\newtheoremstyle{puu}
{3pt}%
{3pt}%
{\rm}
{}
{}
{.}
{ }
{\thmnumber{({\bf #2}) }\thmname{\textsc{#1}}\thmnote{#3}}

\theoremstyle{pu}
\newtheorem{mmpr}[minidef]{Proposition}
\newtheorem{mlem}[minidef]{Lemma}
\newtheorem{mth}[minidef]{Theorem}
\newcommand{\mpr}[1]{\begin{mmpr} #1 \end{mmpr}}

\theoremstyle{puu}

\newtheorem{mrem}[minidef]{Remark}

\newcommand{\murem}{\noindent{\mystufffont{Remark.}} }

\title{The decomposition matrices of the Brauer algebra \\ 
over the complex field 
}
\author{Paul P Martin}
\date{}
\begin{document}
\maketitle
\setlength{\parskip}{1ex}



\newcommand{\Grp}{{\mathbf{ Grp}}}
\newcommand{\Mat}{{\mathbf{ Mat}}}
\newcommand{\Fld}{{\mathbf{ Fld}}}
\newcommand{\Rng}{{\mathbf{ Rng}}}
\newcommand{\Br}{{\mathbf{ Br}}}
\newcommand{\specht}{{\mathcal S}}  
\newcommand{\res}{\mbox{Res}} 
\newcommand{\BS}{Brauer-Specht}

\newcommand{\CDMI}{\cite{CoxDevisscherMartin05}}
\newcommand{\CoDMI}{CoxDevisscherMartin05}
\newcommand{\CDMi}{CoxDevisscherMartin05}
\newcommand{\CDMII}{\cite{CoxDevisscherMartin06}}
\newcommand{\CDMii}{CoxDevisscherMartin06}
\newcommand{\CDMIII}{\cite{CoxDevisscherMartin08}}
\newcommand{\CDMiii}{CoxDevisscherMartin08}
\newcommand{\QQ}{\delta}
\newcommand{\zring}{{\scriptscriptstyle \Z}}
\newcommand{\Psio}{o} 
\newcommand{\bb}{{\mathbf b}} 
\newcommand{\Psiob}{\bb} 
\newcommand{\Pow}{P}  
\newcommand{\Gammad}{{\Gamma_{\!\delta}}} 
\newcommand{\proj}{\mbox{Proj}} 
\newcommand{\valley}{valley}
\newcommand{\sis}{p} 
\newcommand{\sing}{s}
\newcommand{\Dgroup}{{\mathcal D}}
\newcommand{\Dgroupp}{\Dgroup_{+}}
\newcommand{\Dmat}{D}
\newcommand{\Cmat}{C}
\newcommand{\Reg}{Reg}
\newcommand{\Ap}{A^+} 
\newcommand{\ed}{e_{\delta}}
\newcommand{\ReRB}{rim-end removable box}
\newcommand{\Deltan}{\Delta_}
\newcommand{\Deltab}{\Delta} 
\newcommand{\Db}[2]{\Deltab_{#1}(#2)} 
\newcommand{\Dbt}[2]{\Deltab_{#1}(#2)'} 
\newcommand{\Pn}{P_} 
\newcommand{\Pb}[2]{P_{#1}(#2)} 
\newcommand{\Pbt}[2]{P_{#1}(#2)'} 
\newcommand{\Lb}[2]{L_{#1}(#2)} 
\newcommand{\Lbt}[2]{L_{#1}(#2)'} 
\newcommand{\Gal}{G_{alc}}    
\newcommand{\Ge}{G_{even}} 
\newcommand{\Go}{G_{odd}}  
\newcommand{\oonn}{(1)+(1)}
\newcommand{\mib}[2]{#1/#2#1}
\newcommand{\adjac}{\;{}^{\triangleleft}_{\triangleright}\;}
\newcommand{\StDM}[3]{( #1({#2}) : \Db{}{#3} )}
\newcommand{\StDMt}[3]{( #1({#2})' : \Dbt{}{#3} )}
\newcommand{\SDM}[2]{[ #1 : #2 ]}       
\newcommand{\LbDM}[2]{\SDM{#1}{\Lb{}{#2}}}

\newcommand{\lal}{\lambda / \alpha\lambda} 
\newcommand{\RST}{reflexive-symmetric-transitive}
\newcommand{\Td}{d} 
\newcommand{\TLd}{TL-diagram}
\newcommand{\onetwo}{\underline{12}} 
\newcommand{\KL}{Kazhdan--Lusztig}
\newcommand{\pKLp}{parabolic \KL\ polynomial}

\newcommand{\lambdaa}{a}

\newcommand{\Hyp}{{\mathbb H}}

\newcommand{\footnot}[1]{}
\newcommand{\fotnote}[1]{}

\section{Introduction}

For each field $k$, natural number $n$ and parameter
$\delta \in k$, the Brauer algebra $B_n(\delta)$ is a finite
dimensional algebra, 
with a basis of pair partitions of the set $\{1,2,...,2n \}$
\cite{Brauer37}.
Indeed  there is a $\Z[\delta]$-algebra $B_n^{\zring}$
(for $\QQ$ indeterminate), free of finite rank as a 
 $\Z[\delta]$-module, that passes to each Brauer algebra by 
the natural base change;
and a collection of modules  
$\{ \Delta^{\zring}(\lambda) \}_{\lambda\in\Lambda^n}$ 
for this algebra that are
 $\Z[\delta]$-free modules of known rank,  
so that
\[
\Delta^k(\lambda) = k \otimes_{\Z[\delta]} \Delta^{\zring}(\lambda)
\]
are $B_n(\QQ)$-modules, 
and that there is a choice of field $k$ extending $\Z[\delta]$ for which 
$\{ \Delta^k(\lambda) \}_{\lambda}$ is a complete set of simple modules.
Accordingly
we are presented with the following tasks 
in studying the representation theory of $B_n(\QQ)$:
\\
(1) There are finitely many isomorphism 
classes of simple modules --- index these.
\\
(2) Describe the blocks (the \RST\ closure of the relation on the index
set for simples given by ${\lambda} \sim {\mu}$ 
if simple modules 
$\Lb{}{\lambda} $ and $ \Lb{}{\mu}$ 
are composition factors of the same
indecomposable projective module). 
\\
(3) Describe the composition multiplicities of indecomposable
projective modules
(which follow from the composition multiplicities for the 
$\Delta^k(\lambda)$ (see  for example 
\cite[\S 16]{CurtisReiner90},\cite[\S 1.9]{Benson95})).

Over the complex field, (1) was effectively solved in \cite{Brown55},
and (2) in \CDMI\ (see references therein for other important
contributions). 
Here we solve (3).

\medskip


\newcommand{\Lambdand}{\Lambda^{n,\QQ}}

The layout of the paper is as follows.
For each $n,\delta$
we wish to compute the Cartan decomposition matrix $\Cmat$ given by
$
\Cmat_{\lambda \mu} = \SDM{\Pb{}{\lambda}}{\Lb{}{\mu}}
$
where $\{ \Pb{}{\lambda} \}_{\lambda\in\Lambdand}$
and $\{ \Lb{}{\lambda} \}_{\lambda\in\Lambdand}$
are complete sets of imdecomposable projective and simple modules
respectively. 
We firstly recall some organisational results to this end.
We construct the modules $\Db{}{\lambda}$, such that
projective modules are filtered by these, with well-defined
composition multiplicities denoted $\StDM{P}{\lambda}{\mu}$;
and that $\Cmat=\Dmat \Dmat^T$, where 
$\Dmat_{\lambda,\mu} = \StDM{P}{\lambda}{\mu} 
= \SDM{\Db{}{\mu}}{\Lb{}{\lambda}}$
(what might be called the 
$\Deltab$-decomposition matrix).
Then we construct an inverse limit for the sets $\{ \Lambdand \}_n$
and show that the Cartan decomposition matrices (and the
$\Dmat$s) for all $n$ can be obtained by projection from a
corresponding limit.

Next we give an explicit matrix $\Dmat$ for each $\QQ$
(this construction takes up the majority of the paper). 
And finally we prove, in Section~\ref{ss:main th}, that it is the
limit 
$\Deltab$-decomposition matrix. 


It is probably helpful to note that 
the original route to the solution of the problem 
was slightly different. It 
proceeded from a conjecture, following \cite[\S1.2]{MartinWoodcock03}, 
that $\Dmat$ would consist of 
evaluations of \pKLp s 
for a certain reflection group given in, and parabolic 
determined 
by, our joint work in \CDMII.
This is essentially correct, as it turns out,
and without this idea we would not have had
a candidate for $\Dmat$, the form of which then drives the proof of the
Theorem. 
However the proof does not, in the end, lie entirely within the
realms of  Kazhdan-Lusztig theory and alcove geometry. 
Accordingly we do not use this
framework, but instead a more general one within which the proof
proceeds uniformly. With regard to the alcove geometry we restrict
ourselves to incorporating some key ideas; and beyond that just
a few remarks, where it seems helpful to explain
strategy.

We  return to discuss our parabolic Kazhdan-Lusztig
polynomial solution  
in a second part to the paper: section~\ref{ss:pKLp} and thereafter.


As the derivation of our main result is somewhat involved, we end here
with a brief preview of the result itself.
For each fixed  
$\QQ\in\Z$, the rows and columns of the limit
$\Delta$-decomposition matrix $D$ may be indexed by $\Lambda$,
the set of all integer partitions. 
This matrix may be decomposed, of course, as a direct sum of matrices for the
limit blocks. 
In this sense we may describe the blocks by a partition of $\Lambda$. 
As we shall see, there is a map for each block to the set 
$\Pow_{even}(\N)$ of subsets of $\N$ of even degree.
Under these maps all the block summands of $D$ (and for all $\QQ$)
are identified with the same matrix. Thus we require only to give 
a closed form for the entries of this
matrix. 
The closed form is given in Section~\ref{ss:decomp data},  
but an indication of its structure 
is given by a
truncation to a suitable finite rank.
Such a truncation is given  in Figure~\ref{fig: big pKL}
(the entries in this matrix  
encode polynomials that will be used
later, and which must be evaluated at 1 to give the decomposition numbers;
the {\em blank} entries evaluate to zero, and all other 
entries evaluate to 1). 

This paper is a contribution toward 
a larger project, with Cox and De Visscher, aiming to
compute the decomposition matrices of the Brauer algebras over fields
of finite characteristic. This is a very much harder problem again
(it includes
the representation theory of the symmetric groups over the same fields
as a sub-datum --- see \CDMII), and so it is appropriate to present 
the characteristic zero case separately.

\section{Brauer diagrams and Brauer algebras}

We mainly base our exposition on 
the notations and terminology of \CDMI,
as well as key results from that paper.
 For 
self-containedness,  however, we review the  notation here.
Our hypotheses are slightly more general than in \CDMI,
however many of the proofs in \CDMI\ go through
essentially unchanged (as we shall
indicate, where appropriate). 
We shall also  
make use of a  categorical formulation of the Brauer
algebra (a subcategory of 
the partition algebra category of
\cite[\S7]{Martin94}).


\mdef For $n\in \N$
we write $S_n$ for the symmetric group, and 
$
\Nset{n} \; := \; \{1,2,..,n \}
$ 
and 
$\Nset{n'} \; := \; \{1',2',..,n' \}$ (and so on).
For $S$ a set we write 
$P(S)$ for the power set and 
$J_S$ for the set of pair-partitions of $S$.
We define $J_{n,m} = J_{\Nset{n} \cup \Nset{m'}}$.
For example, in $J_{n,n}$  let us define
\eql(Uij)
U_{ij} = \{ \{1,1' \}, \{2,2'\},...,\{i,j\},\{i',j'\},...,\{n,n'\}\}
\eq
\[
(ij) =  \{ \{1,1' \}, \{2,2'\},...,\{i,j'\},\{i',j\},...,\{n,n'\}\}
\]


\mdef An $(n,m)$-Brauer diagram is a representation of a pair
partition of a row of $n$ and a row of $m$ vertices, arranged on the
top and bottom edges (respectively) of a rectangular frame. 
Each part is drawn as a line, joining the corresponding pair 
of vertices, in the rectangular interval.
We identify two diagrams if they represent the same partition.
It will be evident that these diagrams can be used to describe 
elements of $J_{n,m}$. 
For example, from $J_{6,6}$:
\[
U_{24}  \;\; = \;\;\;  \raisebox{-.43in}{
\includegraphics[width=1.52in]{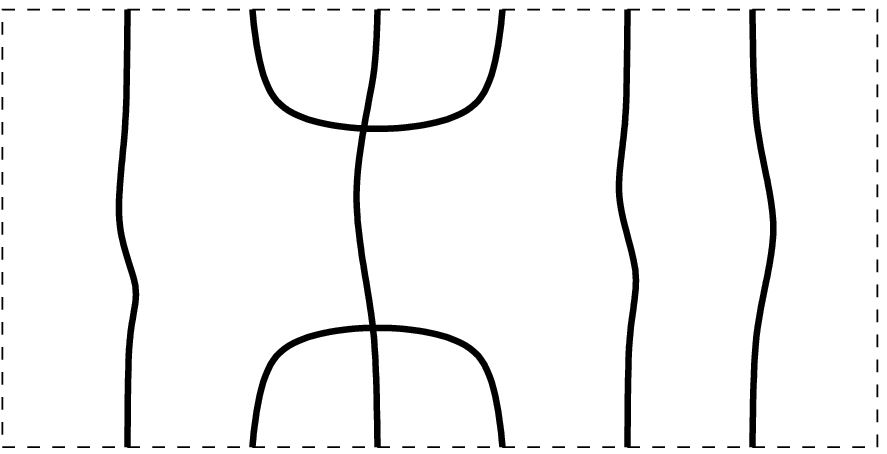} } 
\]

We then define a map
\[
J_{n,m} \times J_{m,l} \rightarrow \N_0 \times J_{n,l} 
\rightarrow \Z[\QQ] J_{n,l}
\]
as follows. Firstly juxtapose the diagrams so that the two sets of $m$
vertices meet.
This produces a diagram for an element $d$ of $J_{n,l}$ 
(the pair
partition of the vertices on the exterior of the combined frame);
together with some number $c$ of closed loops.
The final image is then $\QQ^c d$. 

We denote this composition by  $\circ$.
For $k$ a ring and $\QQ \in k$
we 
have a $k$-linear  
category 
with object set $\N_0$, hom-sets $\{ kJ_{n,m} \}_{n,m \in \N_0}$,
and composition $k$-linearly extending $\circ$.
We denote this category by
$\Br^k_{\delta}$, or just  $\Br^{}_{\delta}$ if
$k$ is fixed. 
(Here we allow $k=\Z[\QQ]$ or any suitable base change.)


\newcommand{\homB}{\hom_B} 
\newcommand{\homBr}{\Br}
\newcommand{\homBro}[1]{\Br^{\underline{#1}}}

\mdef 
Write  
$\homBr(m,n)$
for the set of $(m,n)$-Brauer diagrams;
$\homBr^{\leq l}(m,n)$ for the subset with $\leq l$ propagating lines; 
 $\homBr^l(m,n)$ for the subset with $l$ propagating lines; 
and  $\homBro{l}(m,n)$ for the subset of these in which
none of the $l$ propagating lines cross.
Write  $1_r$ for the identity diagram in $\homBr(r,r)$. 

Note that the category composition defines a bijection:
\eql(1_lxl,l) \label{1_lx(l,l)}
\homBro{l}(m,l) \; \times \; \homBr^{l}(l,l) \; \rightarrow \; \homBr^l(m,l)
\eq

Define a product
\[
 \otimes: \homBr(m,n) \times \homBr(r,s) \; \rightarrow \homBr(m+r,n+s)
\]
by placing diagrams side by side.
Hence define an injection
adding propagating lines $\{\{m+1,n+1' \},...,\{m+r,n+r' \}\}$:
\[
i_{m+1,m+r}: \homBr(m,n) \hookrightarrow \homBr(m+r,n+r)  
\] \[ \qquad
D \; \mapsto \;\; D \otimes 1_r
\]

The Brauer algebra $B_n(\delta)$ over $k$ is the free $k$-module 
with basis $\homBr(n,n)$ and the category composition
(i.e. 
replacing each closed loop formed in composition by a factor $\delta$).


\mrem{
The fully `integral' version is the case $k = \Z[\delta]$. 
From here there are thus two aspects to the base change to a field:
the choice of $k$ and the choice of $\delta$.
More precisely this is the choice of $k$ equipped with the structure
of $\Z[\delta]$-algebra.
Thus we have possible intermediate steps: base change to $k[\delta]$
($k$ a field);
base change to $\Z$ (a $\Z[\delta]$-algebra by fixing 
$\delta=d\in \Z$). 
Each of these ground rings is a  
principle ideal domain and hence a Dedekind domain,
and hence amenable to a $P$-modular treatment
(see  for example 
\cite[\S 16]{CurtisReiner90},\cite{Benson95}).}



\section{\BS\ modules} 
Here we construct the integral representations (in the sense of
\cite{Benson95}) that we shall need. (These
base change, entirely transparently, 
to the {\em standard} modules of \CDMI.)

\mdef \label{de:regdec1}
For any ring $k$ and $\delta\in k$, 
we have,
as an elementary consequence of the composition rule, 
a sequence of $B_n(\delta)$-bimodules:
\eql(reg seq1)
k\homBr(n,n) = k\homBr^{\leq n}(n,n) 
  \supset  k\homBr^{\leq n-2}(n,n) \supset  k\homBr^{\leq n-4}(n,n)
  \supset ...\supset  k\homBr^{1/0}(n,n)
\eq
Note that the $i$-th section of the sequence~(\ref{reg seq1}) has basis 
$\homBr^{n-2i}(n,n)$. For $n-2i=l$ we have
\eql(reg dec2)
k \homBr^{l}(n,n) \; \cong \;\;
\bigoplus_{w \in \homBro{l}(l,n)} \;\; k\homBr^{l}(n,l) \; w
\eq
as a left module; where all the summands are isomorphic to 
$k\homBr^{l}(n,l)$.

Fixing a ring $k$, 
it will be evident that   $\homBr^l(m,l)$ is a basis for a
left-$B_m(\delta)$ right-$kS_l$ bimodule, where
the action on the left  is via the category composition,
 quotienting by $k\homBr^{\leq l-2}(m,l)$.


\mpr{ \label{pr:Brexact1}
Fix any ring $k$. The free $k$-module    $k\homBr^l(m,l)$ 
(which is a left $B_m(\QQ)$-module by the action in (\ref{de:regdec1}))
is 
a projective right $kS_l$-module; and hence the functor
\[
k\homBr^l(m,l) \otimes_{kS_l} -  \; : \;
      kS_l\! -\!\mod \;\; \rightarrow \;\;\; B_m(\delta)\! -\!\mod
\]
between the categories of left-modules
is exact.
}
\proof{ $k\homBr^l(m,l)$ is a direct sum of copies of the regular 
 right $kS_l$-module. 
\Qed}

\mdef \label{de:specht1}
Let $\Lambda_n = \{ \lambda \vdash n \}$,
the set of integer partitions of $n$. 
Let $\Lambda$ be the set of all integer partitions;
and 
\[
\Lambda^n = \Lambda_{n} \cup  \Lambda_{n-2} \cup \ldots \cup \Lambda_{0/1} 
\]
For $\lambda \vdash l$ let  $\specht(\lambda)$ denote the 
corresponding  $kS_l$-Specht module (see e.g. \cite{JamesKerber81}),
and define
\[
\Delta_m(\lambda) = k\homBr^l(m,l) \otimes_{kS_l} {\mathcal S}(\lambda)
\]
as the image of this Specht module  
under the  functor in (\ref{pr:Brexact1}). 
\\
We may write $\Delta_m^k(\lambda)$ for $\Delta_m(\lambda)$ 
if we wish to emphasise the ring, or
$\Delta_m^{\QQ}(\lambda)$ ($\QQ\in k$) if $k$ is fixed as a field,
to fix it as a $\Z[\QQ]$-algebra. On the other hand, where unambiguous
we may just write $\Delta_{}(\lambda)$. 
We shall adopt analogous conventions for projective and simple modules.


\mpr{ \label{pr:projDfilt1}
Fix $n$ and suppose $k$ is such that  left regular module  
${}_{k S_l}kS_l$ is filtered by 
$\{ {\mathcal S}(\lambda) \}_{\lambda\in\Lambda_l}$ 
for all $l \leq n$.
Then 
the left regular module ${}_{B_n}B_n$ is filtered by 
$\{ \Delta_n (\lambda) \}_{\lambda\in\Lambda^n}$.
In particular Brauer algebra projective modules
over $\C$ (any $\QQ$) are filtered by 
$\{ \Delta_n (\lambda) \}_{\lambda\in\Lambda^n}$.
}

\proof{Note first that if a module $M$ is filtered by a set $\{ N_i \}_i$,
and these are all filtered by a set  $\{ N'_j \}_j$,
then $M$ is filtered by  $\{ N'_j \}_j$.
By (\ref{de:regdec1}) the set $\{ k \Br^l(n,l) \}_l$ gives (via the
action therein) a left-$B_n$ filtration of $B_n$.
By Prop.~\ref{pr:Brexact1} each factor itself has a filtration by $\Delta$s
under the stated condition.
For the last part, simply note that $\C S_l$ is semisimple, and
each projective
$P_n(\lambda)$ a direct summand of  ${}_{B_n}B_n$.
\Qed
}


\mpr{ \label{pr:basis} 
{\rm \cite[Lemma~2.4]{\CDMi}}
Let $b(\lambda)$ be a basis for $ \specht(\lambda)$.
Then 
\[
b_{\Delta_m(\lambda)} = 
\{ a\otimes b \; : \; (a,b) \in \homBro{l}(m,l) \times b(\lambda) \}
\]
is a basis for $\Delta_m(\lambda)$. 
}
\proof{
This is a set of generators 
by (\ref{1_lx(l,l)}).
On the other hand this set passes to a basis (of the image) under
the surjective multiplication map
(using 
from 
\cite{JamesKerber81}
that $\specht(\lambda)$ 
is a left ideal), so it is $k$-free. \Qed
}
\medskip


\mdef \label{de:base1}
We mention explicitly the following low rank cases,
which form the bases for inductions later on.
We have $B_0(\delta) \cong B_1(\delta) \cong k$.
For $B_2(\delta)$ we have 
$\Delta_2(\emptyset)$, $\Delta_2(2)$, $\Delta_2(1^2)$,
each of rank 1. 
For $\delta=0$ we have, over $\C$,
\[
\Delta_2(2) \stackrel{\sim}{\rightarrow} \Delta_2(\emptyset)
\]
Thus we may regard  $\Delta_2(2)$, $\Delta_2(1^2)$ as the inequivalent
simple modules, and $P_2(2)$ is the self-extension of  $\Delta_2(2)$,
while $P_2(1^2)=\Delta_2(1^2)$. 


\subsection{Globalisation functors}
Here we 
define certain functors that will allow us,
in Section~\ref{ss:char}, to manipulate composition
muliplicity data for all $n$ simultaneously.

\mdef For $n+m$ even 
the $k$-module   
$k \homBr(n,m)$ is an algebra bimodule. Thus there is a functor
between left-module categories 
\[
k\homBr(n,m) \otimes_{B_m} \! - \; : \;
   B_m\! -\!\mod \;\; \rightarrow \;\; B_{n} \! -\!\mod
\]
Let us write $F$ for the functor 
$k\homBr(n-2,n) \otimes_{B_n} \! - \;$; and $G$ for the functor 
$k\homBr(n,n-2) \otimes_{B_{n-2}} \! - \;$.


\mpr{ \label{pr:Fexact1}
Suppose either $n>2$ or $\delta$ invertible in $k$.
Then \\ 
(I) the $k$-space $k\homBr(n-2,n)$ is projective as a right
$B_n$-module;
and indeed
\[
k\homBr(n-2,n) \; \cong \;\;\; e \; ( k\homBr(n,n) ) 
\]
as a right $B_n$-module,
for a suitable idempotent $e \in k\homBr(n,n)$
(see the proof for an explicit construction of $e$). 
\\
(II) Functor $F : B_n\! -\!\mod \rightarrow B_{n-2}\! -\!\mod$ is exact; 
$G$ is a right-exact right-inverse to $F$.
}
\proof{(We prove a left-handed version. The right-handed follows
  immediately.) 
As a left module
\[
B_3 \cong \bigoplus_{3 \; copies} k\homBr(3,1)
\]
so $k\homBr(3,1)$ is projective. 
Since $k\homBr(n,n-2)$ is a left $k \homBr(n,n)$-module by the category
composition, 
the natural $k$-linear extension of the 
injection
$i_{4,n}: \homBr(3,1) \hookrightarrow \homBr(n,n-2)  \qquad (n>2)$ 
allows us to induce to $ k \homBr(n,n) \; i_{4,n} (k\homBr(3,1) )$,
which is therefore also left projective. This is a submodule of 
$k \homBr(n,n-2)$ by construction; but considering for example 
the `herniated' form of a diagram in $\homBr(n,n-2)$ as in (a) below:
\[
(a) \; \includegraphics[width=1.2in]{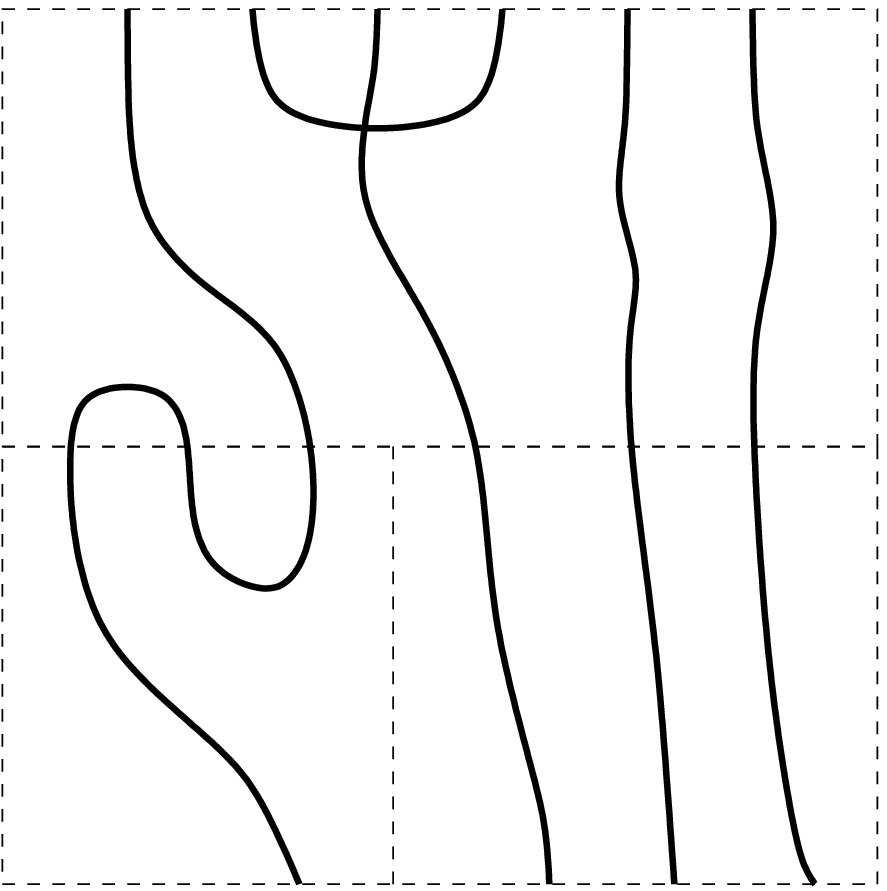} 
\qquad \qquad
(b) \; \includegraphics[width=1.2in]{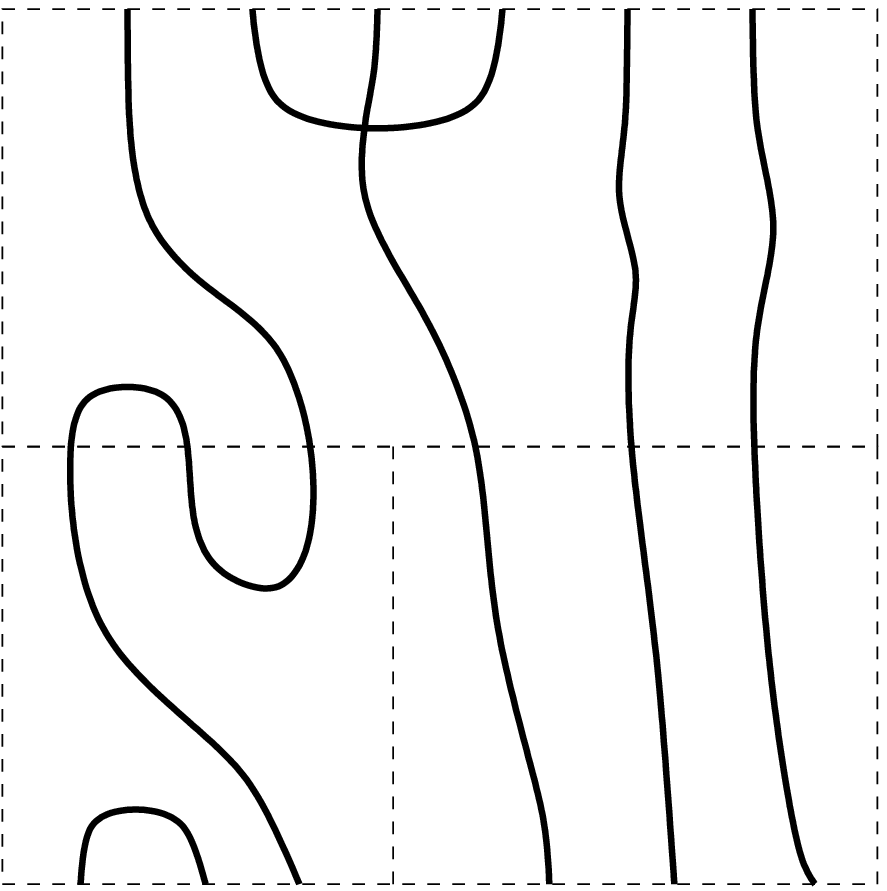} 
\]
we deduce that every diagram appears in the submodule and hence
\[
k \homBr(n,n-2) \;\; = \;\;\; k \homBr(n,n) \;\; i_{4,n} (\homBr(3,1) )
\]
is left projective.
The (left-handed version of the) claimed isomorphism is indicated in
the passage to figure~(b) above
(in particular this shows that a suitable choice for $e$ in case $n>2$ is 
$e=U_{23} U_{12}$). 
In case $\delta$ invertible in $k$ one sees directly that $k\homBr(2,0)$ is
left projective.  
\Qed
}


\mdef \label{de:snup}
The first section in (\ref{reg seq1}) obeys
\[
k \homBr^{\leq n}(n,n) / k \homBr^{\leq n-2}(n,n) \;\; \cong \; k S_n
\]
Thus each $S_n$-module induces an identical $B_n$-module, where the
action of any diagram with fewer than $n$ propagating lines is by 0.

\medskip

Via Proposition~\ref{pr:basis} and the various definitions:

\mpr{ \label{pr:GG1}
For $\lambda \vdash l$ and regarding $\specht(\lambda)$ as a
$B_l$-module as in (\ref{de:snup}), we have 
\[
\Delta_{2m+l}(\lambda) \cong \; G^m \specht(\lambda)
\]
}

\mdef
In particular (unless $n=2$ and $\delta=0$) the category 
$B_{n-2} \! -\!\mod$ fully embeds in 
$B_{n} \! -\!\mod$ under $G$,
and this embedding takes $ \Delta_{n-2}(\lambda)$ 
to  $ \Delta_n(\lambda)$. 

The embedding allows us to consider a formal limit module category
(we take $n$ odd and even together),
from which all $B_{n} \! -\!\mod$ may be obtained by localisation. 


%

By construction 
\mpr{\label{pr:simp1}  The set
$\{ \head ( \Delta_n(\lambda) ) \; | \; \lambda \vdash n,n-2,... \}$
is a complete set of simple modules
for $B_n(\delta)$ over any field $k$. 
}

\proof{To show that $\head\! ( \Delta_n(\lambda) )\;$ is simple, 
the only case not covered by applying
  Prop.~\ref{pr:Fexact1} to Prop.~\ref{pr:GG1} 
(or indeed by \CDMI)
is 
$\Delta_{2m}^{\QQ=0}(\emptyset)$ ($m>1$).
Here apply right exact functor $G^{m-1}$ to 
\[
0 \rightarrow \Delta_2(2) \stackrel{\sim}{\rightarrow}
\Delta_2(\emptyset) \rightarrow 0
\]
and use that $G^{m-1} \Delta_2(2)$ has simple head.
Completeness follows from Prop.~\ref{pr:projDfilt1}. \Qed
}

However regarded as a list this construction may give rise to multiple
entries, depending on $k$ and $\delta$. 
Over the complex field there is no overcount with $\delta \neq 0$,
and with $\delta=0$ just the element $\lambda=\emptyset$ should be
excluded (as shown by the case treated above).  
\\
This completes task~(1) over $\C$.

%


\newcommand{\ind}{{\mbox{Ind}}} 

\mpr{ \label{pr:ind1} \label{pr:ResG=Ind1}
{\rm \cite[Lemma~2.6,Prop.2.7]{\CDMi}}
Let $\ind$ and $\res$ denote the induction and restriction
functors associated to the injection 
$B_n(\QQ) \hookrightarrow B_{n+1}(\QQ)$.
\\
(i) We may identify the functors $\res \; G - \; = \ind -$. 
\\
(ii)
Over the complex field we have short exact sequence
\[
0 \rightarrow \bigoplus_{\mu \triangleleft \lambda} \Delta_{n+1}(\mu) 
  \rightarrow \ind \; \Delta_n(\lambda)
  \rightarrow \bigoplus_{\mu \triangleright \lambda} \Delta_{n+1}(\mu)
  \rightarrow 0
\]
(recall $\mu \triangleleft \lambda$ if $\mu$ is obtained from $\lambda$
by removing one box from the Young diagram).}

\proof{(i) Unpack the definitions.

(ii) Note from (i)  and Prop.~\ref{pr:GG1}
that it is enough to prove the equivalent result for
  restriction.
Use the diagram notation above. Consider the restriction acting
on the first $n$ strings. We may separate the diagrams out into those
for which the $n+1$-th string is propagating
(which span a submodule, since action on the first $n$ strings cannot
change this property), 
and those for which it
is not. The result follows by comparing with diagrams from the
indicated terms in the sequence, using the induction and restriction
rules for Specht modules.}

\subsection{Characters and $\Delta$-filtration factors}\label{ss:char}

\mdef \label{de:d=0caveat}
Over the complex field the modules 
$\{ \Delta_n (\lambda) \}_{\lambda\in\Lambda^n}$
are pairwise 
non-isomorphic except precisely in the case $n=2,\QQ=0$ in (\ref{de:base1}). 
If $\QQ\neq 0$ the heads are also distinct,
so there is a unique expression for any character in terms of
$\Delta$-characters. This means that the $\Delta$-filtration
multiplicities for projectives,
denoted $(P_i: \Db{n}{\lambda} )$, 
are also uniquely defined. 
The set $\{ \Pb{n}{\lambda} \}_{\lambda}$
of isomorphism classes of indecomposable projectives inherits
its labelling scheme from the simples in the usual way.

For the case $\QQ=0$, when $n=2$ the isomorphism means that these
multiplicities are not uniquely defined (we could simply discard one
of the isomorphic modules to make them so). 
For all other $n$, however,
provided we asign $\Delta_n(\lambda)$ as the top section of
$P_n(\lambda)$,
then the non-isomorphism of $\Delta$s removes this ambiguity.
In particular, the sectioning of projectives in the block of
$\emptyset$ up to $\lambda\vdash 4$ is indicated by
\[
P_4(2) = \Delta_4(2) /\!/ \Delta_4(\emptyset)
\qquad\qquad
P_4(31) = \Delta_4(31) /\!/ \Delta_4(2)
\]
(this is an easy direct calculation).
In this sense we may treat $\QQ=0$ as a degeneration of the more
general case, and treat the multiplicities $(P_i:\Delta_n(\lambda))$
as uniquely defined throughout. 
We do this hereafter. 


\mdef Recall from  
Proposition~\ref{pr:GG1}
\[
G \Delta_{n}(\lambda) = \Delta_{n+2}(\lambda)
\]
By Prop.~\ref{pr:simp1}
every $B_n$-module character can be expressed 
as a not necessarily non-negative combination of $\Delta$-characters:
\[
\chi(M) = \sum_{\lambda} \alpha_{\lambda}(M) \; \chi(\Delta(\lambda))
\qquad (\alpha_{\lambda}(M) \in \Z )
\] 
If in addition a module $M$ has a $\Delta$-filtration then this is a
non-negative combination and 
(with the caveat mentioned in (\ref{de:d=0caveat}))
\[
(G M : \Delta_{n+2}(\lambda)) = \left\{ \begin{array}{ll} 
   ( M : \Delta_{n}(\lambda)) & |\lambda| <n+2 \\
    0   & |\lambda| =n+2
\end{array} \right.
\]


The functor $G$ evidently takes projectives to projectives.
It also preserves indecomposability, so 
\[
G P_n(\lambda) = P_{n+2}(\lambda)
\]
Combining these we see that 
\[
(P(\lambda) : \Delta(\lambda))=1
\]
and otherwise
\eql(eq:unitri)
(P(\lambda) : \Delta(\mu))=0 \mbox{ if }  |\mu| \geq |\lambda|
\eq
Since these multiplicities depend on $n$ only through the range of possible
values of $\lambda$, 
for each $\QQ$ (here with $k=\C$)
there is a semiinfinite matrix $D$ with rows and columns indexed by
$\Lambda$ such that
\[
\StDM{P}{\lambda}{\mu} = \Dmat_{\lambda,\mu}
\]
for any $n$. 
In our case this `standard' decomposition matrix also determines the
Cartan decomposition matrix $\Cmat$ (see e.g. \cite[\S 1.9]{Benson95}). 
That is 
$\Dmat_{\lambda,\mu} = \StDM{P}{\lambda}{\mu} 
= \SDM{\Db{}{\mu}}{\Lb{}{\lambda}}$,
so that $\Cmat = \Dmat \Dmat^T$.
In particular there is an inverse limit of blocks that is a partition
of $\Lambda$. 

\label{para:C=DD}\label{eq:C=DD} 

Equation(\ref{eq:unitri}) says that  
the matrix $\Dmat$ is lower unitriangularisable. 
From this we have
\mpr{\label{pr:proj1}
If $P$ is a projective module containing $\Delta(\lambda)$ 
with multiplicity $m$ and no
$\Delta(\mu)$ with $|\mu| > |\lambda|$, then $P$ contains 
$P(\lambda)$ as a direct summand with multiplicity $m$. \Qed
}


The induction functor takes projective modules to projective modules,
and has a behaviour with regard to standard characters determined by 
Prop.~(\ref{pr:ind1}). 
From this we see that 

\mpr{ \label{pr:proj2} \label{pr:pweight}
For $e_i$ a removable box of $\lambda$,
\[
\ind \; \Pb{}{\lambda-e_i} \cong \Pb{}{\lambda} \bigoplus Q
\]
where $Q=\oplus_{\mu} \Pb{}{\mu}$ a possibly empty sum with no $\mu
\geq \lambda$. 
}

\noindent
Proof: By Prop.\ref{pr:proj1} 
a projective module is a sum of indecomposable projectives
including all those with labels maximal in the dominance order 
of its standard factors. Now use (\ref{pr:ind1}). 
\Qed

\mrem{
From the definitions we have
\[
F \Delta_{n}(\lambda) = \left\{ \begin{array}{ll} 
    \Delta_{n-2}(\lambda) & |\lambda| <n \\
    0   & |\lambda| =n
\end{array} \right.
\]
}

\mdef As we shall see shortly, the 
Young diagram labelling scheme we have
for the various indecomposable modules,
which is natural in light of (\ref{de:specht1}),
is the transpose of 
the labelling that it is convenient to work with in
describing the blocks. For this reason it is convenient to define
\[
\Dbt{n}{\lambda} = \Db{n}{\lambda^T}
\]
and similarly for simples and projectives.


\section{Blocks}
\newcommand{\MiBS}{minimal $\delta$-balanced skew}
\newcommand{\MIBS}{MiBS}
\newcommand{\MaBS}{maximal $\delta$-balanced subpartition}
\newcommand{\dBS}{$\delta$-balanced subpartition}
\newcommand{\lad}{\leftarrow^{\delta}}
\newcommand{\lu}{\lambda/\mu}

We now assemble the results we shall need on the blocks of the Brauer
algebras. 
These include important results from \CDMI, \CDMII, \CDMIII\  
and
extensions thereof.
The Young diagram inclusion partial order 
$(\Lambda,\subset)$ restricts to  
a partial order on each block
(any such construction evidently survives the inverse limit). 
By construction this order has a transitive reduction, 
that is, a directed graph that describes
the limit of Hasse diagrams. 
This graph is  key to our main result, and we describe it here.
For example we
endow the implicit definition of graph edges above (and in \CDMI)
with an explicit contruction that we shall need.

\subsection{$\delta$-balance}


Recall that the content $c(b)$ of a box $b$ in a Young diagram is
$c(b)=$ column position - row position.
In \CDMI\ we explain
how it is that the block
structure comes to depend on the relative content of the labelling
Young diagrams. 
It will be convenient now to cast the appropriate content condition
for blocks in various forms. 

\mdef The  $\delta$-{\em charge} of a box in a Young diagram is
\[
chg(b) = \delta-1 -2c(b)
\]
(cf. the conjugate function $ch(b)$ 
used in \CDMI). 

As for content, the lines of constant $\delta$-charge run parallel to the main
diagonal. 
The key difference from content is that the line of $\delta$-charge 0 
for given $\delta$ is no
longer (unless $\delta=1$) the main diagonal itself.
That is, the $\delta$-charge-0 main diagonal is shifted from the
ordinary main diagonal of the Young diagram. 
(Indeed for $\delta$ even there are no boxes with charge 0, so the
charge 0 line lies `between' diagonal runs of charge +1 and charge -1
boxes.) 

In the present setting, the point 
is that  
$\mu \subset \lambda$ is in the same block only if 
$\lambda^T/\mu^T$ consists of $\pm$charge pairs of boxes 
\CDMII.
(We give a precise statement shortly.)

For example, with $\delta=2$ the skew $(2^2)/(1^2)$ contains $\pm 1$,
so potentially (and in fact)
we have $L(2^2) \sim^{\delta=2} L(2)$.


\mdef \label{de:mibs1}
A Young diagram, or indeed any skew, 
can be considered as a planar graph all of whose faces are square. Its
geometrical dual graph is obtained by drawing a vertex for each 
square face 
and drawing an edge between a pair of vertices whenever the
corresponding pair of squares has a common edge.
A skew is called a {\em chain} if its dual graph is a chain.
A skew chain that is removable from a Young diagram is sometimes
called a { rim} of that diagram. 
Here
a {\em rim} is any skew that is a chain (i.e. not necessarily
removable from a given Young diagram). 

Two rims are {\em $\delta$-opposite} if there is a 
rotation by $\pi$ (hereafter called a
$\pi$-rotation) of the
plane about a point on the  
$\delta$-charge-0 main diagonal that takes one into the other.
\\
(Evidently this rotation is the same as reflection in the vertical
line defined by the point of rotation; followed by reflection in the
horizontal line defined by this point.)

\begin{figure*}
\[
\includegraphics[width=1.9in]{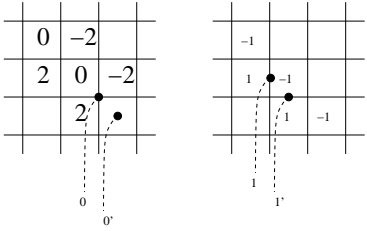}
\]
\caption{\label{fig:pipoints}Possible $\pi$-rotation points.}
\end{figure*}

Note that any such $\pi$-rotation is necessarily about a point positioned
as shown in one of the cases in Figure~\ref{fig:pipoints}.

\begin{figure*}
\[
\includegraphics[width=1.6in]{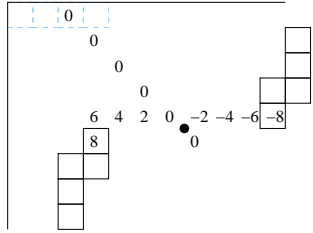}
\]
\caption{\label{fig:pipoints2} A $\pi$-rotation in case $\delta=5$.}
\end{figure*}

Note further that such a rotation has the effect of exchanging 
boxes in specific pairs, that are $\pm$charge pairs.
See Figure~\ref{fig:pipoints2} for an
example (rotation of rims about the black dot shown).
In this case the position of the charge-0 diagonal corresponds to
$\delta=5$. 

\mdef \label{de:MiBS}
A {\em minimal $\delta$-balanced skew} (\MIBS) is a skew that is a
 $\delta$-opposite pair of rims
such that no row of the skew is fixed by the associated  $\pi$-rotation.


\murem{Partition $\mu$ is a \MaBS\ of $\lambda$ (as in \CDMI)
if and only if $\lambda^T/\mu^T$ 
is a \MIBS.}
(Proposition~\ref{pr:cov1} below will serve to confirm this.)
The explicit geometrical form of the construction of \MIBS\ above 
(in contrast to the implicit construction
in the definition of \MaBS\ given in \CDMI)
will be crucial in what follows.

There are several
examples of \MiBS s  shown in Figure~\ref{fig:minskew eg}. 

\mdef Define a relation $(\Lambda, \lad )$   
by $\mu \lad \lambda$ if $\lambda/\mu$ is a 
 minimal $\delta$-balanced skew. 
Define $(\Lambda , <^{\delta})$ as the  partial order
that is the transitive closure of this relation. 


{\mlem{ \label{le:o12}
Possible $\pi$-rotation points for a MiBS are of the forms
 shown in Figure~\ref{fig:pipoints}.
In case-$0'$ there can be no intersection of the skew with the row or
column containing the point.
In case-1 there can be no intersection of the skew with the row
 containing the point.
Hence in either of these cases the skew is disconnected.
\Qed
}}

Define a partial order on the set of boxes 
occuring in  Young diagrams by 
$b' > b$ if $b'$  lies below {\em and} to the right of the top-left-hand
corner of $b$ (and $b' \neq b$).

\mlem{ \label{le:o2}
(Pinning Lemma)
Let $\pi_x$ be a rotation as above, and $b,b'$ two boxes comparable in
the above order,
then 
\[
b' > b \qquad \implies \qquad \pi_x(b) > \pi_x(b')
\]
\Qed
}


\mpr{\label{pr:cov1}
(I) If $\mu \subset \lambda$ and $\lambda / \mu$ a MiBS, then there is no
 $\mu \subset \mu' \subset \lambda$ such that $\mu'/\mu$ is a MiBS.
\\
(II) The relation $(\Lambda, \lad )$ is the cover 
(transitive reduction)
of the partial order $(\Lambda , <^{\delta})$.

}


\proof{(I): Let $\pi_0$ be the rotation fixing $\lambda/\mu$ and suppose
  (for a contradiction) that $\pi_{\gamma}$ fixes 
$\gamma=\mu'/\mu  \subset \lambda/\mu$.

The positive charge part of $\lu$ is connected, so there exists 
$b'\in\lu$ adjacent to $b \in \gamma$. 
Thus $\pi_0(b')$ lies in $\lu$ adjacent to $\pi_0(b)$. 
Since $\gamma$ is a skew over $\mu$, we have $b' \not\leq b$ and hence
(since adjacent) $b' > b$.
Thus $\pi_0(b) > \pi_0(b') $ by Lemma~\ref{le:o2}.

Suppose that $\pi_0 = \pi_{\gamma}$. Then $\pi_0(b') < \pi_{\gamma}(b)$,
contradicting that $\gamma$ is a skew over $\mu$. 
Thus  $\pi_0 \neq \pi_{\gamma}$.

Now, since    $\pi_0 \neq \pi_{\gamma}$, 
$\pi_0$ fixes no pair 
$b,\pi_{\gamma}(b)$ in $\gamma$. Thus for example no charge appears
more than once in $\gamma$, while all the charges appearing in
$\gamma$ appear twice in $\lu$. Thus in particular $\lu$ is
connected. 
Note that the rotation point of $\pi_0$ is necessarily half a box down
and to the right of $\pi_{\gamma}$. It then follows from
Lemma~\ref{le:o12} that $\gamma_+$ and $\gamma_-$ are disconnected
from each other. 

Let $c$ be the lowest charge box in $\gamma_+$. 
The box $\pi_0(\pi_{\gamma}(c))$ is below and to the right of it. Thus
there is a box of $\lu$ to its immediate right.
There cannot be a box of $\lu$ above it 
(since $\gamma$ is a skew over $\mu$)
so there is a box of $\lu$ to the right of  $\pi_0(\pi_{\gamma}(c))$.
But the $\pi_0$ image of {\em this} is to the left of $\pi_{\gamma}(c)
\in \gamma$, contradicting the $\gamma$ skew over $\mu$ property.

Claim (II) follows from (I) since
$\mu \subset \lambda$ is a necessary condition for 
$\mu <^{\QQ} \lambda$
so any failure of the MiBS relation to be a transitive reduction 
implies the existence of a $\mu'$ contradicting (I).
\Qed
}


{\mth{\label{th:bmap} {\rm \cite[Theorem~6.5]{\CoDMI}}
If  $\lambda/\mu$ is a minimal $\delta$-balanced skew then 
\[
\Hom(\Delta^{\delta}_n(\lambda^T) , \Delta^{\delta}_n(\mu^T) ) \neq 0
\]
\Qed
}}


Write $\Lambda^{\sim\delta}$ for the \RST\ closure of the partial order 
$(\Lambda , <^{\delta})$.
Write $[\lambda]_{\delta}$ for the
$\Lambda^{\sim\delta}$-class of $\lambda\in\Lambda$.

\mpr{{\rm \cite[Corollary~6.7]{\CoDMI}}
The relation 
$\Lambda^{\sim\delta}$ gives the (transposed) block relation
for $B_n(\QQ)$ over the complex field. 
\Qed
}

\mdef For any $n$, we write $\proj_{\lambda}-$ for the projection
functor on the category 
$B_n(\QQ)-\mod$ onto the block associated to the class 
$[\lambda]_{\QQ}$ 
(i.e. the block containing $\Deltab_n^{\QQ}(\lambda^T)$). 

\mdef
Let $G_{\delta}(\lambda)$ be the 
$\lambda$-connected component of $(\Lambda, \lad )$.
This may thus be thought of as a directed acyclic graph. 
We call this the {\em block graph}.


\subsection{The block graph}

The structure of the graphs $G_{\QQ}(\lambda)$ 
will be crucial for the statement and proof of the main Theorem.
We can describe it as follows.

\mdef 
Let $\Pow_{even}(\N) \subset \Pow(\N)$ denote the set of subsets of
$\N$ of even order. 
Define
a directed graph, $\Ge$, with
vertex set $\Pow_{even}(\N)$; 
and labelled edges:
\[
a \stackrel{\alpha}{\rightarrow} b \qquad \mbox{ if } \qquad
     a\setminus b=\{ \alpha \}, \;\;   b\setminus a=\{ \alpha+1 \}
  \qquad ( \alpha\in\N )
\]
\[
a \stackrel{12}{\rightarrow} b \qquad \mbox{ if } \qquad
     a\setminus b=\emptyset, \;\;   b\setminus a=\{ 1,2 \}
\]
See Figure~\ref{fig:valley graph}.
(There is a corresponding graph $\Go$ with vertices given by subsets of
$\N$ of odd order. 
The {\em toggle map} between the vertex sets given by toggling the
presence of 1 so as to make an odd set even is readily seen to pass to
a graph isomorphism (the edge labels 1 and 12 are interchanged).) 


We shall shortly construct an isomorphism 
$G_{\QQ}(\lambda) \cong \Ge$ for each $\QQ,\lambda$. 
For now we note that the case
$G_2(\emptyset)$ takes a relatively simple form.
The vertex map, $o_2 : [\emptyset]_2 \rightarrow \Pow_{even}(\N)$, 
is as follows. First
 draw the main diagonal on the Young diagram, as in 
these three examples from $[\emptyset]_2$:
\[
\includegraphics[width=2.3in]{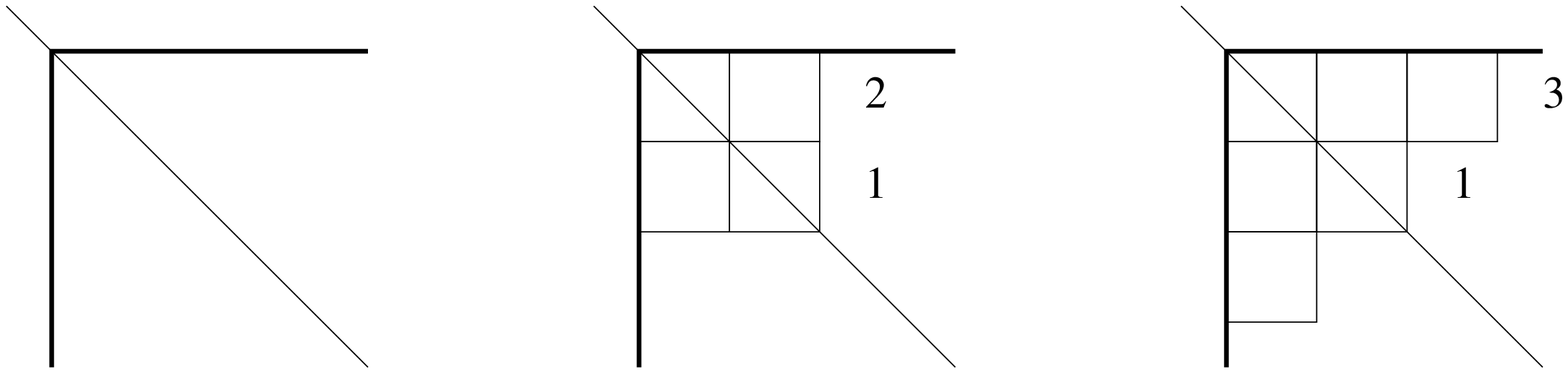}
\]
then count the number of boxes wholly or partly to the right of the
diagonal in each row, and write down the subset of these numbers that
are positive. 
Thus our examples become 
$\emptyset - \{ 2,1 \} - \{ 3,1 \} \cdots$. 
Comparing with (\ref{de:MiBS}) we readily see that 
$o_2$ passes to an isomorphism
$G_{2}(\emptyset) \cong \Ge$.


To generalise this
it is useful to give an alternative statement which emphasises
the geometrical nature of the block condition,
following \CDMII. 

Suppose $\lambda/\mu$ a \MiBS.
Note that if we suspend,  for intermediate steps, 
the dominance requirement (the requirement to
work with partitions rather than arbitrary compositions) 
then we can build $\lambda$ from $\mu$ 
by a sequence of transformations on pairs of
rows. Each transformation extends two rows: adding part of one row, and 
the corresponding opposite charges in the other row.
The no-row-fixed condition 
of (\ref{de:MiBS}) ensures that it is always pairs of rows
(as opposed to a single row) that are involved. 
For each row in question one takes the leading edge of the row in
$\mu$ and performs the two reflections mentioned in (\ref{de:mibs1}). 
The vertical reflection 
(i.e., in a horizontal line) simply swaps the two
rows. The other reflection takes this leading edge as far beyond the
charge-0 diagonal as it was short of it beforehand.
From these remarks it will be evident that this transformation can be
reformulated as in (\ref{de:shift-embed1}) et seq..

\mdef
\murem
Alternatively $\lambda$ can be built by a sequence of transformations
manipulating columns in pairs. The difference is firstly that, unless
we transpose, the intermediate stages are neither partitions nor
compositions (they are `transpose compositions'); and secondly that it
is possible in some cases to require a manipulation on a single
column, rather than a pair; and thirdly that the no-row-fixed
condition must still be imposed.
In light of this we use here the rows-in-pairs version.



\mdef Define a partial order $(\R^{\N},\geq)$ by $v\geq w$ if 
$v_i \geq w_i$ for all $i$.
(Write $v > w$ if  $v\geq w$ and $v \neq w$.)

\mdef \label{de:shift-embed1}
For $\delta\in\R$ 
define
\[
\rho_{\delta}  = - \frac{\delta}{2}(1,1,...) -(0,1,2,...) \in \R^{\N}
\]
For $\Z^f$ the subset of finitary elements of $\Z^{\N}$ define
\eqal(x)
e_{\delta} :  
\Z^f  & \hookrightarrow & \R^{\N}
\\
\lambda & \mapsto & \lambda   
   + \rho_{\delta}
\eqa
In other words,
since $\Lambda \hookrightarrow \Z^f$, 
we have, for each $\delta$, 
embedded our index set $\Lambda$ into a Euclidean space.
Thus our blocks $[\lambda]_{\QQ}$ 
now correspond to collections of points in this
space. 

Example:
\[
e_2 (\emptyset) =(0,0,0,0,...) -(1,1,1,1,...)-(0,1,2,3,...) \; = \; (-1,-2,-3,-4,...)
\]


\mdef \label{de:dom}
Note that all the image points $e_{\delta}(\Lambda)$ 
are strictly descending sequences.
We call such sequences {\em dominant}. 
Indeed  all the image points $e_{\delta}(\Lambda)$ 
are {\em strongly} descending sequences,
meaning that $v_i - v_{i+1} \geq 1 $ for all $i$.
We write $\Ap$ for the set of strongly decreasing sequences.

Considering for a moment the magnitudes of terms in a sequence in $\Ap$,
we see that each magnitude occurs at most twice, i.e. in a sequence of
form 
$
(...,x,...,-x,...).
$ 
We call such a $\pm x$ pairing a {\em doubleton}.
Define a map 
\[
\Reg: \Ap \rightarrow \Ap
\]
such that $\Reg(v)$ is obtained from $v$ by removing the doubletons.

For example
\[
\Reg(1,-1,-3,-4,-5,-6,...)  = (-3,-4,-5,-6,...)
\]
(note in this case that the input is $e_2((2,1))$ 
while the output is $e_{6}(\emptyset)$, that is, the $\Reg$ map 
can increase $\QQ$);
\[
\Reg(4,3,1,0,-1,-5,-6,...)  = (4,3,0,-5,-6,...)
\]


\mdef \label{de:sing-set}
For $\lambda \in \Lambda$
write 
$\sis_{\QQ}(\lambda)$ for the set of pairs of rows $\{ i,j \}$ such that
 $(\lambda+\rho_{\delta})_j = - (\lambda+\rho_{\delta})_i$
(i.e. $e_{\delta}(\lambda)_j = -e_{\delta}(\lambda)_i $). 
Write 
$s_{\QQ}(\lambda)$ for the {\em singularity} of 
$e_{\delta}(\lambda)$: 
\[
s_{\QQ}(\lambda)  =  | \sis_{\QQ}(\lambda) |
\]


\mdef We say a sequence $v \in \R^{\N}$ is {\em regular} if no 
two terms have the same magnitude.
Let $\R^{Reg}$ denote the set of regular sequences.
Define a map
\[
o : \R^{Reg} \cap \Ap \;  \rightarrow \Z^{\N}
\]
as follows.
In the $i$-th term, $|o(v)_i|$ is the position of $v_i$ in
the magnitude ordering of the set of numbers appearing in $v$.
The sign of $o(v)_i$ is the sign of $v_i$,
unless $v_i=0$ in which case the sign is chosen so as to make an even
number of positive terms.
\\
(Remark: this sign choice in case $v_i=0$ is simply for definiteness.
The definition of the 
function we eventually use (constructed next) will make it independent
of this convention.) 


\mdef 
If $v$ is a descending signed permutation of $(-1,-2,-3,...)$ 
then we define $v |_+ \in \Pow(\N)$ as follows.
First take the subset of terms of $v$ that are positive. 
Then, if this set is of odd order, 
toggle the presence of 1 in this set so as to make it even.

Define
\beqa 
\Psio_{\delta} : \Lambda & \rightarrow & P(\N)
\\ \nonumber
\lambda & \mapsto &   o (  \Reg ( \ed(\lambda) )  ) |_+
\eqa

\mdef
Examples: 
$\emptyset \; \mapsto \;\; e_2(\emptyset)=(-1,-2,-3,...) \;\; 
              \mapsto \; \emptyset$

$\qquad (3,3) \mapsto (2,1,-3,-4,...) \mapsto \{ 1,2 \}$

$(3,3,3,1) \mapsto (3,2,1,-2,-4,-5,...)  \mapsto \{ 1,2 \}$

$(4,3,3,1) \mapsto (4,2,1,-2,-4,-5,...)  \mapsto \{ 1 \}
              \stackrel{toggle}{\mapsto} \emptyset$


\mlem{ \label{lem:bich}
Fix $\QQ\in\Z$ and $\lambda\in\Lambda$. That is, fix a class
$[\lambda]_{\QQ} \subset \Lambda$.
Then the restriction 
$o_{\QQ}: [\lambda]_{\QQ} \rightarrow \Pow_{even}(\N)$ 
is a bijection.
}
\proof{The construction of the inverse map (call it
  $o_{\QQ}^{\lambda}$) is straightforward. \Qed}


{\mth{\label{th:graph isom0}
For each $\delta$, $\lambda$, 
the map $\Psio_{\delta}$ passes to an isomorphism
\[
 G_{\delta}(\lambda) \cong  \Ge
\]
(via $\Go$ and the toggle map in case $o_{\delta}(\lambda)$ of odd order).
}}
\medskip

Lemma~\ref{lem:bich} shows that  $o_{\delta}$ restricts to a
bijection on vertex sets. 
The next few paragraphs 
build up to a proof (in (\ref{th:gg4}))
of the graph isomorphism.

{\mpr{
Fix a block, i.e. a pair $(\QQ, [ \lambda ]_{\QQ} )$.
If $(v,w)$ is an edge in $\Ge$ with label $\alpha$ then 
the corresponding pair 
$(\mu,\lambda)=(o^{\lambda}_{\delta}(v), o^{\lambda}_{\delta}(w))$ 
gives $\lambda/\mu$ a \MiBS.
}}

\noindent
This is just a useful restatement of part of 
Theorem~\ref{th:graph isom0}.


\subsection{Geometrical aspects of the block graph}\label{ss:geom3}

\mdef
A Euclidean space together with a collection of hyperplanes defines a
reflection group --- the group generated by reflection in these
hyperplanes. 
Note that 
{
\newcommand{\lambdav}{v}
\[
(ij) : (\lambdav_1,\lambdav_2,...,\lambdav_i,...,\lambdav_j,...)
           \mapsto (\lambdav_1,\lambdav_2,...,\lambdav_j,...,\lambdav_i,...)
\]
\[
(ij)_- : (\lambdav_1,\lambdav_2,...,\lambdav_i,...,\lambdav_j,...)
           \mapsto (\lambdav_1,\lambdav_2,...,-\lambdav_j,...,-\lambdav_i,...)
\]
are 
reflection group actions on $ \R^{\N}$. 
Write $\Dgroup$ for the group generated by these (all $i < j$). 
Write $\Dgroup \lambdav$ for the orbit of a point 
$\lambdav \in \R^{\N}$ under the action of  $\Dgroup$.
Write $\Dgroupp$ for the subgroup $\langle (ij) \rangle_{ij}$. 
}

\mdef
Note that $\Dgroup$ does not preserve the image
$e_{\delta}(\Lambda)$, for any $\delta$.
Indeed the closure of the dominant region (in the sense of
(\ref{de:dom}))
is a fundamental region for the   $\Dgroupp$ action on  $ \R^{\N}$. 
This region is bounded by the reflection hyperplanes 
$\{ (i \; i\! +\! 1)\}_{i \in\N}$ (as is 
the region of ascending sequences).
Although the blocks are not precidely $\Dgroup$-orbits
(we will see that in a suitable sense)
\[
\mbox{orbit $\cap$ dominant = block}
\]


Comparing the definitions of \MiBS\ (\ref{de:MiBS}), $\ed$
and $(ij)_-$ 
we see  that 
{\mlem{ \label{le:prod-com-ref}
If $\lambda/\mu$ is a \MiBS\ then 
$\ed(\lambda)$ can be obtained
from $\ed(\mu)$ by a sequence of one or more transformations by
$(ij)_-$s,
extending rows in pairs of $\delta$-balanced part-rows.
Specifically
\[
\ed(\lambda) = \; \left( \prod_{ij} (ij)_-  \right) \; \ed(\mu)
\]
where the product is over pairs of rows in the skew, from the outer
pair to the inner pair.
\Qed }}{

Note also 
that no subset of this product, applied to $\ed(\mu)$, 
results in a dominant weight.

It follows that the $\Dgroup$ action on $\lambda$, via this construction, 
at least traverses the block $[\lambda]_{\QQ}$.  
In  \CDMII\ it is shown that it intersects no other block.

}

\newcommand{\GG}{{\bf G}}

\mdef  
For $v \in \R^{\N}$ define
\[
V(v) = \Dgroup v \cap \Ap
\]
The partial order $(\R^{\N},\leq)$ restricts to a  
partial order $(V(v),\leq)$.
The latter (unlike the former) has a unique transitive reduction.
This reduction thus defines a directed acyclic graph, denoted
$\GG(v)$. 


\mpr{ \label{pr:gg1} {\rm \cite[Prop.7.1]{\CDMiii}}
For $\lambda\in\Lambda$ the map $e_{\delta}$ 
restricts to a bijection 
$[\lambda]_{\delta} \rightarrow V(\lambda+\rho_{\delta})$;
and this bijection extends to a graph isomorphism 
$$
G_{\delta}(\lambda) \cong \GG(\lambda+\rho_{\delta}).
$$
}

\noindent {\em Proof:}
By \cite[Th.5.2]{\CDMii} 
we have that $e_{\delta}$ defines a bijection between
$[\lambda]_{\delta}$ and $V(\lambda+\rho_{\delta})$. 
Note that $\mu \subset \nu \in \Lambda$ if and only if 
$e_{\delta}(\mu) < e_{\delta}(\nu)$. Thus, 
restricting this to $[\lambda]_{\delta}$, the graphs are covers 
(transitive reductions) of
isomorphic partial orders. These covers thus agree on 
arbitrarily large finite sub-orders,
and hence agree.
\Qed

Note that $v$ is regular if and only if every sequence in $\Dgroup v$ is
regular. 

\mpr{ \label{pr:gg2} {\rm \cite[Prop.7.2]{\CDMiii}}
For $v \in \Ap$ the map $Reg$ restricts to a bijection 
$V(v) \rightarrow V(Reg(v))$; and this bijection extends to a graph
isomorphism 
$$
\GG(v) \cong \GG(Reg(v))
$$
}

\noindent {\em Proof:}
The set of doubletons is an invariant of the elements of $V(v)$,
and there is a unique way of adding these into an element of
$V(Reg(v))$ that keeps the sequence decreasing. 
Thus the restriction of $Reg$ here has an inverse, i.e. 
 the set map is a bijection. 
Now suppose $t,u \in \Ap$ and $a \in \R$ such that 
$$
s=(t_1,t_2,...,t_i,a,t_{i+1},...)
\qquad
s'=(u_1,u_2,...,u_j,a,u_{j+1},...)
$$
are in $\Ap$.
Then $t<u$ if and only if $s<s'$. The $Reg$ map can be built from
pairs of such moves, so $t<u$ if and only if $Reg(t) < Reg(u)$,
which establishes the graph isomorphism.
\Qed

\begin{figure}
\[
\includegraphics[width=2.3in]{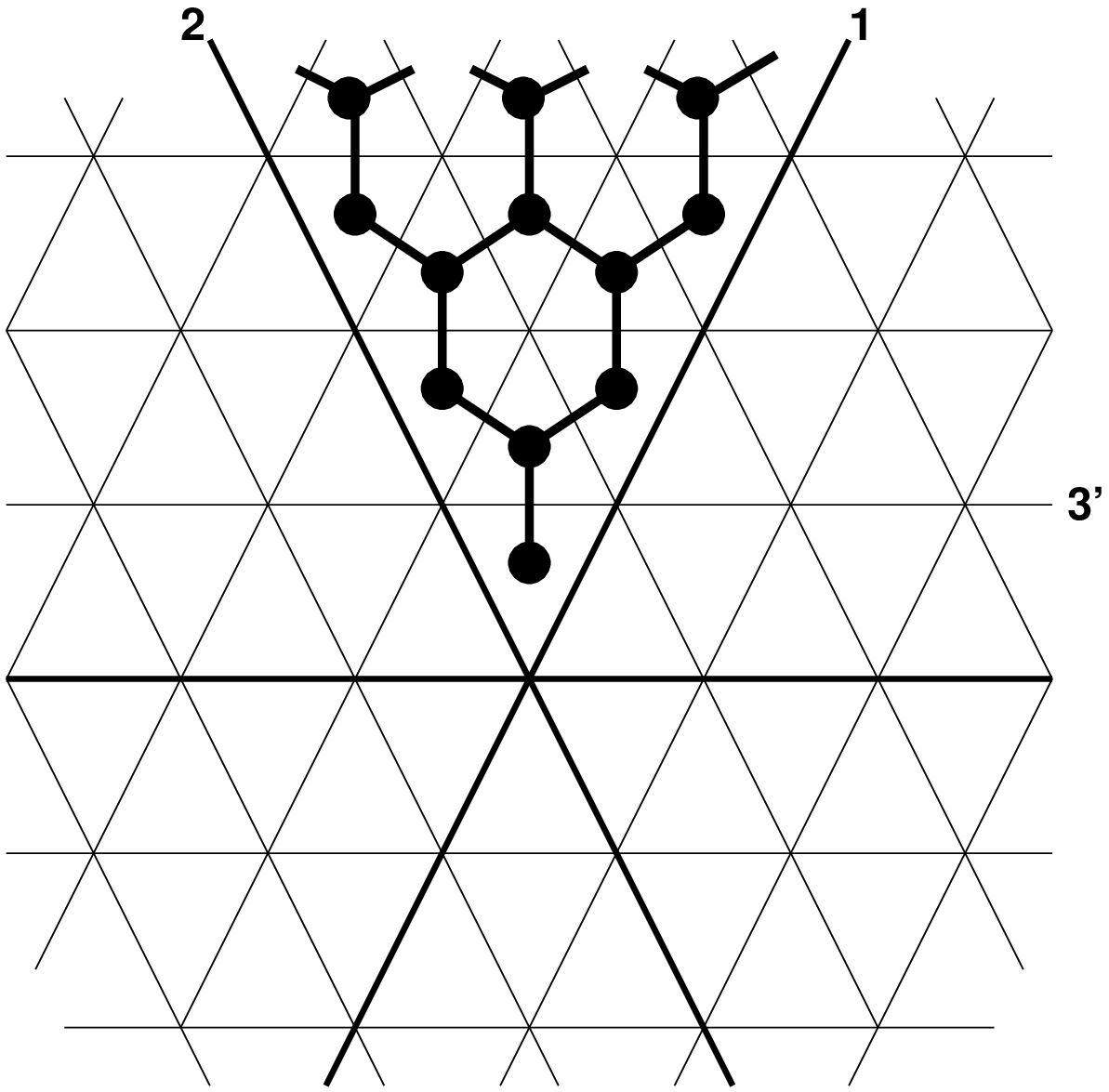}
\]
\caption{\label{A2dual} 
Example of a dominant dual graph: case  
affine-${A_2}/A_2$}
\end{figure}

\begin{figure}\caption{\label{fig:valley graph}\label{fig: odd2x}
The beginning of the  graph $\Ge$, with edge labels.
(Vertex labels have been written in an obvious shorthand.)
}
\[
\includegraphics*[viewport=-3 -3 623 900 , width=5in]{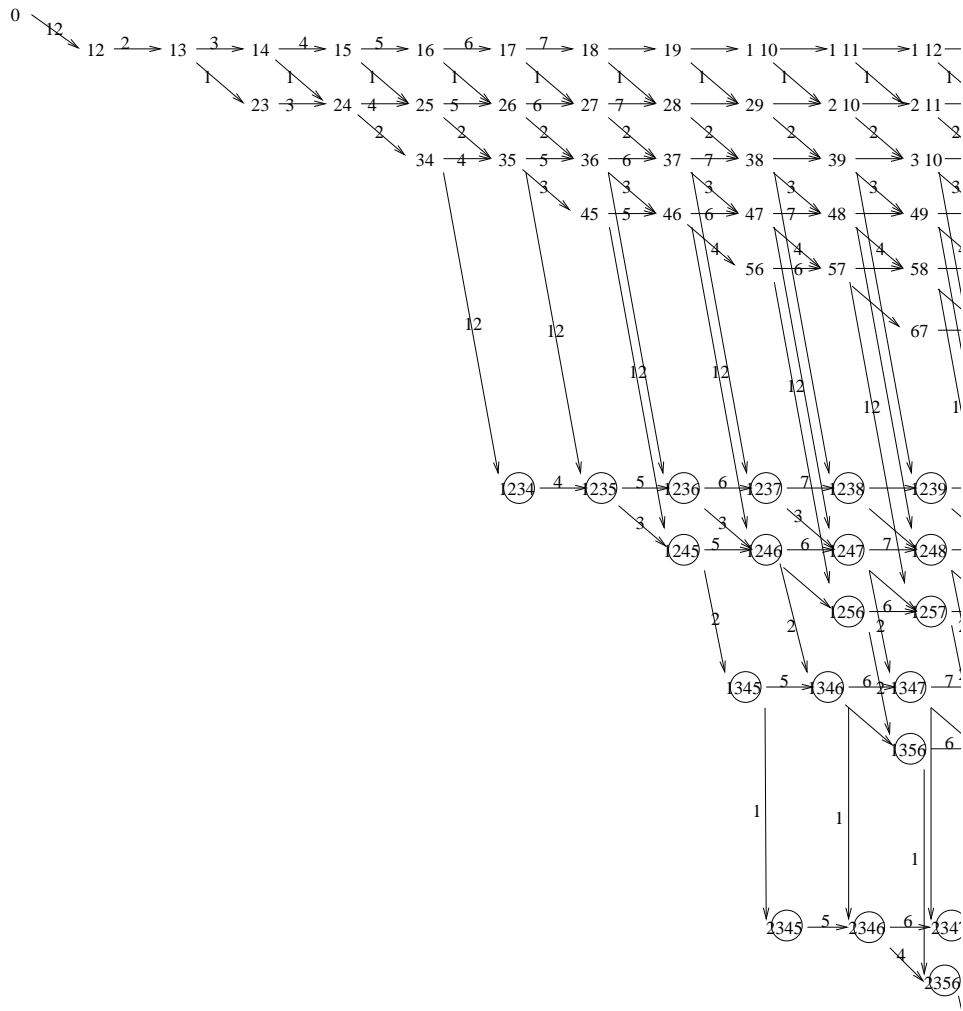}
\]
\end{figure}


\newcommand{\Dual}{D} 

\mdef To any Euclidean space $V$ and set of hyperplanes 
$\Hyp$ we may associate a
      {\em dual graph} $\Dual(\Hyp)$. 
This has a vertex for each connected component
      of the space with the hyperplanes removed (called an `alcove');
and an edge whenever the closures of two alcoves intersect in a
defining subset of a hyperplane (called a `wall').

If the set of hyperplanes is closed (under the reflections they
define) it may be generated by a minimal set 
defined by the walls bounding a single alcove
\cite{Humphreys90} (or see Section~\ref{ss:Hum}).
This minimal set of hyperplanes is thus in bijection with the edges
out of the dual graph vertex for the chosen `fundamental' alcove.
We have then two different enhancements of $\Dual(\Hyp)$ 
to include {\em edge labels}: {\em left edge labelling}
associates to each edge $(a,b)$ the hyperplane defined by $(a,b)$;
 {\em right edge labelling} 
requires the choice of a prefered alcove $C'$ and
associates to  $(a,b)$ the wall of $C'$ in the same reflection group
orbit as the wall $\overline{a}\cap\overline{b}$ defined by $(a,b)$.


Given a pair of a closed set of reflection hyperplanes and a closed
subset
$\Hyp_+$ 
(a parabolic), a 
{\em dominant dual graph} 
is the intersection of the dual graph with a
fundamental chamber 
(a connected component of the space with just the subset removed). 
For example Figure~\ref{A2dual}
shows the dominant dual graph for affine-$A_2$ 
(generated by the hyperplanes 1,2 and 3' shown)
over the subset 
corresponding to $A_2$ (generated by the hyperplanes 1 and 2).
If (as in the example) $\Hyp_+$ is maximal \cite{Humphreys90}
then only one alcove in each chamber has a subset of walls defining
$\Hyp_+$, 
and then by default one chooses the fundamental alcove to be the one
such in the fundamental chamber.

We write $\Gal$ for the dominant
dual graph of our
reflection group action $\Dgroup$ above (with parabolic $\Dgroupp$)
corresponding to the choice of 
$S_{\Dgroupp} = \{  (i \; i\! +\! 1) : i \in \N \}$ as reflection hyperplanes
bounding the fundamental chamber, and (to make contact with the 
given notion of dominance) such that descending sequences lie in the
fundamental chamber; 
and of 
$\{ (12)_- \} \cup S_{\Dgroupp}$ as reflection hyperplanes bounding
the fundamental alcove.
\\
(Figure~\ref{fig:valley graph} shows
a graph isomorphic to $\Gal$, using an isomorphism we shall explain
next.)  


{\mlem{  \label{pr:gg3} {\rm \cite[Cor.7.3]{\CDMiii}}
If  $v\in\R^{\N}$ is regular then it  lies within an alcove; and
$V(v)$ consists of a point within each dominant alcove.
Thus $\GG(v) \cong \Gal$.
\Qed
}} 

A convenient example of a regular $v$ is $e_2(\emptyset)$. In light of
the lemma we may use the orbit of  $e_2(\emptyset)$ to label dominant
alcoves. 
In particular  $e_2(\emptyset)$ itself lies in the fundamental alcove. 
\fotnote{
--- WE SHOULD NOW BE ABLE TO PROMOTE THE LATER THEOREM ---
BUT DO WE WANT TO?
}
By considering the effect of simple reflections in this case,
such as 
\[
(15)_- (4,3,-1,2,-5,...) = (5,3,-1,-2,-4,...)
\]
we see:


{\mlem{ \label{lem:Gal=Ge}
The map from 
$V(e_2(\emptyset))$
to  subsets of $\N$ of even order which
discards all negative entries 
coincides with the final step in 
$o_2:[\emptyset]_2 \rightarrow \Pow(\N)$ and 
extends to a graph isomorphism $\Gal \cong \Ge$.
}}

{\mrem{
The relationship between the $\Dgroup$ action between adjacent
vertices in $\Gal$ and the edge labels in $\Ge$ is not, 
perhaps, transparent  
in this isomorphism, and we shall not need it explicitly
for the computation of decomposition matrices. 
It is useful in the discussion of \pKLp s, however.
We shall resturn to describe it in the second part of the paper.
}}


\mth{ \label{th:gg4}
For all $\delta$, $\lambda$ 
the map $o_{\QQ}$ passes (via $o_{2}^{\emptyset}$) to an isomorphism
$$
G_{\delta}(\lambda) \cong \Gal
$$
}
\proof{By (\ref{pr:gg1}), (\ref{pr:gg2}),  
(\ref{pr:gg3})
 and (\ref{lem:Gal=Ge}).
\Qed}

This is a remarkable result, since the right hand side does not depend
on $\lambda$ or even $\delta$. 

\section{Decomposition data}\label{ss:decomp data}
In this section we prepare the structures needed in the statement of
the main result.
The idea comes from solving for parabolic
Kazhdan--Lusztig polynomials for the $\Dgroup/\Dgroup_+$ system 
(a highly non-trivial exercise). 
However
 the {\em proof} of the main result requires a more general approach, so
we do not emphasise the Kazhdan--Lusztig theory aspect at this stage.
(See later.)

\subsection{Hypercubical decomposition graphs}\label{ss:hyperDG}

\mdef 
Let $\bb : \Pow(\N) \rightarrow \{ 0,1 \}^{\N}$ denote the natural bijection.
For example:
\[
\bb: \{ 1, 3, 5 ,6 \} \mapsto 101011
\]
(if $a$ is finite we omit the open string of 0s on the right).
\\
Define $\bb_{\delta} : \Lambda \rightarrow \{ 0,1 \}^{\N}$ by 
$\bb_{\delta}(\lambda) = \bb(o_{\delta}(\lambda))$. 


\mdef
A generalisation of Brauer diagrams is to allow singleton vertices.
A
vertex pairing in such a  diagram {\em covers} a vertex 
if the
pair lie either side of it. 
A {\em \TLd} 
($TL$ as in Temperley--Lieb)
is here a diagram drawn  in the positive quadrant of the plane, 
consisting of a collection of vertices  drawn on the horizontal part
of the boundary
(countable by the natural numbering from left to right); 
together with a 
collection of {\em non-crossing} arcs drawn in the 
positive quadrant, each terminating  in two of the vertices,
such that no vertex terminates more than one arc,
and no arc covers a singleton vertex. 
An example is: $\;$
\raisebox{-.081015121in}{%
\includegraphics[width=.7681in]{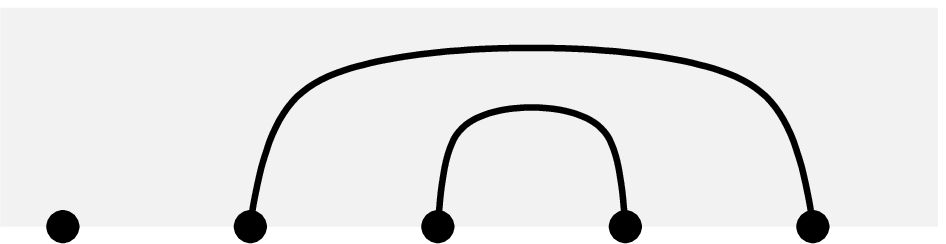}
}
\\
It will be convenient to 
label each arc by the 
associated pair of numbered vertices. 

\murem As with a Brauer diagram, it is the vertex pairings 
(and here singletons)
rather than the precise routes of the arcs that are important.


\mdef \label{de: TL algo}
Each binary sequence $b$ has a \TLd\ $d(b)$ constructed as follows.
\\
1. Draw a row of
vertices, one for each entry in $b$ (up to the last non-zero entry).  
\\
2. For each binary subsequence 01 draw an arc connecting the
corresponding vertices. 
\\
3. Consider the sequence obtained by ignoring the vertices paired in
2. 
For each subsequence 01 draw an arc connecting these vertices 
(it will be evident that this can be done without crossing). 
\\
4. Iterate this process until termination (it will be evident that it
terminates, since the sequence is getting shorter). 
\\
5. Note that this process terminates  either in the empty sequence
or in a sequence of 1s then 0s (either run possibly empty). 
Finally connect the run of 
vertices binary-labelled 1 in adjacent pairs (if any) from the
left. Leave the remaining vertices as singletons.


\begin{figure}
\includegraphics[width=5in]{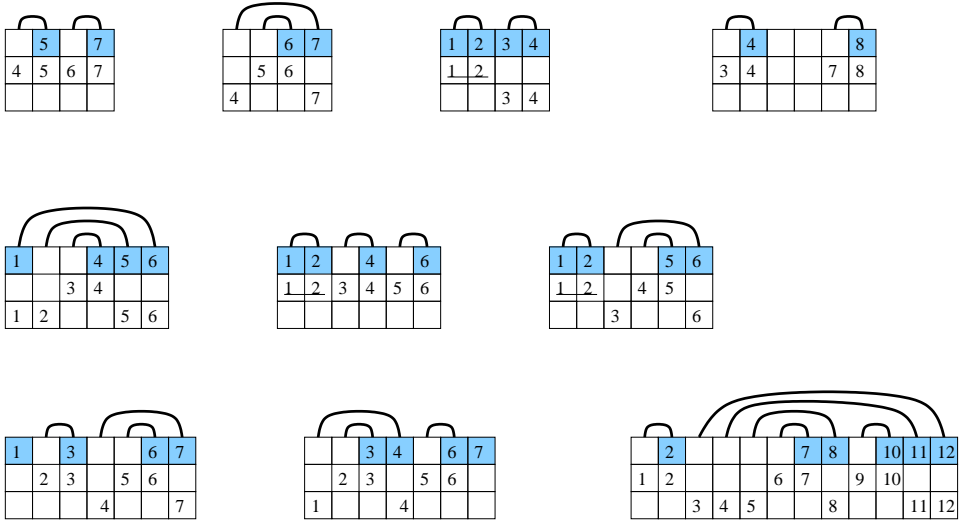}
\caption{\label{algo-check11} Examples for the map from 
sequences 
  to \TLd s, and to sets of pairs. 
In each case the
sequence for a set $a \in \Pow(\N)$ is indicated 
in the first
  (shaded) row of boxes. The second row shows the set of pairs of
  numbers $\Gamma_a$ extracted from the TL construction. 
The third row shows the further pairs added to obtain the set $\Gamma^a$.}
\end{figure}

Example: $d(10011) \; = \;\;\;
\includegraphics[width=.7681in]{xfig/egTL10011.eps} \;\;\;\;\;\;\;$
A number of examples are shown in Figure~\ref{algo-check11}. 

\mdef For $a \in P(\N)$ 
we write $\Gamma_{a}$ for the list of arcs (i.e. pairs) corresponding
to 01 subsequences,
and an initial 11 subsequence (i.e. if there is one in the 12-position); 
and  $\Gamma^{a}$ for the list of all arcs.

In particular, for example, 
\[
\Gamma_{1356} = \{ \{2,3 \} , \{ 4,5 \}  \}
\qquad\qquad\qquad
\Gamma^{1356} = \{ \{2,3 \} , \{ 4,5 \} , \{ 1,6 \} \}
\]
See Figure~\ref{algo-check11} for more examples. 
We may write $\Gamma_{\delta,\lambda}$ for
$\Gamma_{o_{\delta}(\lambda)}$,
and  
 $\Gamma_{\delta}^{\lambda}$ for  $\Gamma^{o_{\delta}(\lambda)}$.


\mdef \label{de:hyper DG}
A {\em hypercubical directed graph} is a rooted directed graph isomorphic to
the subset partial order  on some set $S$.
There is a notion of {\em parallel} edges 
(edges corresponding to deleting the same element of $S$). 
The edges coming out of the top vertex are called {\em shoulder}
edges, and every edge is parallel to one of these. 
\\
There is an obvious association with the notion of the
(geometrical) hypercube or hypercuboid, 
 i.e. the $\{0,1 \}$-span of any linearly
independent collection of vectors in a space. The notion of parallel
edges comes from this.


\mdef \label{de:hyp}
Each $a \in \Pow(\N)$ 
defines a hypercubical directed graph $h^a$,  as follows. 
The vertices are binary sequences
(these should be considered as 
identified with elements of $\Pow(\N)$ by the bijection, 
but it is convenient to treat them as binary sequences for the construction).
Firstly $a$ defines a binary sequence $\bb(a)$ and hence a \TLd\ 
$d( \bb(a) )$.
The top sequence in $h^a$ is the defining sequence $\bb(a)$.
There is an edge out of this corresponding to  
each completed arc in the \TLd\ $d( \bb(a) )$. The sequence at the
other end of a given edge is obtained from the original by replacing 
$01 \rightarrow 10$ 
(or $11 \rightarrow 00$)
at the ends of this arc. 
Indeed every parallel edge in the hypercube follows this transformation rule. 


\begin{figure}
\[
\includegraphics[width=2in]{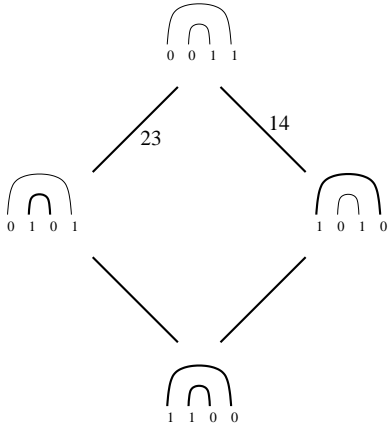}
\]
\caption{\label{fig:seg} Hypercube $h^{3 \, 4}$
(showing the TL arcs used in the construction).
}
\end{figure}
There is an example in Figure~\ref{fig:seg}
(and 
an example starting from given $\QQ$ and $\lambda$ 
in Section~\ref{ss:hypex}).

\mdef Note from the construction
that these hypercubes are multiplicity-free. That is, no two vertices 
have the same label.

Since fixing a block $[\lambda]_{\QQ}$ establishes a bijection between
$\Pow_{even}(\N)$ and  $[\lambda]_{\QQ}$ the construction for $h^a$ also
defines a hypercubical directed graph $h_{\QQ}(\mu)$ for each pair
$(\QQ,\mu) \in \Z \times \Lambda$,
obtained by applying $o_{\QQ}^{\mu}$ to the vertices. 
That is, abusing notation slightly,
\[
o_{\QQ}( h_{\QQ}(\mu)) = h^{o_{\QQ}(\mu)}
\]


\mdef
We label each edge of the hypercube (i.e. each direction) by 
$\{ \alpha, \alpha' \}$,
where $\alpha, \alpha'$  
are the positions of the ends of
the arc associated to this edge. 

If label 
$\alpha'=\alpha+1$ for an 01-arc, we may just label the edge by $\alpha$.
If $\{ \alpha, \alpha' \} = \{ 1,2 \}$
for a 11-arc we may just label the edge by $12$.
Note that these $\alpha$-edges and 12-edges 
in particular then coincide with edges  
of $\Ge$, although other edges do not.

\mdef \label{pa:isMiBS}
It follows from the construction 
and Theorem~\ref{th:graph isom0}
that if a vertex of 
some hypercube $h_{\delta}({\tau})$
is $\bb_{\delta}(\lambda)$ for some $\lambda$, then a vertex
beneath it down an $\alpha$ or 12-edge is  
$\bb_{\delta}(\mu)$ for some $\mu$ a 
\MaBS\
of $\lambda$. 


\mdef \label{para:decomp data}
Note that
we have assigned a hypercube to each 
appropriate binary
sequence and hence to each vertex of $\Ge$. Thus for any given block
$[\lambda ]_{\QQ}$ we have asigned a hypercube to each partition in
the block. The vertices in this hypercube then correspond to 
partitions in the same block (the defining one, together with one of
each of some collection below the defining one).
In this way we can use the hypercubes to determine,
for each $\QQ$, a matrix (of almost
all 0s, and some 1s), with rows and columns labelled by partitions.
The 1's in any given row are given  by the vertices of the hypercube 
associated to the partition labelling that row.
\\
In light of this interpretation we shall write $h_{\QQ}(\mu)_{\nu} =1$ 
if $\nu$ appears in $h_{\QQ}(\mu)$, and $=0$ otherwise.
\\
We will see in Theorem~\ref{th:decomposi}   
that the resultant matrix 
gives our block decomposition matrix.

\begin{figure}
\includegraphics[width=6.07in]{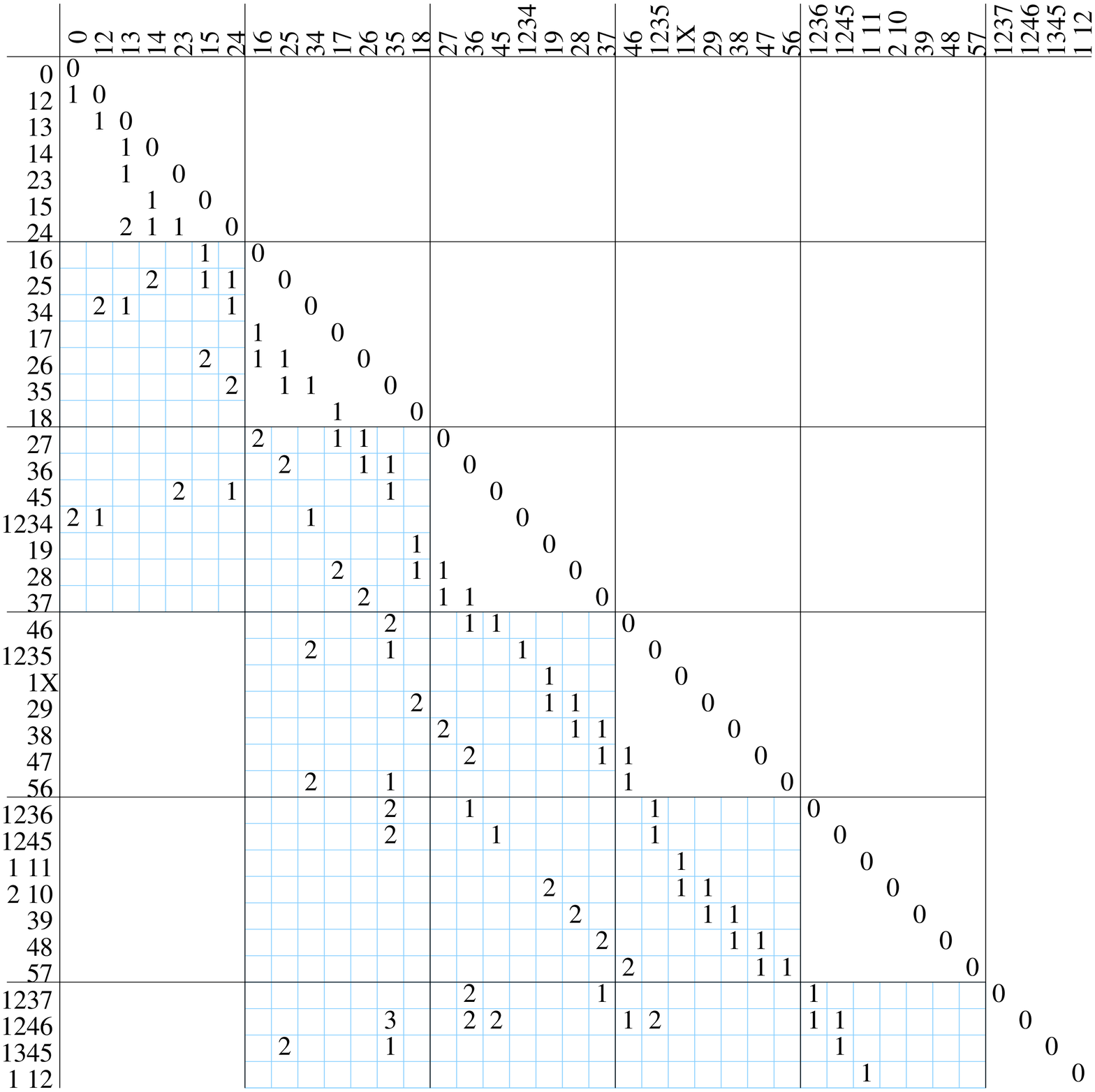}
\caption{\label{fig: big pKL} Table encoding array of 
polynomials in the $\Ge$ labelling scheme
(every non-zero polynomial is of form $v^i$, and the entry shown is $i$).}
\end{figure}

It will also be useful to consider an intermediate 
encoding, between the hypercube and the constant matrix row, in which we
record the {\em depth} $i$ of each entry in the hypercube, 
by writing $v^i$ ($v$ a formal parameter) instead of 1 in the
appropriate position.
(Thus this polynomial version evaluates to the decomposition matrix at
$v=1$.) 
The first few vertices of this form are shown in Figure~\ref{fig: big pKL},
using the $\Pow(\N)$ labelling scheme. 



\subsection{Hypercubical decomposition graphs: tools and examples}
\label{ss:hypex}

\begin{figure}
\[
(a) \!\!\!\!\!\!
\xymatrix{
 && 776532  \ar[dll]_{(35)_-} \ar[d]_{(26)_-} 
    \ar[drr]^{(16)_- (25)_- (34)_- (26)_- (35)_-}   \\
775522 \ar[d] \ar[drr] &&  766531 \ar[dll] \ar[drr] &&  652210 \ar[dll] \ar[d] \\
765521 \ar[drr]  && 642110 \ar[d] && 552200 \ar[dll]  \\
 && 542100
}
\;\;\;\;\;\; (b) \hspace{-.1in}\!\!
\xymatrix{
 && 1356 \ar[dll]^{23} \ar[d]^{45} \ar[drr]_{\underline{16}} \\
1256 \ar[d] \ar[drr] &&  1346 \ar[dll] \ar[drr] && 35 \ar[dll] \ar[d] \\
1246 \ar[drr]  && 25 \ar[d] && 34 \ar[dll]  \\
 && 24
}
\]
\[
(c)
\xymatrix{
 && (6,5,3,1,-2,-4,...) \ar[dll]^{(35)_-} \ar[d]^{(26)_-} 
                       \ar[drr] \\
(6,5,2,1,-3,...) \ar[d] \ar[drr] && (6,4,3,1,-2,-5,...) \ar[dll] \ar[drr] 
      && (5,3,-1,-2,...) \ar[dll] \ar[d] \\
(6,4,2,1,-3,...) \ar[drr]  && (5,2,-1,-3,...) \ar[d] 
 && (4,3,-1,-2,...) \ar[dll]  \\
 && (4,2,-1,-3,...)
}
\]
\caption{\label{fig:below:} Three labellings of the same hypercube
in case $\QQ=2$:
(a) partition labelling; 
(b)  $P(\N)$ labelling; 
(c) descending sequence labelling.
}
\end{figure}
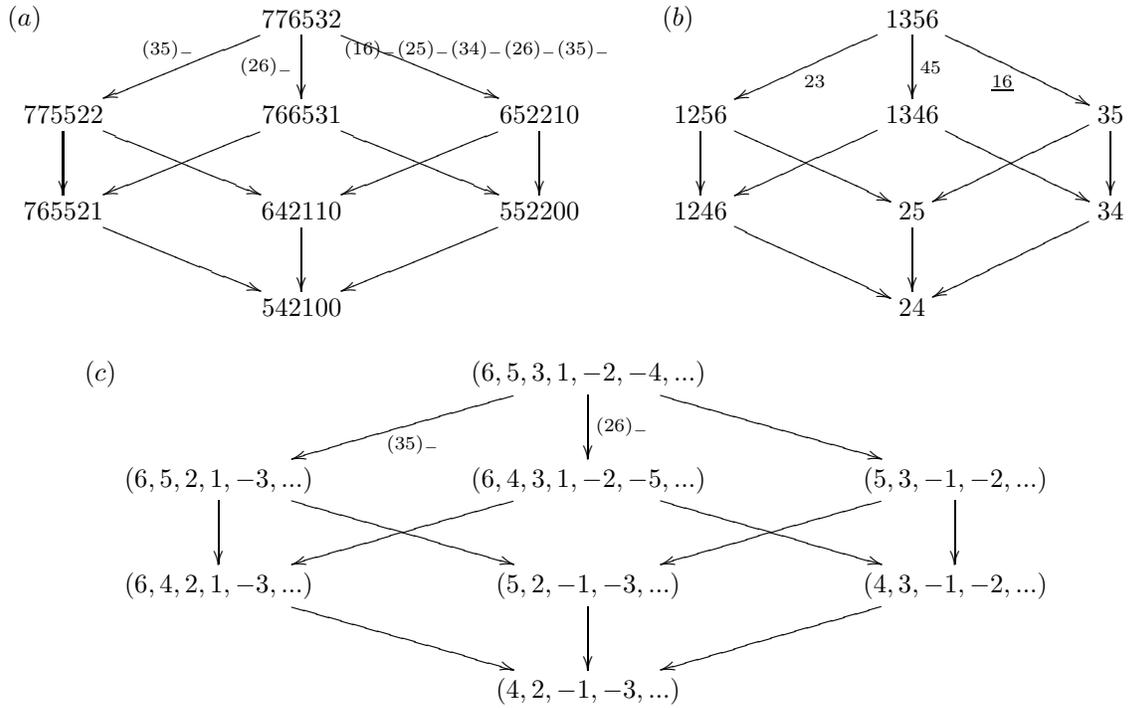

\begin{figure}
\[
\includegraphics[width=2.3in]{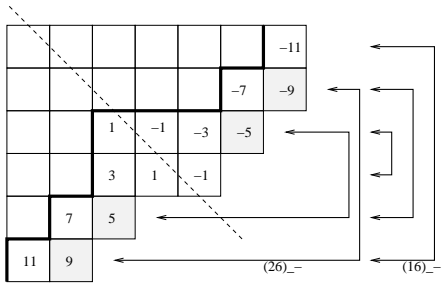}
\]
\caption{\label{fig:fol} 
Explicit reflections on $\lambda =776532$ in case $\QQ=2$.}
\end{figure}

\mdef Here is a concrete example of $h_{\QQ}(\lambda)$ with $\QQ=2$. We take 
$\lambda=(7,7,6,5,3,2)$ so 
\[
\lambda+\rho_2 = (6,5,3,1,-2,-4,-7,-8,...)
\]
giving $\Psio_2(\lambda)= \{ 1,3,5,6 \}$ and hence
$\Gamma_{\QQ}^{\lambda}=\{ \{ 2,3 \},\{ 4,5 \}, \{ 1,6 \} \}$. 
The specific hypercube 
(with  
integer partitions at the vertices) is thus (a) 
in Figure~\ref{fig:below:}.
In the figure
we have recorded both the $\alpha$-action and the specific
reflection group action required to achieve it on each edge
(for the shoulder layer).
The version in (b) shows the $\Ge$ vertex labels.
The version in (c) shows the $\rho_{\QQ}$-shifted vertex labels.
Figure~\ref{fig:fol} shows the explicit reflections and composite
reflection in the shoulder.
Note that the composite can be built as 
five dominance preserving but not all commuting reflections. 


\mdef
Keeping the same $\QQ, \lambda$, 
now consider $\lambda-e_i$ in case $i=4$. 

This gives 
$(7,7,6,4,3,2) \leadsto (6,5,3,0,-2,-4,-7,...) 
\stackrel{\Psio_2}{\leadsto} \{ 1,3,5,6 \}$ 
(by the toggle rule). 
This means that the hypercube $h_{\QQ}(\lambda-e_i)$  is isomorphic to that for $\lambda$
above, so in particular the $\alpha$-actions 
(the formal edge labels) are the same.
Note also that 
the specific reflections (realising these $\alpha$-actions)
in the shoulder of $h_{\QQ}(\lambda-e_i)$ 
are  the {same} as for $\lambda$. 

\murem 
We show in Section~\ref{ss:embed} that so long as $e_i$ 
does not `separate' a MiBS (in the sense of (\ref{de:separates}))
this
holds true in general. That is the hypercubes are isomorphic and the
reflections needed to move through the hypercube are the same.


\medskip

A more complicated example is given in Figure~\ref{fig: big h}.
We conclude this Section with some tools for manipulating these
hypercubes, that we shall need later.


\newcommand{\bump}[1]{{\hat #1} }
\newcommand{\bumpo}[1]{{\check #1} }

\mdef 
Let $b=(b_1,b_2,...)$ 
be a binary sequence, and $\alpha$ a natural number.
Then $\bump{\alpha} b$ is the sequence obtained from $b$ by inserting
01 in the $\alpha,\alpha+1$ positions 
(i.e. so that this pair become the elements in the $\alpha$ and $\alpha+1$
positions in the sequence,
with any 
terms at or above these positions in $b$ 
bumped two places further up in $\bump{\alpha} b$). 
\\
Similarly
$\bumpo{\alpha} b$ is the sequence obtained from $b$ by inserting
10 in the $\alpha,\alpha+1$ positions. 
\\
Examples: $\bump{2} 01 = 0011$, $\bumpo{2} 01 = 0101$.


\newcommand{\double}[1]{(1, #1 )} 

\mdef \label{de:ai act}
Let $h$ be a hypercube (i.e. the $\{0,1 \}$-span of any linearly
independent collection of vectors), and $\alpha$ a vector outside
the span of $h$ 
(or an operator that can otherwise be considered to shift all the
vertices of $h$ by the same amount). 
Then by 
$\alpha h$ we mean the translate of $h$ determined by $\alpha$, and by
$\double{\alpha} h$ we mean the new hypercube which
contains $h$ and a translate of $h$ by $\alpha$ together with the
edges in the $\alpha$ direction. 

More specifically, if $h$ is a hypercube whose vertices are binary
sequences, all of which have 01 (or all 11) in the $\alpha,\alpha+1$
positions, then  
 $\alpha h$ is the hypercube defined from
$h$ by modifying this $01 \rightarrow 10$ 
(respectively $11 \rightarrow 00$). 
In this case  $\double{\alpha} h$ is the hypercubical union of $h$ and
$\alpha h$. 

If the bumped sequence
$\bump{\alpha} \bb_{\QQ}(\lambda)$ makes sense, then by 
$\bump{\alpha} h_{\QQ}(\lambda)$ we understand the corresponding
vertex-modified hypercube
(insert 01 at the same position in every vertex binary sequence,
and modify any edge labels affected by this accordingly).
Note that this is not a hypercube of form $h_{\QQ}(\mu)$,
but a subgraph of somesuch. 
Similarly define 
$\bumpo{\alpha} h_{\QQ}(\lambda) $
(and note that $\bumpo{\alpha} h_{\QQ}(\lambda)
= \alpha \bump{\alpha} h_{\QQ}(\lambda)$).
Note that $\bumpo{\alpha} h_{\QQ}(\lambda)$ is another hypercube not 
 of form $h_{\QQ}(\mu)$. 
However
\eql(eq:hyp lemma)
\double{\alpha} \bump{\alpha} h_{\QQ}(\lambda) = h_{\QQ}(\mu)
\qquad
\mbox{ where } \bb_{\QQ}(\mu) = \bump{\alpha} \bb_{\QQ}(\lambda)
\eq
This is simply a restatement of part of the definition~(\ref{de:hyp}),
that will be useful later.


\section{Embedding properties of $\delta$-blocks in $\Lambda$}\label{ss:embed}

In this section we consider how the block graphs embed in $\R^{\N}$
and hence how the embeddings of the different block graphs relate to
each other. 
The result (\ref{pr:ind1}) means, loosely speaking, that the usual 
{\em metrical} structure on $\R^{\N}$ 
has relevance in representation theory. This, together with the
embedding results we develop here, will allow us to pass information
between blocks. 


\mdef
Suppose $w \in \Dgroup$ such that 
$w e_{\QQ}(\lambda) = e_{\QQ}(\mu)$.
When $\QQ$ is fixed we may write $w.\lambda$ for  $\mu$.
Also if $\lambda$ is a vertex of $\Ge$ or $G_{\QQ}(\mu)$ and $\alpha$
is the label on an edge out of $\lambda$ we write $\alpha\lambda$ for
the vertex at the other end.


\mdef \label{de:adjgraph}
The isomorphism implicit in Theorem~\ref{th:gg4}
between any pair of block graphs 
$G_{\QQ}(\lambda)$ and $G_{\QQ}(\lambda')$ defines a pairing of each
vertex in $G_{\QQ}(\lambda)$ with the corresponding vertex in 
$G_{\QQ}(\lambda')$.
A pair of block graphs is {\em adjacent} if they have the same
singularity, and every such pair of vertices is adjacent
as a pair of partitions.

\mdef \murem 
If $\lambda,\lambda'$ are adjacent partitions in the same 
$\Dgroup$-facet (in
the alcove geometric sense) then the corresponding pair of graphs are
adjacent, since the same reflection group elements serve to traverse
these graphs \CDMII, 
and reflection group elements preserve adjacency of 
partitions. We shall need to show adjacency of
a more general pairing of graphs. 

\mdef \label{para:fi}
For given $\lambda$, 
if $\lambda'=\lambda-e_i$ in (\ref{de:adjgraph}) above we write 
\[
f_i : [\lambda]_{\QQ} \; \rightarrow \;  [\lambda-e_i]_{\QQ}
\]
for the restriction of the graph isomorphism to vertices. 
(Strictly speaking $f_i$ depends on $\lambda$ too, but we suppress 
this for brevity.)


\mdef \label{de:separates}
Fix $\QQ$ and suppose $\lambda\in\Lambda$ has a removable box $e_i$.
Suppose that $\mib{\lambda}{\alpha}$ is a MiBS containing $e_i$.
Write $\pi_{\alpha}$ for the $\pi$-reflection fixing this MiBS. 
Then note that $\pi_{\alpha}(e_i)$ is an addable box of $\alpha\lambda$.
If $\mib{\lambda}{\alpha} \setminus \{ e_i, \pi_{\alpha}(e_i) \}$
is not a MiBS (of $\lambda-e_i$) we say that $e_i$ {\em separates} 
$\mib{\lambda}{\alpha}$.

\mdef Examples: crosses show boxes that separate; ticks show boxes
that do not:
\[
\includegraphics[width=2.83in]{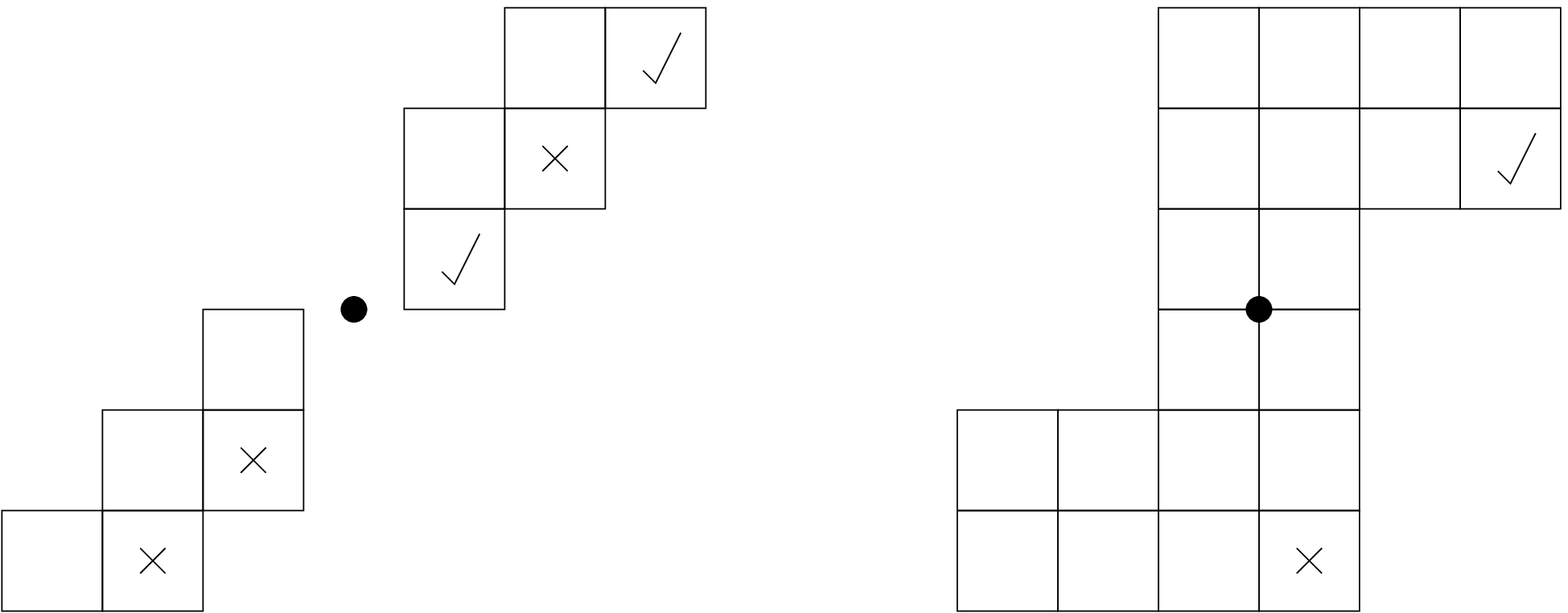}
\]


\mlem{\label{le:ch3}(Charge-row lemma)
Fix any $\QQ$. 
If a row $i$ of partition $\lambda$ ends in a box with charge $c$ we have
\[
(\lambda+\rho_{\QQ})_i = -\frac{c}{2}+\frac{1}{2}
\]
}

\mlem{ \label{le:combinME}
Fix $\QQ$ and suppose $\lambda\in\Lambda$ has a removable box $e_i$
such that singularity 
$\sing_{\QQ} (\lambda) =\sing_{\QQ} (\lambda -e_i)$. Then
\\
(I) $\Psio_{\QQ} (\lambda) = \Psio_{\QQ} (\lambda -e_i)$;
\\
(II) There does not exist a weight $\lambda-e_i-e_{i'}$ $\QQ$-balanced
with $\lambda$. 
\\
(III) There does not exist a weight $(\lambda-e_i)+e_i+e_{i'}$ $\QQ$-balanced
with $\lambda-e_i$. 
}
\proof{
Write $x$ for $(\lambda+\rho_{\QQ})_i$. 
That is 
\eql(eq:lemmy1)
\lambda+\rho_{\QQ} \sim (... ,w, \underbrace{x}_i  
 ,y, ... ),
\qquad
\lambda+\rho_{\QQ}-e_i \sim (...,w , \underbrace{x-1}_i 
, y,... )
\eq
with $w>x$ and $y<x-1$.
From this we see
that $x-1$ cannot occur in
$\lambda+\rho_{\QQ}$
(else it would occur twice in
$\lambda+\rho_{\QQ}-e_i$, contradicting the descending
property of the latter);
and similarly $x$ cannot appear in $\lambda+\rho_{\QQ}-e_i$.
\\
Note also that for $\lambda-e_i-e_{i'}$ to be $\QQ$-balanced
with $\lambda$ we would have to have (for $x \geq 1$)
\eql(eq:lemmy1.2)
\lambda+\rho_{\QQ}-e_i-e_{i'} \sim (... , \underbrace{x-1}_i  
,..., \underbrace{-x}_{i'} 
, ... )
\eq

We now split into two cases, depending on whether 
$\sis_{\QQ} (\lambda) = \sis_{\QQ} (\lambda -e_i)$. 
\\
(A) If $\sis_{\QQ} (\lambda) = \sis_{\QQ} (\lambda -e_i)$ :

(I)
The argument depends on the value of $x$. 
We split into  subcases (i-v).
\\
(i) If $x-1>0$: then $-(x-1)<0$ cannot appear in either sequence
(suppose it appears in the $j$-th position,
then $\{i,j \} \in \sis_{\QQ}(\lambda+\rho_{\QQ}-e_i)$
contradicting hypothesis (A));
\\
and similarly $-x$ cannot appear in either 
(else again  $\sis_{\QQ}$ changes between them).
\\
It follows that $x$ appears in $Reg(\lambda+\rho_{\QQ})$ 
and $x-1$ in the corresponding position in 
$Reg(\lambda+\rho_{\QQ}-e_i)$;
and that these sequences otherwise agree.

Suppose then that $x$ is, say, the $l$-th smallest magnitude entry in
$Reg(\lambda+\rho_{\QQ})$. 
If there is a smaller magnitude entry it's magnitude is smaller
than $x-1$, by the argument following Equations(\ref{eq:lemmy1})
and the argument above.
Since all these other entries are the same for the other sequence,
$x-1$ is  the $l$-th smallest magnitude entry in
$Reg(\lambda+\rho_{\QQ}-e_i)$.
Thus $\Psio_{\QQ}$ is unchanged.

\noindent
(ii) If $x=1$: 
then we have
$
\lambda+\rho_{\QQ} \sim (... ,w>1, \underbrace{x=1}_i   ,y<0, ... )
$.
We note that $-x=-1$ still cannot appear in either sequence
(else $\sis_{\QQ}$ changes). 
Thus 1 in $Reg(\lambda+\rho_{\QQ})$, 
respectively 0 in $Reg(\lambda+\rho_{\QQ}-e_i)$, is the smallest
magnitude entry. 
If there are an even number of other positive entries then this entry
does not contribute to $\Psio_{\QQ}$ in either case (in the former by
the toggle rule, and in the latter by the definition of the $o$-map).
If there are an odd number of other positive entries then this entry
contributes to $\Psio_{\QQ}$ in both cases (similarly).
Thus $\Psio_{\QQ}$ is unchanged.

\noindent
(iii) If $x=0$: 
then we have
$
\lambda+\rho_{\QQ} \sim (... ,w>0, \underbrace{x=0}_i   ,y<-1, ... )
$
and this time 
the hypothesis determines 
that $-(x-1)=1$ cannot appear in either sequence.
Thus 0 in $Reg(\lambda+\rho_{\QQ})$, 
respectively -1 in $Reg(\lambda+\rho_{\QQ}-e_i)$, is the smallest
magnitude entry.
 If there are an even number of strictly positive entries then this entry
does not contribute to $\Psio_{\QQ}$ in either case. 
If there are an odd number of  positive entries then this entry
contributes an element 1 to $\Psio_{\QQ}$ in former cases 
(by the definition of the $o$-map);
the entry -1 does not contribute in the latter case, but 
there is an element 1 by the toggle rule.
Thus $\Psio_{\QQ}$ is unchanged.

\noindent
(iv) If $x=1/2$: 
then we have
$
\lambda+\rho_{\QQ} \sim (... ,w>1, \underbrace{x=1/2}_i   ,y<-1, ... )
$.
Evidently there is no -1/2 in the former or 1/2 in the latter,
so the terms in the $i$-th position are the smallest magnitude terms
in their respective sequence, with all else equal.
Again by the toggle rule $\Psio_{\QQ}$ is unchanged.

\noindent
(v) If $x<0$: 
then 
neither $-x$ nor $-(x-1)$ can appear in either sequence
(else hypothesis (A) is violated much as before). 
The argument is then much as in (i).

\medskip

(II) 
For $x \geq 1$, by equation~(\ref{eq:lemmy1.2}) $\QQ$-balance here
would require
$-x+1$ in the $i'$-position in $\lambda+\rho_{\QQ}$,
and this is already disallowed under hypothesis (A).
\\
(The case $x=1/2$ does not arise; and the cases $x\leq 0$ are similar
to the above, with the order of $i,i'$ reversed.)

(III) By the rules of balance $e_{i'}$ cannot be in the same row as
$e_i$, so $(\lambda+\rho_{\QQ}+e_{i'})_i = (\lambda+\rho_{\QQ})_i=x$.
This would require that in the balance partner  
$(\lambda+\rho_{\QQ}-e_{i})_{i'} =-x$,
but this is already disallowed under hypothesis (A). 
\medskip
\\
(B) If $\sis_{\QQ} (\lambda) \neq \sis_{\QQ} (\lambda -e_i)$ :

(I)
Write $x$ for $(\lambda+\rho_{\QQ})_i$ as before. Then 
from Equation(\ref{eq:lemmy1})
we see firstly that $-x$  occurs in
$\lambda+\rho_{\QQ}$ and $1-x$ occurs in $\lambda+\rho_{\QQ}-e_i$
(if neither occurs then $\sis_{\QQ}$ does not change between them;
if only one occurs then $\sing_{\QQ}$ changes);
\\
of course it follows immediately that $1-x,-x$ occur (and are
adjacent) in both;
\\
secondly, by the same argument as above $x-1$ does not occur in
$\lambda+\rho_{\QQ}$.

In computing $\Psio_{\QQ}$ we discount the $\pm x$ pair in
$\lambda+\rho_{\QQ}$
and the $\pm(x-1)$ pair in  $\lambda+\rho_{\QQ}-e_i$.
The discrepancy is thus now a $1-x$ in  $\lambda+\rho_{\QQ}$ compared
to a $-x$ in $\lambda+\rho_{\QQ}-e_i$.
But if $1-x$ is the $l$-th largest magnitude entry in $\lambda+\rho_{\QQ}$
then $-x$ is  the $l$-th largest magnitude entry in
$\lambda+\rho_{\QQ}-e_i$,
with all else equal, so $\Psio_{\QQ}$ is unchanged.

(II) 
By equation~(\ref{eq:lemmy1.2}) $\QQ$-balance would require
$-x+1$ in the $i'$-position in $\lambda+\rho_{\QQ}$
as before. Although this is not disallowed here, it forces the $-x$
to lie in the next (that is, the $i'+1$) position. 
This would force a second $-x$ in the
same position in  $\lambda+\rho_{\QQ}-e_i-e_{i'}$,
which would thus not be descending --- a contradiction.

(III) Since $(\lambda+\rho_{\QQ}-e_i)_i=x-1$ we would require 
 $(\lambda+\rho_{\QQ}+e_{i'})_{i'} = 1-x$ for balance. Thus 
 $(\lambda+\rho_{\QQ})_{i'} = -x$. But we have already seen that 
$\lambda+\rho_{\QQ}$ contains both $1-x,-x$, so this would require 
$\lambda+\rho_{\QQ}+e_{i'}$ containing $1-x$ in two positions  --- a contradiction.

\Qed
}


\mlem{\label{le:do1}
Fix $\QQ$ and 
suppose $\sing_{\QQ} (\lambda) =\sing_{\QQ} (\lambda -e_i)$ as before. 
Suppose $\lambda$ has an edge down labelled $\alpha$,
i.e. $\mib{\lambda}{\alpha}$ is a \MIBS;
and let  $w$ be the product of commuting reflections such that 
$w e_{\QQ}(\lambda) = e_{\QQ}(\alpha\lambda)$,
as in Lemma~(\ref{le:prod-com-ref}).  
Then \\
(I) $w e_{\QQ}(\lambda-e_i)$ is dominant; \\
(II) $w e_{\QQ}(\lambda-e_i) = e_{\QQ}(\alpha(\lambda-e_i))$; \\
(III) $\alpha(\lambda-e_i) \adjac
         \alpha\lambda$ (i.e. they are adjacent).
}

\proof{
(I) We split into two cases:

\noindent
If $e_i$ does not intersect $\mib{\lambda}{\alpha}$ then 
$w e_{\QQ}(\lambda-e_i) $ is the same as 
$w e_{\QQ}(\lambda)$ everywhere except in row $i$:
$w e_{\QQ}(\lambda-e_i) = w e_{\QQ}(\lambda) -e_i$. 
Since $\lambda-e_i$ is dominant, $\lambda_i > \lambda_{i+1}$,
but $(\alpha\lambda)_i = \lambda_i$ in this case, and  
$(\alpha\lambda)_{i+1} \leq \lambda_{i+1}$, so 
$(\alpha\lambda)_i > (\alpha\lambda)_{i+1}$, so 
$\alpha\lambda -e_i$ is dominant, so 
$e_{\QQ}(\alpha\lambda -e_i) = w e_{\QQ}(\lambda-e_i) $ is dominant. 
\\
 If $e_i$  intersects $\mib{\lambda}{\alpha}$ then 
$\pi_{\alpha}(e_i)$ is addable to $\alpha\lambda$ as noted in
 (\ref{de:separates}). That is 
$e_{\QQ}( \alpha\lambda + \pi_{\alpha}(e_i)) = w e_{\QQ}(\lambda-e_i)$
is dominant. 

(II)
Firstly note that 
$\Psio_{\QQ}(\lambda-e_i)=\Psio_{\QQ}(\lambda)$ by
Lemma~\ref{le:combinME}, so $\alpha(\lambda-e_i)$ makes sense.
Similarly we have 
$\Psio_{\QQ}(\alpha(\lambda-e_i))=\Psio_{\QQ}(\alpha\lambda)$
(since both are equal to the formal set $\alpha\Psio_{\QQ}(\lambda)$).

Since $ w e_{\QQ}(\lambda-e_i)$ is dominant (by (I)) in the $\Dgroup$-orbit of
$\lambda-e_i$ there is some $\mu \in [\lambda-e_i]_{\QQ}$ such that 
$ w e_{\QQ}(\lambda-e_i) = e_{\QQ}(\mu)$.
Since it is adjacent to
$e_{\QQ}(\alpha\lambda)$
and has the same singularity, then by Lemma~(\ref{le:combinME})
(applied appropriately)
$\Psio_{\QQ}(\mu) = \Psio_{\QQ}(\alpha\lambda)$.
That is, $\mu = \alpha(\lambda-e_i)$. 

(III) Follows immediately from (II). 

\Qed
}


\mlem{\label{le:do2}
Fix $\QQ$. Suppose 
$\sing_{\QQ}(\lambda)=\sing_{\QQ}(\lambda-e_i)$ as before,
and $\alpha\lambda/\lambda$ is \MIBS\ (i.e. $\alpha$ is an edge up from
$\lambda$). 
Then there is a reflection group element $w$ such that 
$w. \lambda = \alpha\lambda$ (so $w. \alpha\lambda = \lambda$)
and $w .(\lambda -e_i)$ is dominant; 
whereupon $w .(\lambda -e_i) = \alpha(\lambda-e_i)$. 
}

\proof{
Suppose $w. (\lambda -e_i)$ is dominant. 
Then it is $\mu \in [\lambda-e_i]_{\QQ}$ adjacent to 
$w.\lambda=\alpha\lambda$
with the same singularity, hence the same $\Psio_{\QQ}$ by 
Lemma~(\ref{le:combinME}).
Thus it is enough to show that  $w. (\lambda -e_i)$ is dominant. 

Given that $w.\lambda$ is dominant, any {\em failure} of dominance 
of  $w. (\lambda -e_i)$ must
involve the $i$-th row itself
being shorter than row-$i+1$ in $w.(\lambda-e_i)$ 
(i.e. row-$i+1$ intersects the MiBS);
or a row with which row-$i$ is paired in $w$
($j$, say) being longer than row-$j-1$ in $w.(\lambda-e_i)$. 
We must consider the cases: 
(A) $e_i$ lies `behind' the skew (i.e. it's image under the
$\pi$-rotation $\pi_{\alpha}$ that fixes $\alpha\lambda/\lambda$ 
extends some row of the skew);
or (B) not.

(A)
In this case the failure would have to be
that the image of $e_i$ under the $\pi$-rotation broke dominance,
i.e. extended beyond the row above it.

Suppose  $e_i$ is behind other than the last row of the skew.
Then there is a box of the skew immediately to its right and one
immediately below it. The $\pi$-rotation images of these are behind
and above the image of $e_i$, so $w. (\lambda -e_i)$ is dominant.

On the other hand, suppose $e_i$ is behind the last row of the skew.
For example:
\[
\includegraphics[width=1.3in]{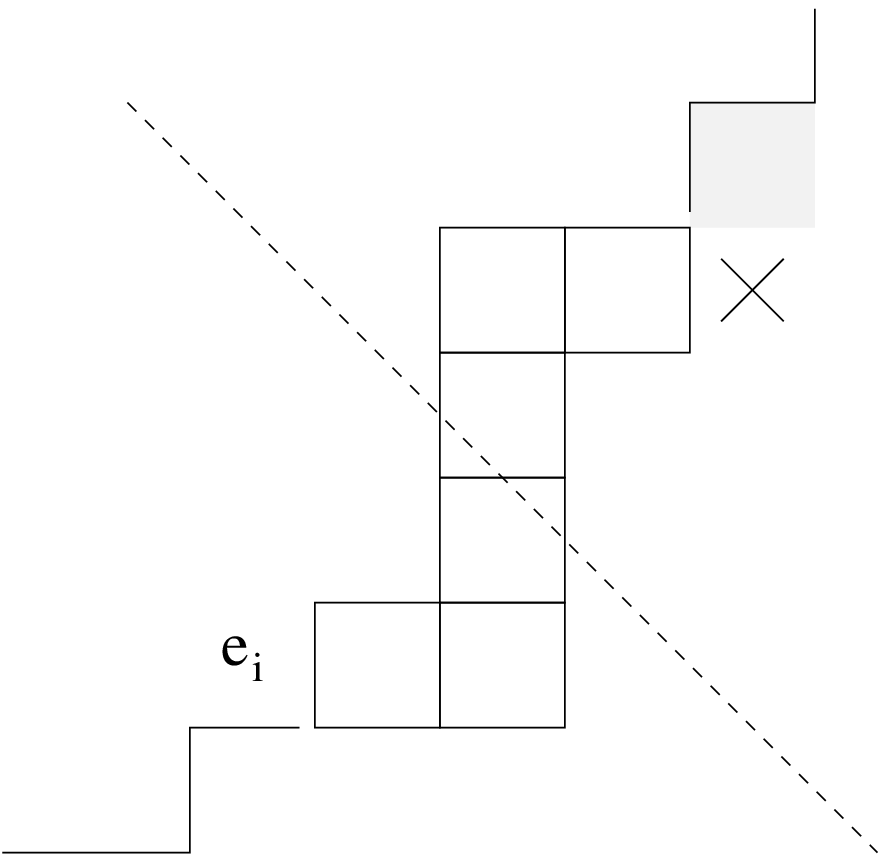}
\]
(the box $\pi_{\alpha}(e_i)$  is marked $\times$). 
Here $w.(\lambda-e_i)$ is dominant unless the box above 
$\pi_{\alpha}(e_i)$
is missing from $\lambda$. 
But if this is missing then this row and the $i$-row are a
singular pair in $\lambda-e_i$. Neither row can be in a singular pair
in $\lambda$ so this contradicts the hypothesis. 

(B)
If the $i$-th row is not moved by $w$ then the failure would have to be
that the skew $\alpha\lambda/\lambda$ includes a box directly under
$e_i$. But in that case a $\QQ$-balanced box to $e_i$ 
given by $\pi_{\alpha}(e_i)$  is directly to
the left of the skew, and we have a setup something like the
following:
\[
\includegraphics[width=2.3in]{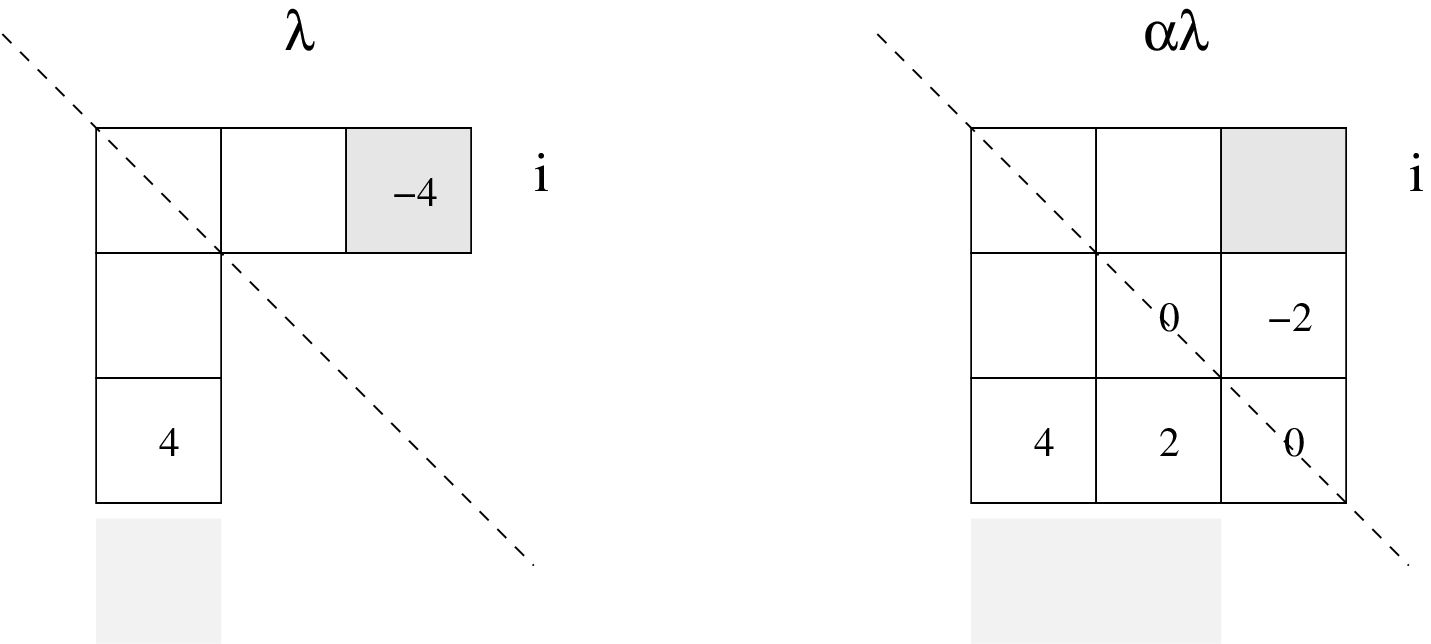}
\]
(the  $\QQ$-balanced box is the box marked 4).
If there is no box below the $\pi_{\alpha}(e_i)$  in $\lambda$ then  row-$i$
is not in a singular pair in $\lambda$, and row-$i$ and the row containing
the $\pi_{\alpha}(e_i)$  are a singular pair in $\lambda-e_i$,
thus $\sing_{\QQ}(\lambda) \neq \sing_{\QQ}(\lambda-e_i)$ 
so we can exclude this.
If there is a box below the 
$\pi_{\alpha}(e_i)$  in $\lambda$ then this row and row-$i$
are a singular pair in $\lambda$, and row-$i$ and the row containing
the $\pi_{\alpha}(e_i)$  are a singular pair in $\lambda-e_i$.
In this case, a $w$ which also 
has a factor
acting on the $i$-th and undrawn row 
has the same effect on $\lambda$ as one which does not. 
Its effect on $\lambda-e_i$ is to
restore the box $e_i$ and to add a box in the undrawn row.
This $w.(\lambda-e_i)$ is dominant since the added box is under a box
added in the original skew. 

\Qed
}


Since the block graph is connected we may use Lemmas~\ref{le:do1}
and~\ref{le:do2} to show:
{\mth{ \label{th:embed}
(Embedding Theorem) $\;$ 
If $\sing_{\QQ}(\lambda) = \sing_{\QQ}(\lambda-e_i)$
then $G_{\QQ}(\lambda)$ is adjacent to  $G_{\QQ}(\lambda-e_i)$.
\Qed
}}


\mdef \label{de:isomg}
Lemma~(\ref{le:combinME})(I)    
says that if the partitions $\lambda,\lambda-e_i$ have
the same singularity then they pass to the same point on the block
graph $\Ge$. 
That is 
\[
f_i (\lambda) = \lambda-e_i
\]
and so on. Thus for $\mu\in [\lambda]_{\QQ}$
\[
h_{\delta}(\lambda)_{\mu}  = h_{\delta}(\lambda-e_i)_{f_i(\mu)}
\]


\footnot{{CHAT:

\mdef
For any fixed $\QQ$, 
the graph isomorphism Theorem provides a bijection between any pair of
blocks. In general a pair of blocks related by having adjacent points
($\lambda$ and $\lambda-e_i$ say) will not be pairwise adjacent under
the isomorphism. We are interested in the conditions under which this
does happen. 
Since the reflection group action preserves adjacency of a pair of points, 
a sufficient condition would be that $\lambda$ maps to $\lambda-e_i$
under the isomorphism, and every edge in the block graph can be
realised by the same reflection group element in the two blocks. 
This is not true in general, so the question here is when does it
hold?

Let us suppose that $\sing_{\QQ} (\lambda) =\sing_{\QQ} (\lambda -e_i)$. 
(Since singularity is a block invariant this says the blocks have the
same singularity.) 
We CLAIM that this implies that the pair are isomorphic under the
graph isomorphism. 
Now let us consider some edge out of $\lambda$,
labelled $\alpha$ say. 
By the isomorphism there is an edge with the same label out of
$\lambda-e_i$. There are then various possibilities for the
relationship between $e_i$ and the MiBS $\mib{\lambda}{\alpha}$.

\mlem{
Fix $\QQ$ and suppose $\lambda\in\Lambda$ has a removable box $e_i$.

Then $e_i$ may or may not intersect a MiBS of $\lambda$. 

Cases:

There may be no $\mu$ below $\lambda$ in the block not containing
$e_i$.

There may be a $\mu$ below $\lambda$ in the block not containing
$e_i$, but no MiBS containing $e_i$.

There may be a $\mu$ below $\lambda$ in the block not containing
$e_i$, and a \oonn\ MiBS containing $e_i$.

There may be a $\mu$ below $\lambda$ in the block not containing
$e_i$, and a non-\oonn\ MiBS containing $e_i$.

The singularity 
$\sing_{\QQ} (\lambda) =\sing_{\QQ} (\lambda -e_i)$. 
\\
(I) 
\\
(II) 
\\
(III) 
}

\mdef Next we would like to show that the pairs of vertices in the
isomorphism $h_{\QQ}(\lambda) \cong h_{\QQ}(\lambda-e_i)$
(same singularity) are adjacent throughout. They all have the same
singularity, since they are in two blocks of the same singularity.
It would be enough if we could show that each pair can be reached by
applying the same set of reflections to the top pair (since
reflections preserve adjacency). 
Consider a reflection that takes some vertex $\mu$ down an edge
$\alpha$ of the
hypercube.
This same reflection takes the partner in the other hypercube to a
possibly non-dominant weight adjacent to $\alpha\mu$. 

1. can this really be non-dom?
The structure of this reflection is a nested sequence of row pair
reflections.
It preserves the dominance of $\mu$, and the partner differs from
$\mu$ in only one box in one row.
(I think that WLOG we can consider the partner to be the smaller one?)
It is also a pi-rotation of a pair of rims.
If it misses the differing box altogether then we are certainly OK,
so suppose it involves the differing box. 
It looks to me like if $e_i$ is removable then its pi-image is addable
to the corresp subpartition, and this addition would give the
reflection image of the partner. That is, the image of the partner is
dominant.
This is (surely!?) then also the image we would get according to the
partner's own hypercube...?! Right?! Surely!?...

2. What is the reflection that takes the partner to {\em its}
$\alpha$-image? 
Elsewhere we have claimed that (at least in some case) this is the
same reflection. Why?
If they are in the same facet the I think this is clear (?!).
So what about the different facet case?

}}

\mlem{\label{le:ijx}
Fix $\delta$. 
No pair of weights of form   
$\lambda$ and $\lambda-e_i+e_j$ are  in the same block
(unless $i=j$).
\\
That is, no pair of weights of form 
$\lambda+e_i$ and $\lambda+e_j$ are in the same block
(unless $i=j$).
}
\proof{
Such a pair cannot meet the charge-pair form of the balance condition
\CDMI,
since each of the skews involved has rank 1. \Qed
}


\mlem{ \label{le:ei11}
If $\sing_{\QQ}(\lambda)=\sing_{\QQ} (\lambda -e_i)$ then
for all pairs 
$(\mu,f_i(\mu)) \in [\lambda]_{\QQ}\times [\lambda -e_i]_{\QQ}$
\[
\proj_{\lambda} \ind \; \Dbt{n}{f_i(\mu)} = \Dbt{n+1}{\mu}
\]
\eql(eq:back)
\proj_{f_i(\lambda)} \ind \; \Dbt{n}{\mu} = \Dbt{n+1}{f_i(\mu)}
\eq
}

\proof{Note that the pair $(\mu,f_i(\mu))$ are adjacent by
  Theorem~\ref{th:embed}. 
For any $\nu$
\[
\ind \; \Dbt{}{\nu} = \left( \bigplus_j \Dbt{}{\nu+e_j} \right) \; 
               \bigplus \;  \left( \bigplus_k \Dbt{}{\nu-e_k} \right)
\]
For $\nu=f_i(\mu)$ adjacent to $\mu$, one of these summands is
$\Dbt{}{\mu}$. Specifically either (i) $\mu = \nu+e_l$ (some $l$);
or  (ii) $\mu = \nu-e_l$ (some $l$).

In case (i) other summands are of form $\mu-e_l+e_j$, $\mu-e_l-e_k$.
By Lemma~(\ref{le:ijx}) the former are not in $[\mu]_{\QQ}$, and since 
$\sing_{\QQ}(\lambda)=\sing_{\QQ} (\lambda -e_i)$
we may use Lemma~(\ref{le:combinME})(II) to exclude the latter. 
The other case is similar. \Qed
}

\footnot{{
\proof{
SKIP THE FIRST CASE AND DO THE SECOND!!!?... (maybe we did the first
case somewhere already?)

We start with the defining pair $(\mu,\nu)=(\lambda,\lambda-e_i)$.
We have 
\[
\Deltan{\lambda} = 
 \left( \bigplus_{x} \Deltan{\lambda-e_x}  \right)
 + \left( \bigplus_{y} \Deltan{\lambda+e_y}  \right)
\]
Evidently $\proj_{\lambda-e_i}\Deltan{\lambda}$  contains
$\Deltan{\lambda-e_i}$
from the left-hand sum (and no more, by Lemma~\ref{le:ijx});
and at most some $\Deltan{\lambda+e_y} = \Deltan{(\lambda-e_i)+e_i+e_y}$.
If $\lambda-e_i \sim^{\delta} (\lambda-e_i)+e_i+e_y$ then 
$ (\lambda-e_i)+e_i+e_y / \lambda-e_i$ is a \MiBS,
that is, the charges of $e_i$ and $e_y$ are equal and opposite. 
We next show that this is a contradiction.

Any candidate for $e_y$ must be below and right of the balance partner
($e_{i'}$ say) of $e_i$ in the original skew, and indeed 
 below and right of any box in the skew with the same charge.
But for such a box $e_y$ to be addable to $\lambda$ there must be a box in
$\lambda$ immediately above it, and one immediately left of it.
Only one of these can be in the skew (else one of them is the \ReRB,
rather than $e_i$), but then  the other one would lie `outside' 
(pinning) the skew --- a contradiction.

This is well illustrated by the examples in Figure~\ref{nominmin}.
In the first example in 
this Figure we suppose the unshaded boxes form a  \MiBS\ of $\lambda$,
with indicated removable box $e_i$  (charge 8).
The shaded -8 box indicated is in the 
necessary position for partner to $e_i$ in a
\oonn\ \MiBS\ {\em above} $\lambda$. But this is not addable without the 
other shaded boxes being in $\lambda$, and these would prevent the
unshaded set from being a skew.
The second example is similar. 

For the general $(\mu,\nu)$ case, suppose that $\proj_{\nu} \ind \Deltan{\mu}$ 
contains two summands. Then 

Note that $\mu$ is the image of $\lambda$ under some sequence of
reflections. 
We claim  that $\nu$ is the image of $\lambda-e_i$ under the 
{\em  same} sequence.
To see this note first that the step down from  
$\lambda-e_i \mapsto \alpha(\lambda-e_i)$
(where $\lambda \mapsto \alpha\lambda$ is the step down defining the
original skew) 
is the image of the same reflection
(in general a product of commuting reflections).
If the row containing $e_i$ contains other boxes in the skew 
then this is clear. If not, then the reflection factor associated to
$e_i$ and its balance partner {\em fixes} the relevant rows of
$\lambda-e_i$ (so could be omitted, but need not be).

 NOW WHAT ABOUT STEPPING DOWN AGAIN???!
$\alpha(\lambda-e_i)$ is now the larger weight, with $\alpha\lambda$
 the smaller. The difference is the rotation image of $e_i$, $e_{i'}$.
We can continue if this is contained in a \MiBS\ for 
$\alpha(\lambda-e_i)$.
BUT WHY SHOULD IT BE?!?
Well, it is removable 
(since $\alpha(\lambda-e_i)-e_{i'} = \alpha\lambda$). 
So actually just need to make sure that the \MiBS\ is not \oonn.

NONSENSE! $e_i$ might well not be in a MiBS at all!

Two cases:

(0) $e_i$ not in MiBS. WHAT TO DO HERE?

(I) $e_i$ is in MiBS, in which case 
just need to make sure that the \MiBS\ is not \oonn.

But this cannot happen since if $e_{i'}$ is one of these two boxes
then the other must be a box up and left from $e_i$. This will not be
removable in $\alpha(\lambda-e_i)$ since the presence of $e_i$ 
a \ReRB\ in
$\lambda$ implies the presence of one of the boxes above or left of
$e_i$ in  $\alpha\lambda$ and hence in $\alpha(\lambda-e_i)$.

IT REMAINS TO ARGUE THAT THESE TWO MUST HAVE AN IMAGE IN THE
DEFINING PAIR SETTING --- HENCE GIVING a contradiction.
It is reasonable that they have images, by the reflection
characterisation of the blocks. But how know that neither falls off
the edge of the dominant region?!!

MAYBE can use the f-map somehow?? NOT CLEAR.

\Qed (for 1st part)
}

}}

\section{The Decomposition Theorem}\label{ss:main th}

{\mth{ \label{th:decomposi}
For each $\delta\in\Z$ and $\lambda\in\Lambda$, 
the hypercube $h_{\delta}(\lambda)$ 
gives the  $\lambda$-th row of the 
$[\lambda ]_{\QQ}$-block of the 
global Brauer algebra $\Delta$-decomposition matrix $\Dmat$ over $\C$.
That is
\[
(P_n^{\delta}(\lambda)' : \Delta_n^{\delta}(\mu)' ) 
  = h_{\delta}(\lambda)_{\mu}
\]
for all $n \geq |\lambda|$;
or equivalently
\[
P_n^{\delta}(\lambda)' = \;
  \bigplus_{\mu\in h_{\delta}(\lambda)} \; \Delta_n^{\delta}(\mu)'
\]
(Recall we omit $\lambda=\emptyset$ in case $\delta=0$.
With this caveat Specht and standard modules coincide and we may
interpret the above either as Specht characters,
as required for the Cartan decomposition matrix;  
or as multiplicities in standard filtrations.)
\\
This data determines the Cartan decomposition matrix $\Cmat$ for any finite
$n$ by (\ref{eq:C=DD}). 
}}

\proof{
We prove for a fixed but arbitrary $\delta$, working by induction on
$n$.
The base cases are $n=0,1$, which are trivial
(and $n=2$ for $\delta=0$, which is straightforward).
We assume the theorem holds up to level $n-1$, 
and consider $\lambda\vdash n$.
}


The $\lambda$-th row of $D$ encodes the standard content of projective
module $\Pbt{}{\lambda}$. 
We apply the induction functor to a suitable $\Pbt{}{\lambda-e_i}$ in
level $n-1$ (known by the inductive assumption),
and use Prop.(\ref{pr:proj2}):
\[
\proj_{\lambda} \ind \; \Pbt{}{\lambda-e_i} 
    \cong \Pbt{}{\lambda} \bigoplus Q   
\]
(some $Q$).
Thus the  challenge is to determine 
the $\Delta$-content of 
$\proj_{\lambda} \ind \; \Pbt{}{\lambda-e_i} $  and $Q$.
In general determining $Q$  
can be complicated, but we will show that 
there is always a choice of $ \lambda-e_i$ which makes it tractable. 

Note that if $\lambda$ is at the bottom of its block then the claim
is trivially true. If  $\lambda$ is not at the bottom of its block
then
the binary sequence  
$b_{\delta}(\lambda)$ has at least one 01 (or initial 11) subsequence. 
Thus we can choose $e_i$ to be a removable box from the
skew associated to the corresponding edge $\alpha$ of $h_{\delta}(\lambda)$.
(We sometimes write $\mu=\alpha\lambda$ for the partition at the other end of this
edge, so the skew is $\lambda/\mu=\lambda/\alpha\lambda$.) 
Note that this skew is a \MiBS, by (\ref{pa:isMiBS}). 


The next step depends on 
whether the skew $\lal$ is of form (1)+(1), 
or otherwise.

\subsection{Properties of \MiBS s}

\mdef
We will say that such a skew $\lal$ is {\em boxy} if 
every box in it lies within a $(2^2)$-shape that also lies within the
skew. 
In our case, these are the skews in which the pair of rims fully overlap 
(i.e. run side-by-side).
Thus in our case boxy skews have a terminal  $(2^2)$-shape at each
end, in which the largest magnitude charges reside. 
Note that since no $(2^2)$-shape has a removable box
of largest magnitude charge, 
neither does a
boxy skew (on the other hand every such shape has a removable box of
next-largest magnitude, and one can see that the largest of these is
removable at one end of the boxy skew or the other). 
An example is given in Figure~\ref{fig:minskew eg}(iii).

If a minimal skew is neither of form (1)+(1) nor boxy we shall say that it is
generic. 


\mlem{
Let $\lambda/\mu = \mib{\lambda}{\alpha}$ 
be a minimal $\delta$-balanced skew.
Then there are a pair of boxes in the skew of greatest magnitude charge. 
In case the skew is of shape (1)+(1) 
both of these are removable; 
in the boxy cases (such as  $(2^2)$)
neither are removable (but precisely one of the
next-largest is removable);
and otherwise precisely one of them is removable. }

\proof{ All statements are (by now) clear except the last. For this 
note that if both were removable this would contradict
that  $\alpha\lambda$ is a \MaBS, 
since removing just this pair 
from $\lambda$ would give a larger \dBS; 
while if neither were removable then again this would contradict 
the \MaBS\ property, since removing the 
complement (i.e. the boxes in  
$\mib{\lambda}{\alpha}$
{\em not} in this pair) 
would give a larger \dBS.
\Qed}

\mdef
We call a removable box of largest
magnitude charge (among those removable in the given skew) a 
{\em rim-end removable box}. 
(Since the skew is a (possibly touching) pair of rims, and this box lies
at one of the outer ends.)




\mdef {{\label{exa:get}
Examples of minimal skews are shown in Figure~\ref{fig:minskew eg}. 
\begin{figure}
\[
(i) \; \;
\includegraphics[width=2.099in]{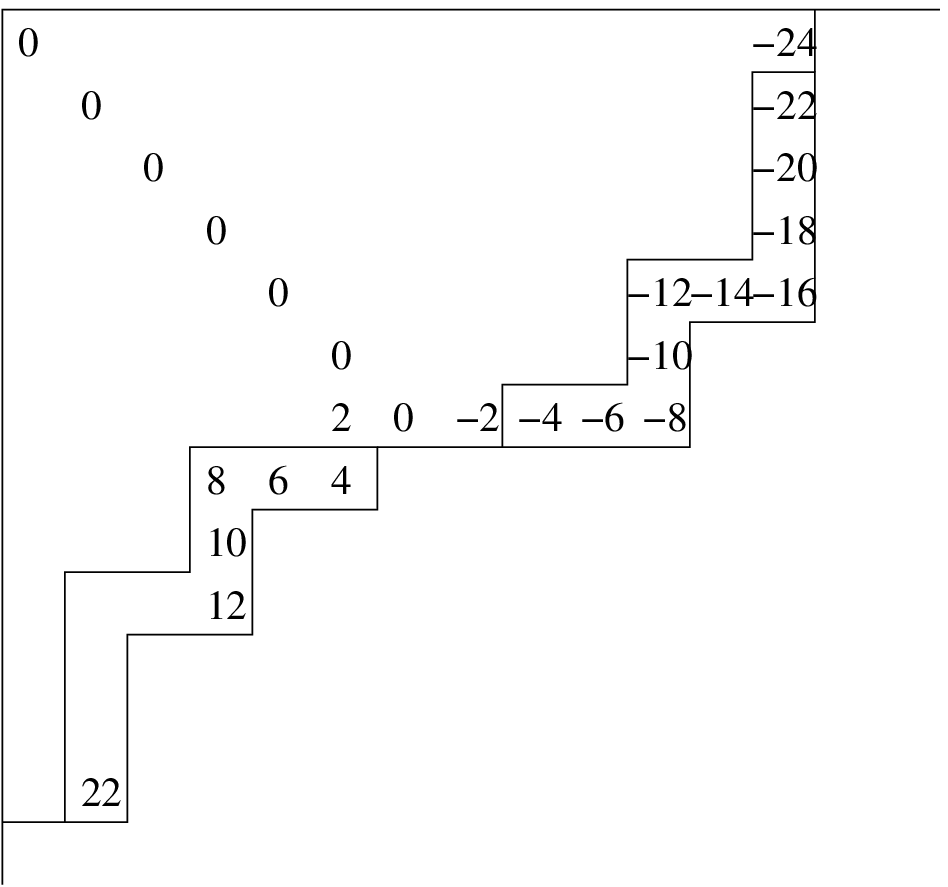}
\hspace{-.01in} (ii) \;\;
\includegraphics[width=2.099in]{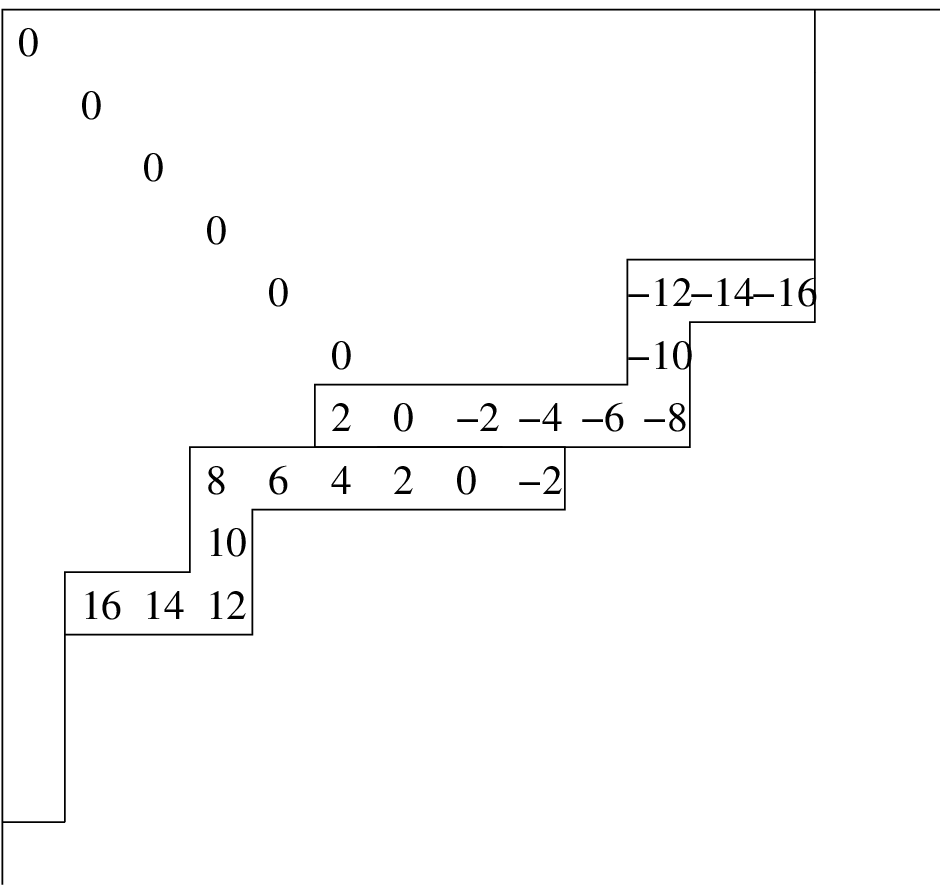}
\hspace{-.26in} (iii) \;\;
\includegraphics[width=1.1in]{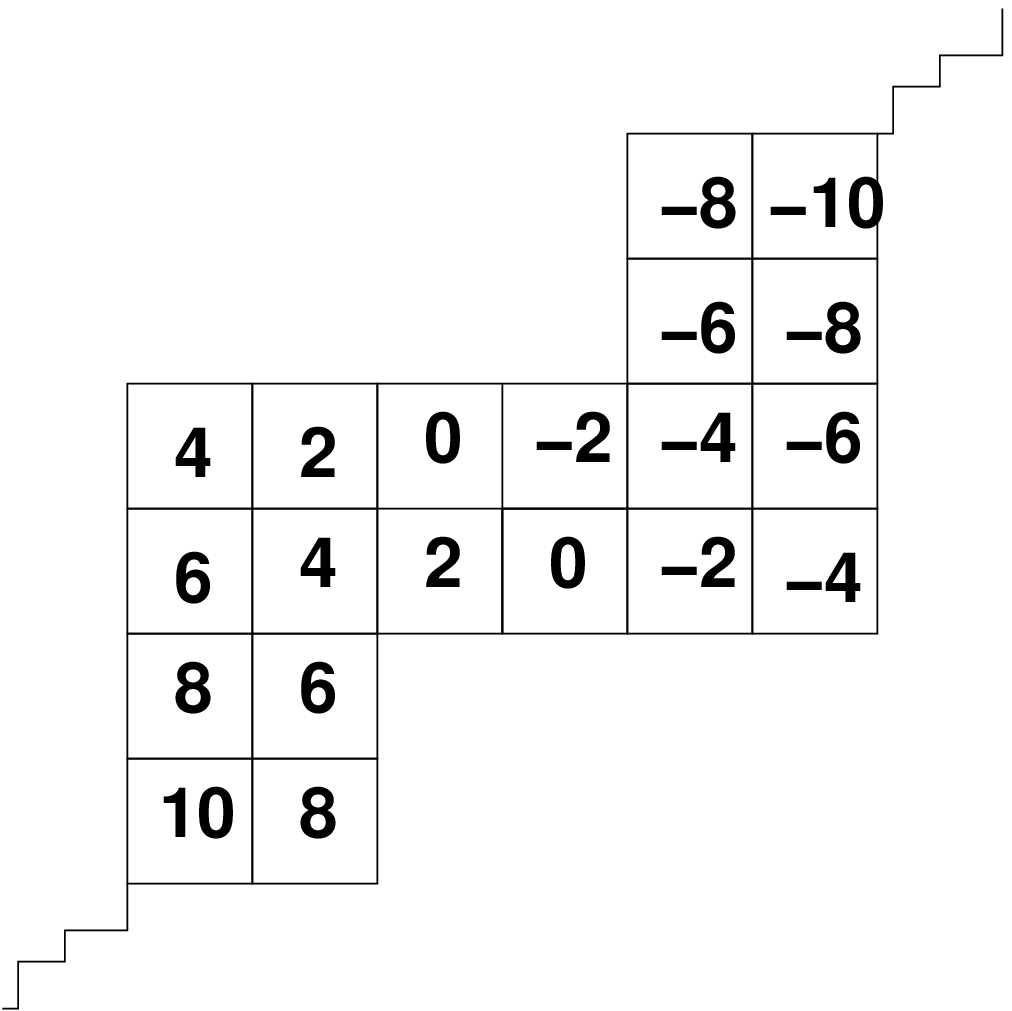}
\]
\caption{\label{fig:minskew eg} Examples of \MiBS s.}
\end{figure}
The rim-end removable boxes (as labelled by charge) in the figure are 
(i) 22; (ii) -16; (iii) 8.
\\
(For $\QQ=1$ example (i) is, in greater detail, 
\[
\lambda + \rho_1 = 
( 25/2, 23/2, 21/2, 19/2, 17/2, 11/2, 9/2, 
-3/2, -9/2, -11/2, -17/2, -19/2, -21/2, ...)
\]
which is five-fold singular (in the sense of (\ref{de:sing-set})), giving
$
\Psio_1(\lambda) = \{2,3 \} 
$
for its valley set.)
}}



\footnot{OLD:

{\mlem{ \label{le:facet-co}
Let $e_i$ be a \ReRB\ 
of $\lambda$.  
 If partitions $\lambda$ and $\lambda-e_i$ have the same 
singularity then  $\Psio_{\QQ}(\lambda) = \Psio_{\QQ}(\lambda-e_i)$.
}}
\proof{
It is clear that for weights differing by one box, if the singularity
is unchanged then either $\sis_{\QQ}$ is unchanged, or one pair is
swapped for another.

(I) Suppose  $\sis_{\QQ}(\lambda)=\sis_{\QQ}(\lambda-e_i)$.
To verify this case condider moving between $\lambda$ and $\lambda-e_i$ via 
$\lambda+\kappa e_i$, $0\leq \kappa \leq 1$.
This journey is too short for any hyperplane to be crossed
(by the strongly descending property), and the
same one cannot be left and joined by changing a single term in the
sequence. 

(II) Suppose
 $\sis_{\QQ}(\lambda) \neq \sis_{\QQ}(\lambda-e_i)$.
To see this case note that the conditions imply
$\lambda\sim(...,x+1,x,...,-x)$ and
 $\lambda-e_i\sim(...,x+1,x,...,-x-1)$, or similar. 
A different pair is discarded in each case by the $Reg$ part of the 
$\Psio_{\QQ}$ map, but since the singletons
retained are adjacent they pass to the same term in the magnitude
ordering $o$. 
\Qed
}

}



\begin{figure}
\includegraphics[width=6in]{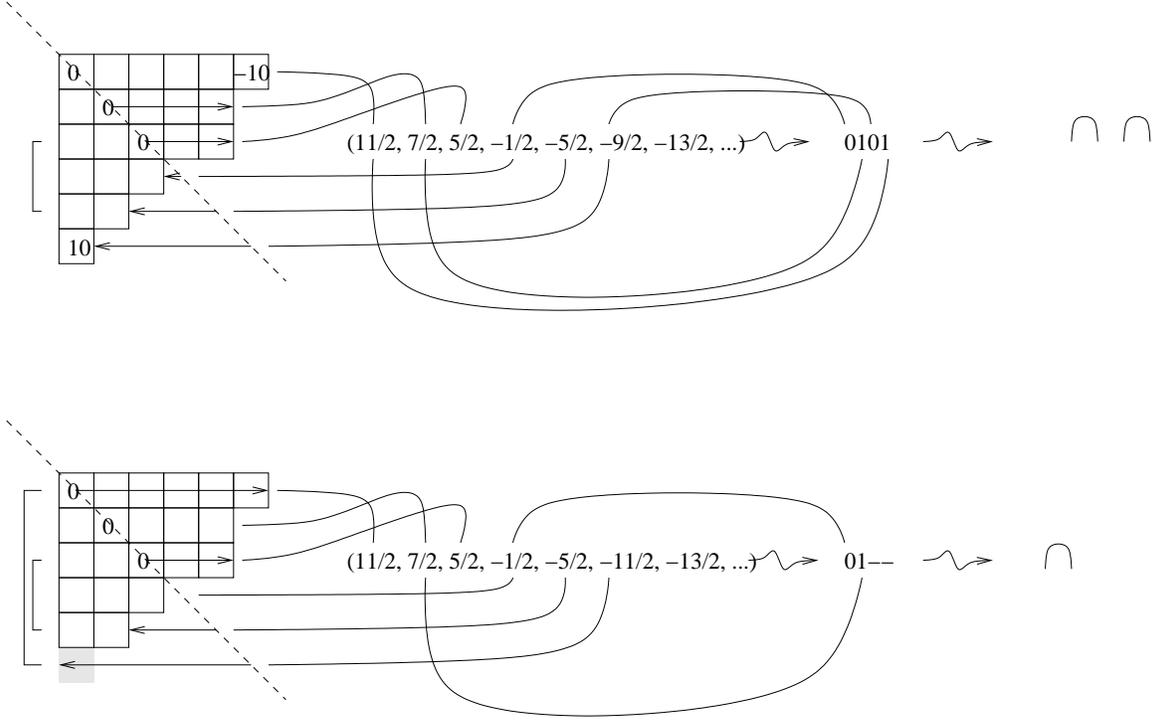}
\caption{\label{epos01} 
{Two examples showing the passage from 
Young diagram $\lambda$, via corresponding (case $\QQ =1$)
descending sequence $\lambda+\rho_1$, 
to binary  sequence and TL diagram. 
The connecting lines indicate the precise passage of data through the process.
The two cases are related in the form 
$\lambda$, $\lambda-e_i$, illustrating a step up in singularity
(the singular pairs of rows in each case are marked on the left).}}
\end{figure}

\subsection{Cases in the inductive step}

{\mpr{\label{prsing}
Fix $\QQ$, and hence an identification between valley sequences and
partitions. 
Pick $\alpha\in\Gammad_{,\lambda}$ and let $e_i$ be a 
rim-end removable box in
$\lal$. 
Then
the   singularities obey
$$
s_{\delta}(\lambda-e_i)= \left\{ 
\begin{array}{ll} 
s_{\delta}(\lambda)+1   & \mbox{ if } | \lal | =2
\\
s_{\delta}(\lambda)     & \mbox{ o/w }    \end{array}\right.
$$
}}

\medskip


\noindent {\em Proof:}
If $|\lal | =2$ we are in the $(1)+(1)$ or $(1^2)$ case,
and the charges in the boxes are (say) $x$ and $-x$. 
Removing $x$ (from row $i$) we get a row ending in charge $x+2$,
giving $(\lambda+\rho_{\QQ})_i = -\frac{x+2}{2}+\frac{1}{2}
=-\frac{x+1}{2}$ (by Lemma~\ref{le:ch3}). 
The row ending in $-x$ has 
$(\lambda+\rho_{\QQ})_j = -\frac{-x}{2}+\frac{1}{2}
=\frac{x+1}{2}$
thus these two rows are now a singular pair. 
\\
Figure~\ref{epos01} gives an example.

Suitable examples of the generic situation 
are given in Example~\ref{exa:get}. 
If the upper end of a rim ends in a row (of length greater than 1),
such as the upper rim in Example~\ref{exa:get}(ii), which ends in -16,
then the end box of this row is removable, but its balance partner is
not.
It follows that 
singularity is unchanged on removing the end-box $e_i$,
since this row becoming part of a singular pair would imply
a removable balance partner. 
(Thus 
$\Psio_{\QQ}(\lambda-e_i) = \Psio_{\QQ}(\lambda) $, 
indeed we remain in the same facet.)
\\
If   the lower end of a rim ends in a column (of length greater than 1),
such as the lower rim in (i), which ends in 22,
then the end-box of this column is removable. 
This time 
$\lambda-e_i$ lies on different hyperplanes to $\lambda$,
but {\em overall} singularity is unchanged. 
\\
(In the particular example $-21/2 \rightarrow -23/2$.)

In the case $(2^2)$  we have 
\[
\includegraphics[width=.4in]{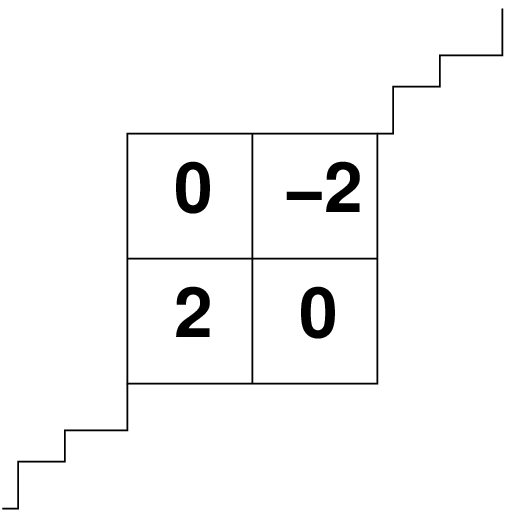}
\mapsto
(...,3/2,1/2, \leq -5/2,...)
\qquad \leadsto \qquad
\includegraphics[width=.4in]{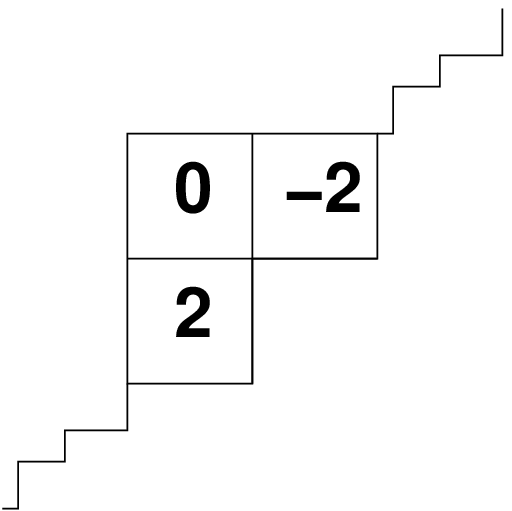}
\mapsto
(...,3/2,-1/2, \leq -5/2,...)
\]
which shows that the singularity does not change. 

For the remaining (boxy) cases there are a couple of analogous
variations to the generic `ends in row/column' cases treated above.
Here we merely illustrate with a couple of examples.
In the case $(2^4)$ we have 
\[
\includegraphics[width=.4in]{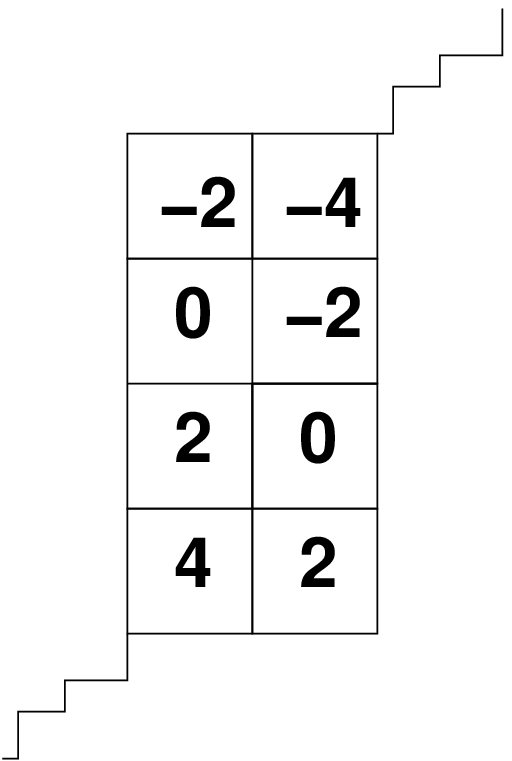}
\mapsto
(...,5/2,3/2,1/2,-1/2, \leq -7/2,...)
\qquad \leadsto \qquad
\includegraphics[width=.4in]{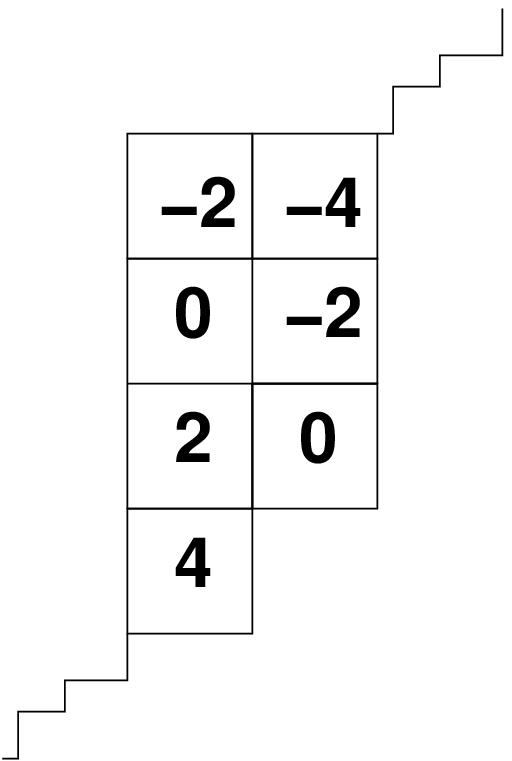}
\mapsto
(...,5/2,3/2,1/2,-3/2, \leq -7/2,...)
\]
which shows that the singularity does not change, 
although the wall does. 
In the case $(3^2)$ we have (similarly embedded, in general)
$(3^2) \mapsto (...,2,1, \leq -3,...)$ $\leadsto$ 
$(32) \mapsto (...,2,0, \leq -3,...)$
which has the same singularity 
(and wall set).
A more typical boxy skew is 
\[
\includegraphics[width=.75in]{xfig/boxTyp-1.eps}
\mapsto
(...,11/2,9/2,7/2,5/2,-5/2,-7/2, \leq -13/2,...)
\]
Removing the removable 8 here changes $-7/2 \rightarrow -9/2$,
giving the same singularity (different wall). 
\Qed


{\mpr{
Fix $\QQ$. 
Pick $\alpha\in\Gammad_{,\lambda}$ and let $e_i$ be a 
rim-end removable box in
$\lal$. 
In the cases in which the skew is neither $(1)+(1)$ nor $(1^2)$ 
\\
(i)
the standard decomposition pattern for $\Pbt{}{\lambda}$ is the 
{\em `translate'} of that for $\Pbt{}{\lambda-e_i}$: 
\[
\StDMt{P}{\lambda}{\mu}  = \StDMt{ P}{\lambda-e_i}{f_i(\mu)} 
\qquad \forall \mu \in [\lambda]_{\QQ}
\]
\\
(ii)
This verifies the inductive step for the main theorem in such cases.
That is, $h_{\QQ}(\lambda) \cong h_{\QQ}(\lambda-e_i)$. 
}}

\medskip


\noindent {\em Proof:}
Consider the  `translation' $\proj_{\lambda} \ind \; \Pbt{}{\lambda-e_i}$
of $\Pbt{}{\lambda-e_i}$.
By  Proposition~\ref{pr:pweight}
\[
\proj_{\lambda} \ind \; \Pbt{}{\lambda-e_i} \; = \; \Pbt{}{\lambda} \oplus Q
\]
with $Q=\proj_{\lambda} Q$ some projective, possibly zero. 
In the cases under consideration (skew neither  $(1)+(1)$ nor $(1^2)$) 
each standard module occuring in 
$\Pbt{}{\lambda-e_i}$ induces precisely one standard module after
projection onto the block of $\lambda$,
by  Lemma~\ref{le:ei11} (noting Proposition~\ref{prsing}). 
More specifically, writing 
\eql(eq:p10)
\Pbt{}{\lambda-e_i} = \; \bigplus_{\mu} \; c_{\mu} \; \Dbt{}{f_i(\mu)}
\eq
(for some multiplicities $c_{\mu}$), using (\ref{para:fi}); then 
\[
\Pbt{}{\lambda} \oplus Q 
=  \proj_{\lambda} \ind \; \Pbt{}{\lambda-e_i} \; 
= \; \bigplus_{\mu} \; c_{\mu} \;\proj_{\lambda} \ind \; \Dbt{}{f_i(\mu)}
= \; \bigplus_{\mu} \; c_{\mu} \; \Dbt{}{\mu}
\]


On inducing again and projecting back to
the block of $\lambda-e_i$, by (\ref{eq:back})
we have
\[ 
\proj_{\lambda-e_i} \ind \; (\Pbt{}{\lambda} \oplus Q) \; 
=  \; \bigplus_{\mu} \; c_{\mu} \; \Dbt{}{f_i(\mu)}
\] 
That is, 
each standard module occuring in 
$( \Pbt{}{\lambda} \oplus Q)$ induces precisely one standard module after
projection onto the block of $\lambda-e_i$. 
Comparing with (\ref{eq:p10}),
it follows that this second `translation'  
may be identified with 
$\Pbt{}{\lambda-e_i}$ again. 
Since this is indecomposable,
the {\em first} translation cannot be split, and hence is precisely
$\Pbt{}{\lambda}$ --- with the same decomposition pattern. 

For the last part use (\ref{de:isomg}). 
\Qed


The remaining cases needed to move between level $n$ and $n-1$ 
are skews of form (1)+(1).


{\mpr{ \label{pr:fin}
Fix $\QQ$. 
Pick $\alpha\in\Gammad_{,\lambda}$ and let $e_i$ be a 
rim-end removable box in
$\lal$. 
Then
in the cases in which the skew is of form $(1)+(1)$ or $(1^2)$
\\
(I)
the sequence 
$ \Psiob_{\QQ}(\lambda) =   
\bump{\alpha} \Psiob_{\QQ}(\lambda-e_i) $. 
Thus, hypercube 
$h_{\QQ}(\lambda) =   \double{\alpha}   \bump{\alpha} h_{\QQ}(\lambda-e_i)$ 
(i.e. has increased `dimension' by +1
compared to $h_{\QQ}(\lambda-e_i)$).
The sequence  
$ \Psiob_{\QQ}(\alpha\lambda) =\bumpo{\alpha} \Psiob_{\QQ}(\lambda-e_i)$ 
(i.e. differs from   
$ \Psiob_{\QQ}(\lambda-e_i) $ 
by insertion of subsequence 10 in the $\alpha$ position). 
\\
(II)
the standard decomposition pattern for $\Pbt{}{\lambda}$ is 
in agreement with the above,
in the sense of the equality in the main theorem:
$(\Pbt{}{\lambda}:\Dbt{}{\mu}) = h_{\delta}(\lambda)_{\mu}$ 
(all $\mu$)
if 
\newcommand{\mup}{\nu}%
$(\Pbt{}{\lambda-e_i}:\Dbt{}{\mup}) = h_{\delta}(\lambda-e_i)_{\mup}$
(all $\mup$). 
}}

\noindent {\em Proof:} 
(I) As shown in the proof of Prop.~\ref{prsing} (or see below), 
removing $e_i$ from $\lambda$ makes that
row part of a singular pair with the row containing the box with
opposite charge. Thus  $ \Psiob_{\QQ}(\lambda-e_i) $ 
differs from $ \Psiob_{\QQ}(\lambda) $ 
in that a pair which contributed an 01 sequence in the latter does not
contribute to the valley sequence in the former
--- i.e.  $ \Psiob_{\QQ}(\lambda-e_i) $ differs
by the removal of this 01 sequence. 
(Figure~\ref{epos01} serves as an example here.) 
It remains to confirm the {\em position} of the removal.
The situation $\lal \sim (1)+(1)$
is well illustrated by the following generic example:
\[
\includegraphics[width=2.16in]{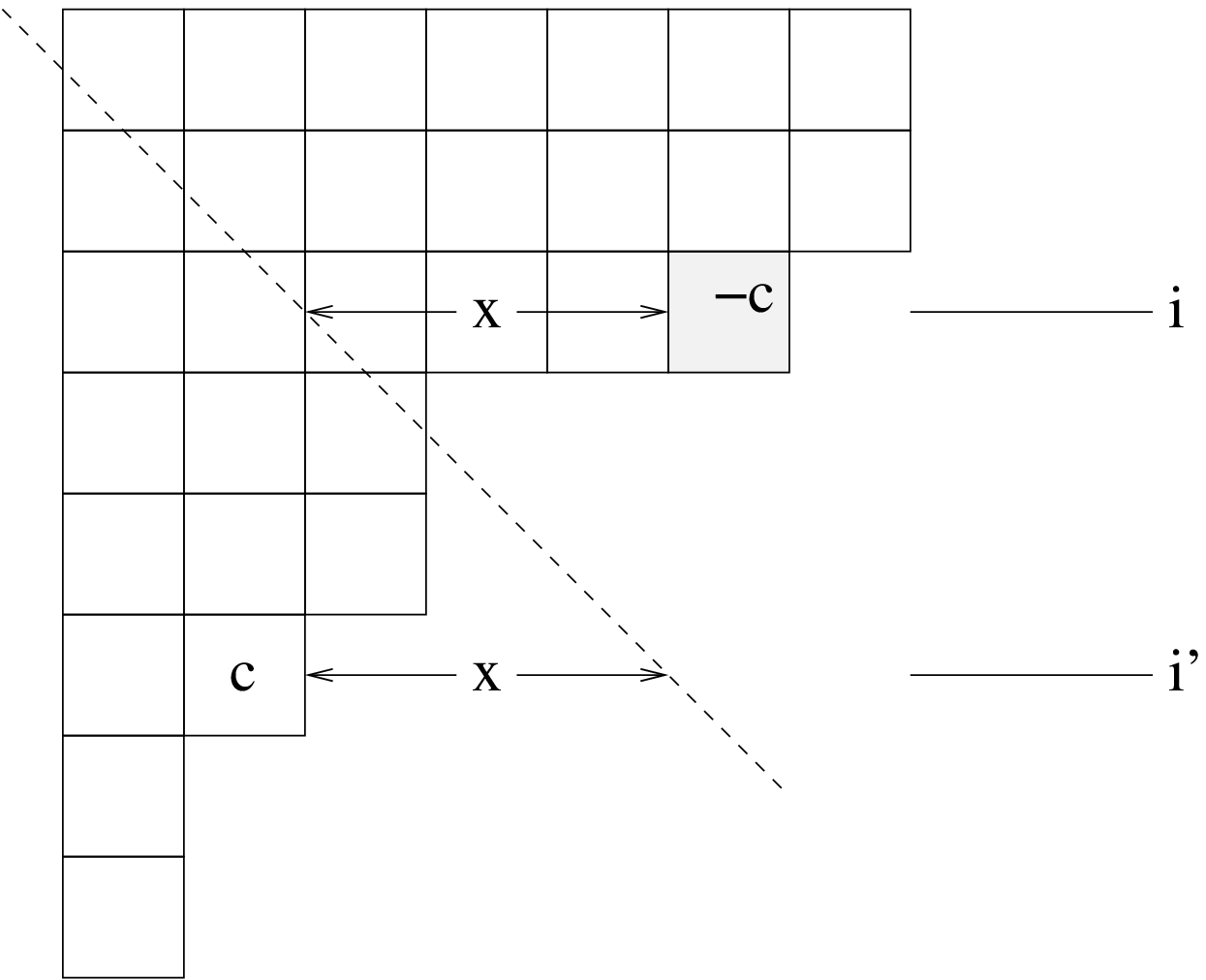}
\]
In general we have 
\[
\lambda +\rho_{\QQ} \sim (...,\underbrace{\underbrace{x+1}_i,...,-x}_{01},...)
\]
Altogether the bracketed pair contribute an 01 in binary as indicated.
The $x+1$ lies at some position $l$, say, in the {\em magnitude} order,
depending on the rest of $\lambda$. 
Confer
\[
\lambda-e_i +\rho_{\QQ} \sim (...,\underbrace{x}_i,...,-x,...)
\]
Here the $x,-x$ are a singular pair, so do not appear in the
magnitude order --- 
to obtain its binary representation from that of
$\lambda$ one deletes the binary pair 01 in the $l-1,l $ 
position. That is, $ \Psiob_{\QQ}(\lambda) =   
\; \bump{{l\!\! -\!\! 1}} \; \Psiob_{\QQ}(\lambda-e_i) $.
Finally 
\[
\alpha\lambda +\rho_{\QQ} =
\lambda-e_i-e_{i'} +\rho_{\QQ} 
\sim (...,\underbrace{\underbrace{x}_i,...,-x-1}_{10},...)
\]
Since the $\alpha$ action on $\lambda$ manifests 
(by definition) as $10 \leftrightarrow 01$ in the
$\alpha,\alpha+1$ position of $\Psiob_{\QQ}(\lambda)$ 
we see that position $l-1=\alpha$ as claimed.
The other assertions follow immediately.


\medskip

(II) 
Applying $\proj_{\lambda} -$ to Proposition~\ref{pr:ind1}(ii) here we get
a short exact sequence
\[
0 \rightarrow \Dbt{}{\lambda-e_i-e_{i'}} \rightarrow 
 \proj_{\lambda} \ind \; \Dbt{}{\lambda-e_i} \rightarrow \Dbt{}{\lambda} 
\rightarrow 0
\]
(non-split, by \cite[Lemma~4.10]{\CDMi}). That is 
\eql(p-dd)
\proj_{\lambda} \ind \; \Dbt{}{\lambda-e_i} \;
 = \Dbt{}{\lambda} + \Dbt{}{ \lambda-e_i - e_{i'}}
 = \Dbt{}{\lambda} + \Dbt{}{ \alpha\lambda}
\eq
(non-split). 
Translating $P_{\lambda-e_i} \; := \; \Pbt{}{\lambda-e_i}$
away from and then back to $\lambda-e_i$
therefore produces 
a projective whose dominating content is two copies of 
$\Dbt{}{\lambda-e_i}$   
(one from each of the summands on the right of (\ref{p-dd})).
Indeed every $\Delta$-filtration factor of $P_{\lambda-e_i}$ 
engenders at most two factors in 
$\proj_{\lambda-e_i} \ind \; ( \proj_{\lambda} \ind \; P_{\lambda-e_i} ) $
(we shall be able to make a precise statement shortly).
Hence, by (\ref{pr:proj1}),
$\proj_{\lambda-e_i} \ind \; ( \proj_{\lambda} \ind \; P_{\lambda-e_i} ) 
    \; = \; P_{\lambda-e_i} \oplus P_{\lambda-e_i}$.
It follows that 
\[
\proj_{\lambda} \ind \; P_{\lambda-e_i} \; = \;  P_{\lambda} 
\]
It remains to show that 
$\StDMt{\proj_{\lambda} \ind P}{\lambda-e_i}{-} = h_{\QQ}(\lambda)$
(given $\StDMt{P}{\lambda-e_i}{-} = h_{\QQ}(\lambda-e_i)$).
%

For each $\Delta_{\mu}$ occuring in
the $\Pbt{}{\lambda-e_i}$ decomposition we will see that the translation is
$\Delta_{\mu} \leadsto \Delta_{\mu+} + \Delta_{\mu-}$
for some pair $\mu+,\mu-$ in the $\lambda$-orbit.
For $\lambda-e_i$ itself we have seen 
in the proof of (I)  
that $b_{\QQ}(\lambda-e_i)$ gives 
 $b_{\QQ}(\lambda)$ and  $b_{\QQ}(\alpha\lambda)$ by inserting 01 
(respectively 10) in the $\alpha$ position.
For other $\mu \in h_{\QQ}(\lambda-e_i)$, 
note that the relevant singular pair of rows in
$\lambda-e_i$, while not contributing to the magnitude order (since
they are singular) are formally permuted (in the $\Dgroup$-action sense)
along with the rest of the rows, in the collection of reflection group
actions that traverse $h_{\QQ}(\lambda-e_i)$. 
Thus they (jointly)
maintain a formal position in the magnitude order, between
two terms that are properly consecutive in this order. 
The difference with $\mu_+,\mu_-$ is that in these one of the pair is
extended by 1, or contracted by one. Thus the singularity is broken,
and the pair appear properly in the order, between the given two
terms, and hence bumping up the larger of the two.
Since $\mu+\rho_{\QQ}$ is just a signed permutation of 
$\lambda-e_i+\rho_{\QQ}$
(and hence just a permutation, as far as the magnitudes are concerned), 
the position of the pair in the
magnitude order, and hence 
the position of the bump in the binary representation, is at $\alpha$,
the same as for $\lambda-e_i$.
That the collection thus engendered overall is $h_{\QQ}(\lambda)$ now
follows directly from Equation(\ref{eq:hyp lemma}).
Indeed, for $\mu \in h_{\QQ}(\lambda-e_i)$, and 
$\bumpo{\alpha}\mu, \bump{\alpha}\mu$ the two partitions associated to
$\mu$ by the doubling 
$h_{\QQ}(\lambda) =   \double{\alpha}   \bump{\alpha} h_{\QQ}(\lambda-e_i)$, 
we have (non-split \cite[Lemma~4.10]{\CDMi})
\[
0 \rightarrow \Dbt{}{\bumpo{\alpha}\mu} \rightarrow 
 \proj_{\lambda} \ind \; \Dbt{}{\mu} \rightarrow \Dbt{}{\bump{\alpha}\mu} 
\rightarrow 0
\]
(From an alcove geometric perspective one may view this argument as
follows: 
Since $[\lambda-e_i]_{\QQ}$ is a 
strictly more singular orbit than
$[\lambda]_{\QQ}$ the reflection group elements moving through
$h_{\QQ}(\lambda-e_i)$ will also serve to move the pair
$\lambda,\alpha\lambda$ through these pairs $\mu+,\mu-$,
thus they remain adjacent above and below $\mu$.)
\\
\Qed


\mdef
Example for Proposition~\ref{pr:fin}: 
$\QQ=1$, computing for $\lambda=4422$ via 
$\lambda-e_2 = 4322$. We have
\[
\includegraphics[width=1.16in]{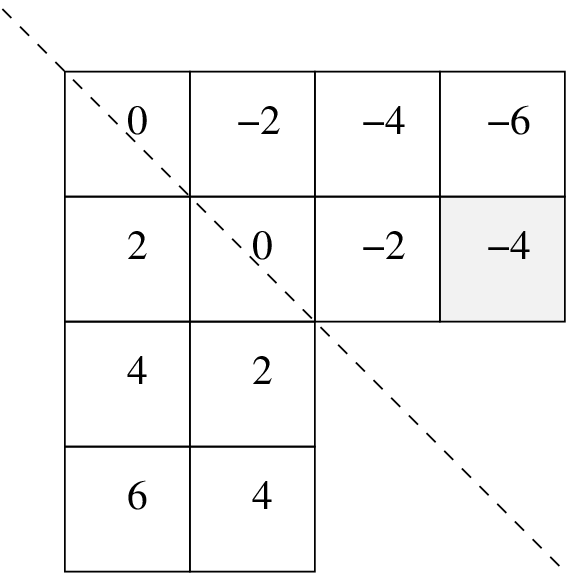}
\]
In particular $e_1(4322)= (7/2,3/2,-1/2,-3/2,...)$ 
so $\Psio_{1}(4322)=\mbox{toggle}(\{ 2\})=\{ 1,2\}$.
By the inductive hypothesis we have
\[
\StDM{P}{4322}{-} \;\; =\;\;  h_{1}(4322)_- \;  
= \raisebox{.31in}{ 
\xymatrix{ 4322 \ar@{-}[dr] \\ & 221} }
\cong 
\raisebox{.31in}{ \xymatrix{ 12 \ar@{-}[dr]_{\underline{12}} \\ & \emptyset} }
\cong 
\raisebox{.31in}{ \xymatrix{ 01 \ar@{-}[dr] \\ & 10} }
\]
Here the first form of the hypercube is in partition labelling;
the second form is in $P(\N)$ labelling (having applied the toggle);
and the last is the untoggled binary representation.
Note that 
we have reverted to the untoggled form at the last since we will be inserting
an 01 subsequence (removing the need for the toggle) at the next step.
Translating off the wall we get 
$4322+221 \rightarrow (4422+4321)+(321+22)$.
In binary this corresponds to 
$01 \rightarrow 0**1 \rightarrow 0101+0011$
and $10 \rightarrow 1**0 \rightarrow 1100+1010$.
These four sequences therefore encode the content of $P_{4422}$. 

Meanwhile 
\[ h_{1}(4422) 
= \raisebox{.31in}{ 
\xymatrix{ & 4422 \ar@{-}[dl]\ar@{-}[dr] 
\\ 4321\ar@{-}[dr] && 321\ar@{-}[dl] \\ & 22} 
}
\cong 
\raisebox{.31in}{ 
\xymatrix{ & 34 \ar@{-}[dl]_{23} \ar@{-}[dr]_{14} 
\\ 24 \ar@{-}[dr] && 13 \ar@{-}[dl] \\ & 12} 
}
\cong 
 \raisebox{.31in}{ 
\xymatrix{ & 0011 \ar@{-}[dl]\ar@{-}[dr] 
\\ 0101\ar@{-}[dr] && 101 \ar@{-}[dl] \\ & 11} 
}
\]
confirming the assertion of the Theorem in this case.
\\
Note how the insertion of a binary pair in the $\alpha$ position,
and action of $\alpha$ on that pair, 
transforms $h_{1}(4322)$ 
to produce $h_{1}(4422)$. 
The effect is (i) to
extend the hypercube by a new generating direction (labelled by $\alpha$);
(ii)
the generating
edge inherited from  $h_{1}(4322)$  changes label from 12 to 14
due to the bump (which illustrates how such non-$\Ge$ edge labels
arise in this contruction).


\medskip

Proposition~\ref{pr:fin} completes the main inductive step for the Theorem. 
\Qed




\newcommand{\E}{{\mathbb E}}
\newcommand{\EN}{\R^{\N}} 
\newcommand{\Valley}{Valley}

\section{Background: parabolic \KL\ polynomials}
\label{ss:pKLp}
In the remainder of the paper we explain where the idea for
hypercubical decomposition graphs comes from.
\newcommand{\ph}{p} 
\newcommand{\PH}{P} 

Associated to each Coxeter system $C$ and parabolic $A$,
acting as reflection groups on a suitable space,  
is an alcove geometry on that space. 
For each such pair $C/A$ there is, therefore, 
an array $\PH = \PH(C/A)$ of \KL\ polynomials
--- one for each ordered pair of alcoves. 
 (Deodhar's recursive formula
\cite{Deodhar87} computes these polynomials in principle. 
However it generally tells us very little about them in practice.)
These polynomials are of interest from a number of points of view.
For example they are often important in representation theory 
(see \cite{Soergel97a,MartinWoodcock03}
and references therein).
So, with the reflection group pair $\Dgroup/\Dgroup_+$ manifesting itself in Brauer
algebra block theory (as we have seen),
one is motivated to compute them in this case.

\subsection{Chamber geometry} \label{ss:Hum}
\newcommand{\AX}{{\mathcal A}}
\newcommand{\CX}{{\mathcal C}}
\newcommand{\CO}{C_0}  

We first need to  
review the notion of chamber geometry. 
In this we follow Humphreys \cite{Humphreys90}. 
(Alcove geometry is a mild generalisation associated to the
group/parabolic pair. Humphreys introduces this in the context of
affine extensions, but it serves equally well in general.)

Let $V$ be a Euclidean space, and $(W,S)$ a Coxeter system
with an action generated by reflections on $V$. 
Let $H_{s}$ be the reflection hyperplane of $s \in S$,
or indeed of any reflection  $s \in W$ generated by these.
For $T$ any subset of $S$ let $[T]$ be the set of reflections
generated by $T$.
Set
\[
\Hyp_T = \bigcup_{t \in [S]\setminus [T]} H_t
\]
A chamber is a maximal connected component of 
$V \setminus \Hyp_{\emptyset}$.  
Write $\CX_W$ for the set of chambers.

The set 
$H_t' = H_t \setminus \Hyp_{\{ t \}}$
(the subset of hyperplane $H_t$ that intersects no other hyperplane)
may similarly be broken up into  connected components.
At most one of these components intersects any given 
chamber closure $\overline{C}$.
If $H_t'$ intersects   $\overline{C}$ in this way it is called a wall
of $C$. 

For any given $C$, 
the set 
$\{ t \; : \; H_t' \cap\overline{C}\neq \emptyset \}$
of $t$s that make up its walls   functions as a
choice of $S$ in $W$
(i.e. they are an equivalent choice of Coxeter generators to 
the original set $S$). 
On the other hand $S$ may or may not determine such a $C$ uniquely.

The choice of a prefered chamber $\CO$  corresponds to 
the choice of a simple system
in $V$, and the associated reflections are simple reflections. 
(Given a non-commuting pair of these, the 
conjugate of one by the other is also a reflection,
but not `simple' in this choice.)

A reflection $s$ in $W$ is simple for chamber $B$ if its hyperplane
$H_s$ makes a wall of $B$
(NB simple for $B$ is not the same as simple, unless $B=\CO$).
For our purposes it will be convenient to think specifically of the
intersection of the hyperplane with the chamber closure 
(i.e. this facet) as the wall
(thus we distinguish the walls of distinct chambers in general, even if
they come from the same hyperplane).


\mdef \label{para:reflectx}
The reflection action of 
$W$ acts to permute  $\CX_W$. 
This action is transitive 
and indeed regular (simply transitive). 
See for example \cite[\S1.12]{Humphreys90}.
\\
Note that $W$ does {\em not}  act transitively on $V$, or 
specifically, on the set of walls.
The walls of $\CO$ are representatives for the $W$ orbits of the set
of all walls. 

Regularity says that we may identify $\CX_W$ with $W$, 
and the action of $W$ with the left-action on itself. 
In particular write
\eql(idW)
A = w_A \CO
\eq
(so we may indentify $\CO$ with $1$). 

Note that it follows from this identification 
that there is another commuting 
action of $W$ on $\CX_W$,
corresponding to the right-action of $W$ on itself. 

Noting the choice of $\CO$, 
define a length function on $\CX_W$: $l_W(A)$ is the number of hyperplanes
separating $A$ from $\CO$. 
(If $W$ is clear from context we shall write simply $l=l_W$.) 


\mdef \label{para:r-action}
We define a digraph $G(W,S)$ with vertex set $\CX_W$ by 
$(A,B)$ an edge if $B=tA$ with $t$ simple for $A$ and $l(B)=l(A)+1$. 

We call $t$ the left-action label of edge $(A,tA)$. 

By (\ref{idW})
the edge $(A,tA)$ may also be written $(w_A \CO, t w_A \CO)$.
The image under $w_A$ of a particular `initial' edge $(\CO, s\CO)$ ($s \in S$)
is $$
(w_A \CO, w_A s \CO)  
=  (w_A \CO , w_A s w_A^{-1} w_A \CO)
=  (A, w_A s w_A^{-1} A)
$$
Using the right-action this can be expressed as
$$
(w_A \CO, w_A s \CO)  = (w_A \CO, w_A \CO s) = (A,As)
$$
We call this $s$ the right-action label of the edge. 
(With this label the graph is essentially the right Cayley graph
 $\Gamma(W,S)$, and $s$ is the `colour' label.)

Evidently $G(W,S)$ is a rooted acyclic digraph, with root $\CO$.

\mdef
Let $v \in \CO$, and let $Wv$ be the $W$-orbit of $v$ in $V$. 
In the same way as above we may associate a graph to this orbit.
It will be evident that this graph is isomorphic to $G(W,S)$,
for any such $v$. 

\subsection{Alcove geometry}

Let $(W',S')$ be a system containing $(W,S)$ as a parabolic subsystem,
with both acting on $V$. 
The chambers of $W'$ are then called alcoves. Thus the alcoves are a
further subdivision of the chambers of $W$. 
Write $\AX=\CX_{W'}$ for the set of alcoves, and 
$X^+$ for the set of alcoves lying in $\CO$. 
Thus $X^+$ is a representative set for the $W$-orbits of $\AX$. 
(In this setting we will call any $v \in \CO$ {\em dominant}.) 

Choose $C'$ a prefered alcove in $\CO$. 
As before, the hyperplanes bounding $C'$ determine $S'$ 
(a superset of $S$, by the inclusion in $\CO$).

The digraph $G(W',S')$ has vertex set $\AX$, and  
$(A,B)$ an edge if $B=sA$ with $s$ simple for $A$ and 
$l_{W'}(B)=l_{W'}(A)+1$. 
This is evidently a rooted acyclic digraph, with root $C'$.
The edges are in correspondence with the set of walls,
and may thus be partitioned into $W'$-orbits,
labelled by the walls of $C'$. 

\newcommand{\Ga}{G_a}

\mdef
Let $\Ga = \Ga(W',W)$ denote 
the full subgraph of $G(W',S')$
with vertex set $X^+$. This is still rooted. 
Thus  
any alcove $A \in X^+$ may be reached from $C'$ by a sequence of
simple reflections, always remaining in $X^+$. 
\\
We shall denote the poset defined by the acyclic digraph $\Ga$ as 
$(X^+,<)$. 

The array $\PH = \PH(W'/W)$ is a (generally semiinfinite) lower unitriangular
matrix, with row and column positions indexed by $X^+$. 
It is natural to organise this data into rows 
(although it is also of interest 
to organise it into columns). 
These rows are thus `finite' (i.e. of finite support),
while the columns are not in general. 

\subsection{The recursion for $\PH(W'/W)$} \label{ss:pKLp def}

The recursion for rows of $\PH$ above the root in the poset
(acyclic digraph) order may be given as follows
(see \cite{Soergel97a} for equivalent constructions). 
Write $\PH = (\ph_{AB})_{A,B \in X^+}$. 
To compute the row $\ph_A$ for alcove $A$ we first 
compute another 
polynomial for each alcove $D$, $\ph'_{AD}$, also denoted $\ph'_A(D)$
as follows. 
(Actually $\ph'_A(D)$ can depend on the choice made next in the computation,
but $\ph_A$ does not and we supress this dependence in notation.)
\\
Pick an edge $(B,A)$ in $\Ga$ ending at $A$ (so $\ph_B$ is known).
For each alcove $D$ 
let $\Gamma_D^{\pm}$ be the set of alcoves $D'$ of $\Ga$ 
such that $(D',D)$ (resp. $(D,D')$) is an edge in the orbit of the edge
$(B,A)$.
(By (\ref{para:r-action}) we can express $(B,A)=(B,Bs)$, $s \in S'$,
whereupon any such $D'$ must obey $(D',D)=(D',D's)=(Ds,D)$
(respectively $(D,D')=(D,Ds)$).%
)
Then 
\eql(eq:pKL recursion)
\ph'_A(D) =  
   \sum_{D' \in \Gamma_D^+} ( v\ph_B(D) +\ph_B(D') )
+  \sum_{D' \in \Gamma_D^-} ( v^{-1}\ph_B(D) +\ph_B(D') )
\eq
(As noted there is at most one edge in the orbit of $(B,A)$ involving 
any alcove $D$. Thus at most one of these sums is non-trivial,
and that contains only one entry.
In particular $(B,A)$ is in its own orbit, so 
$\ph'_A(A) = v^{-1} \ph_B(A) + \ph_B(B) = 1$.) 

To obtain the row of $\PH$ that we want from $\ph'_A$ 
it is then necessary to perform a
subtraction in case the evaluation
$\ph'_A(D)(v=0)$ is non-zero for any $D<A$:
\[
\ph_A = \ph'_A - \sum_{D<A} \ph'_A(D)(v=0) \; \ph_D
\]
(But we shall see that the sum always vanishes in the case we are
interested in. So in our case $\ph_A = \ph'_A$.)

In order to work with this rule in any given alcove geometry 
it is necessary to be able to manipulate the graph $\Ga$ and its edge
orbits efficiently. 
In Section~\ref{S:D/A} we 
set up the requisite machinery for the case $\Dgroup/\Dgroupp$. 



\section{The reflection group action $\Dgroup$ on  $\EN$}
\label{S:D/A}

\newcommand{\axial}{axial} 
\newcommand{\ai}{\alpha}
\newcommand{\aop}[1]{\langle #1 \rangle}
\newcommand{\aopu}[1]{\langle \underline{#1} \rangle}
\newcommand{\Dd}{{\mathcal D}}

Define
$
v_- = o_2(\emptyset) =  -(1,2,3,...) \; \in \EN . 
$
In Section~\ref{ss:geom3} we chose the  alcove
 containing $v_-$ as $C'$, for the reflection group $\Dgroup$.
(We shall refer to  $\Dd v_-$ as the {\em fully-regular} orbit.)

\noindent
Thus our choice of $C'$ corresponds to choosing 
$S_{\Dgroup} = \{ (12)_-, \; (i \; i\! +\! 1) \}_{i \in \N}$
for the Coxeter generating set of $\Dgroup$. 

\medskip


The orbit $\Dd v_-$ consists in the set
of {\em co-even permutations}
(signed permutations of $v_-$ with an even number of positive terms).
By (\ref{para:reflectx}) 
this orbit (and hence each of the others) 
is isomorphic, via the left action of $\Dgroup$ upon it, 
to the (limit) regular representation.
It is easy to check that the action we are using is the left-regular
action.
By (\ref{para:r-action}) 
it is the associated right action that we need to determine in order to
compute (\ref{eq:pKL recursion}).
This commuting right action
corresponds to signed permutations of the {\em entries} in the
sequence, rather than signed permutations of the {\em positions}. 
For example
\[
(4,3,-1,-2,-5,...) (45) = (5,3,-1,-2,-4,...)
\]


\newcommand{\Gee}{G_e} 


Via the isomorphism 
between $V(v_-)$ and  
$\Pow_{even}(\N)$      
we understand left- and right-actions of 
$w \in \Dgroup$ on any $\lambdaa \subset \N$
(noting that $w \lambdaa$, respectively $\lambdaa w$, 
is not necessarily expressible in $\Pow(\N)$,
since it is not necessarily dominant).
When  $\lambdaa w$ {\em is} dominant we shall see now that 
the right-action transformation $\lambdaa \rightarrow \lambdaa w$ is 
expressible in a simple form in $\Pow(\N)$
which 
facilitates computation of the pKLps. 
Let $\Gee$ denote the simple relabelling of $\Gal$ from $\Pow(\N)$
using the above isomorphism.
(We shall shortly be able to identify $\Gee$ with $\Ge$.) 
The following crucial result is routine to show.


{\theo{\label{th:crux}
Let $\lambdaa \subset \N$ of even degree. 
Then there exists an edge $(\lambdaa,\lambdaa (\ai,\ai+1))$
in $\Gee$ iff $\lambdaa\cap\{\ai,\ai+1\} = \{ \ai \}$,
whereupon $\lambdaa (\ai,\ai+1) \cap\{\ai,\ai+1\} = \{ \ai+1 \}$;
and an edge  $(\lambdaa,\lambdaa (12)_-)$
in $\Gee$ iff $\lambdaa\cap\{1,2\} = \emptyset$,
whereupon $\lambdaa (12)_- \cap\{1,2\} = \{1,2\}$.
Every edge is one of these types. 
\Qed
}}

That is, we may associate edge labels corresponding to 
the right-action 
in $\Gee$,
taken from the Coxeter generating set $S_{\Dd}$
(as required by (\ref{para:reflectx})). 
To streamline still further we may write simply $\ai$ 
as `right-action' label for edges of form  
$(\lambda,\lambda (\ai,\ai+1))$ and $12$ for $(\lambda,\lambda(12)_-)$.
This makes explicit the identification with $\Ge$. 
See Figure~\ref{fig:valley graph}.

{\mrem{
The left-action labels are of course different in this regard.
Only elements of form $(ij)_-$ preserve dominance. 
}}





A convenient summary of the above is as follows 
(when we speak of an edge orbit on $\Ge$ we shall mean the orbit
induced by the graph isomorphism with $\Gal$ from the edge orbit
thereon):

{\theo{\label{th:big edge}
Two edges in 
$\Ge$ 
pass to $\Gal$ edges 
in the same $\Dgroup$-orbit (up
to direction) if and only if they have the same  
label.
\Qed
}}





\section{Solving the polynomial recursion}
To give an indication of the
nature of the data set, note that 
a table of the first few \pKLp s 
is encoded in Figure~\ref{fig: big pKL} 
(these first few may even be computed by brute force if desired). 
Now we solve the recursion in closed form.

\subsection{
Hypercubes revisited} 

As we have noted in Theorem~\ref{th:crux}, 
the right-action of $\Dgroup$ takes a particularly simple form 
when between  `dominant' elements,
i.e. between elements expressable as $\lambdaa \subset \N$.
We define $\aop{\ai} \lambdaa  = \lambdaa  (\ai,\ai+1)$ 
to be this action between dominant elements.
{\em 
I.e. only for the appropriate domain. }
(Because the underlying descending sequences consist first of positive
terms of descending magnitude, and then negative terms of ascending
magnitude, we call  $\lambdaa \subset \N$ a \valley\ set, and 
$\aop{\ai}$ a \valley\ edge operator.)

\mdef
We generalise the set of \valley\ edge operators $\aop{\ai}$ as follows.

Operator $\aop{ij}$ has action defined in case one of $i,j$ is in
$\lambdaa$, and swaps it for the other (i.e. 
swaps the side of the valley that each of $i,j$ are on).
\\
Example
\[
\aop{36} 56 = 35
\]
(Thus $\aop{\ai} = \aop{\ai \;\; \ai\! +\! 1 }$.
NB, Throughout this section we shall continue to write simply 
$\ai\lambdaa$ for 
$\aop{\ai \;\; \ai\! +\! 1 } \lambdaa$
where no ambiguity arises.)
\\
Where defined, 
each such operator acts involutively; and, where defined, takes
$\lambdaa$ to $\aop{ij} \lambdaa$ comparable to $\lambdaa$ in the $\Ge$
order. 
\\
Each such operator has the same effect on the given $\lambdaa$ as some 
(strictly descending (or ascending))
sequence of $\aop{\ai}$ edge operators. In our example
\[
56 \stackrel{4}{\rightarrow} 46  
   \stackrel{3}{\rightarrow} 36
  \stackrel{5}{\rightarrow} 35 
\]

Operator  
$\aopu{ij}$ has action defined in case 
both or neither of $i,j$ are in $\lambdaa$, and toggles this state. 
\\
Example
\[
\aopu{16} 1456 = 45
\]
which expands, for example, as
\[
1456 \stackrel{3}{\rightarrow} 1356  
   \stackrel{2}{\rightarrow} 1256
  \stackrel{4}{\rightarrow} 1246
   \stackrel{5}{\rightarrow} 1245
  \stackrel{12}{\rightarrow} 45
\]


{\mrem{
Let $v$ be the fully-regular (FR) image of 
$\lambdaa \in \Pow(\N)$ such that
$\aop{ij}\lambdaa$ is defined. 
Unless $j=i+1$ it does {\em not} follow that the fully-regular image of 
$\aop{ij}\lambdaa$ is 
given by the right-action of $(i,j)$ on $v$.
Note, for example, that the underlying 
descending
sequence of
$\aop{ij} \lambdaa$ is {\em not} 
in general a pair permutation of that of $\lambdaa$. 
}}


\mdef  
Let $S$ be a set of generalised \valley\ edge labels,
and $\lambdaa \in \Pow(\N)$. 
If for each subset 
$S' \subseteq S$ the elements of $S'$ may be applied to $\lambdaa$ in
any order to obtain the same set,
and this set lies below $\lambdaa$ in $\Ge$,
then the {\em dominant hypercube} $hh(\lambdaa,S)$ is the digraph
consisting of this collection of sets (vertices) and edges.  



\mdef In Section~\ref{ss:hyperDG} we defined 
a map $\bb: \Pow(\N) \rightarrow \{0,1 \}^{\N} $ and
a map $\Td$ from binary sequences to \TLd s.
It will be convenient to write ${\cal T}(\lambdaa)=\Td(\bb(\lambdaa))$.
We also defined  $\Gamma_{\lambdaa}$ and $\Gamma^{\lambdaa}$
(for $\lambdaa \in \Pow(\N)$, note).
\fotnote{TO GO:
Consider the set of pairs of numbers involved in the 
$\alpha$-operations associated to the edges in the 
lower hypercubical closure of $\lambdaa$ above. 
One can show from the construction  that this set may be identified with
$\Gamma_{\lambdaa}$. 
(I.e. if $(\lambdaa, \aop{\alpha}\lambdaa)$ 
or  $(\lambdaa, \aop{\onetwo}\lambdaa)$ 
is one of these edges then 
$\{ \alpha, \alpha+1 \} \subset  \Gamma_{\lambdaa}$.)
}
By construction we have
\[
h^{\lambdaa} = hh( \lambdaa , \Gamma^{\lambdaa} )
\]
(with the understanding that if $\{ i,j\}$ appears in
$\Gamma^{\lambdaa}$ and is a subset of $\lambdaa$ then the edge
operator is $\aop{\underline{ij}}$).

\fotnote{ NEEDED?:

\mdef
Let 
$$
[m,n] := \{ m , m+1, ... , n \} \subseteq \Z
$$ 

Suppose that $x$ is the largest integer appearing in 
$\lambda \subset \N$. 
We call the subset of $\N$ up to $x$ the $\lambda$-interval,
$[\lambda]$. We call $\lambda$ an $x$-sequence. 

We will often consider the set of pairs `flattened', i.e. just as a
sequence of numbers. In this sense 
note that  $\Gamma_{\lambda}$ is of even degree,
and a subset of $[\lambda]$.}

\begin{figure}
\[
\xymatrix{
& 1 \; 7\; 8\; 10\; 11
\ar[dl]^{67} 
\ar[d]_{9\; 10} 
\ar[dr]_{(5 \; 8)}
\ar[drrr]^{(4 \; 11)}
\\
1 \; 6\; 8\; 10\; 11 \ar[d]\ar[dr]\ar[drrr]& 
1 \; 7\; 8\; 9\; 11  
 \ar[dl] \ar[dr]
\ar[drrr]
&
1 \; 5\; 7\;  10\; 11 \ar[d] \ar[dl] \ar[drrr] &&
1 \; 4\; 7\; 8\; 10 \ar[d] \ar[dl]\ar[dr]
\\
1 \; 6\; 8\; 9\; 11 \ar[dr] \ar[drrr]& 
1 \;5\; 6\;  10\; 11 \ar[d] \ar[drrr]&
1 \;5\; 7\;  9 \; 11 \ar[dl] \ar[drrr]& 
1 \; 4\; 6\; 8\; 10 \ar[d]\ar[dr]&
1 \; 4\; 7\; 8\; 9 \ar[dl] \ar[dr]& 
1 \; 4\; 5\; 7\; 10 \ar[d]\ar[dl]
\\
& 1 \; 5\;6\;  9\;11 \ar[drrr]&&
1 \; 4\; 6\; 8\; 9 \ar[dr]&
1 \; 4\;5\; 6\;  10 \ar[d] &
1 \; 4\; 5\; 7\; 9 \ar[dl]
\\
&&&&
1 \; 4\;5\; 6\;  9&
}
\]
\caption{\label{fig: big h} The hypercube $h^{1 \; 7\; 8\; 10\; 11}$.}
\end{figure}
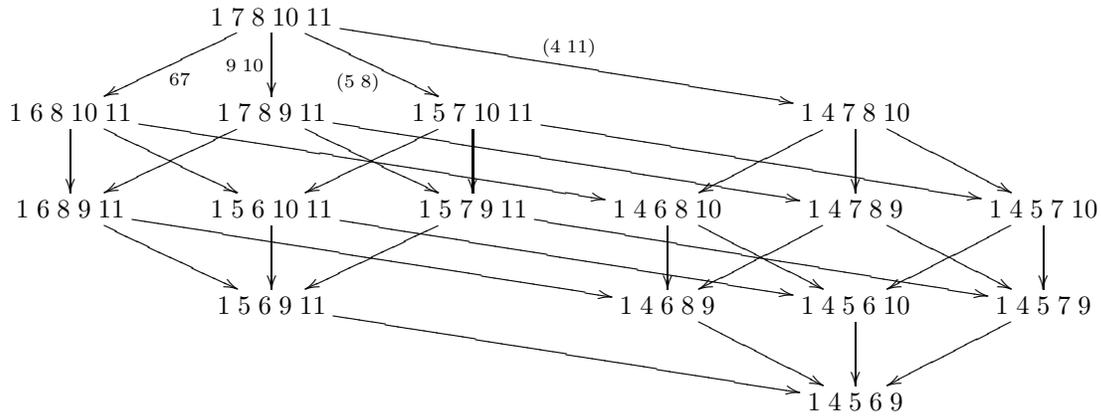

See Figure~\ref{fig: big h} for an example.


\newcommand{\aii}{\ai\! +\! 1}
\newcommand{\aiii}{\ai, \ai\! +\! 1}

{\lem{
\label{lem:TL1}
Suppose $\{ \aiii \} \in \Gamma_{\lambdaa}$ (so $\aop{\ai}\lambdaa < \lambdaa$). 
Let $\{ \ai \} \cup X$, $\{\aii \} \cup Y$ be parts in 
${\cal T}(\aop{\ai}\lambdaa)$ 
($X,Y$ could contain a vertex or be empty). 
Then ${\cal T}(\lambdaa)$ differs from ${\cal T}(\aop{\ai}\lambdaa)$ 
in that  these parts are replaced by 
$\{ \ai, \aii \}$, $X \cup Y$ 
($X\cup Y$ may be empty).
}}
\medskip \\
{\em Proof:}
It is clear that $\{ \ai, \aii \}$ is in ${\cal T}(\lambdaa)$,
so it remains to consider $X,Y$; and to show that all other pairs
agree between  ${\cal T}(\lambdaa)$ and ${\cal T}(\aop{\ai}\lambdaa)$ . 
\\
If $X \cup Y = \emptyset$ then $\ai,\aii$ singletons in
$\aop{\ai}\lambdaa$ and there are no pairs bridging over them,
so no other pair is changed between $\aop{\ai}\lambdaa$ and $\lambdaa$.
\\
If $X=\{i \}, \; Y=\{ j \}$ say, then $j \in \aop{\ai}\lambdaa$ 
(since $\aii \not\in \aop{\ai}\lambdaa$ by construction). 
Suppose $j>\aii$ and $i < \ai$. Then we are in a situation like
\[
\includegraphics[width=4in]{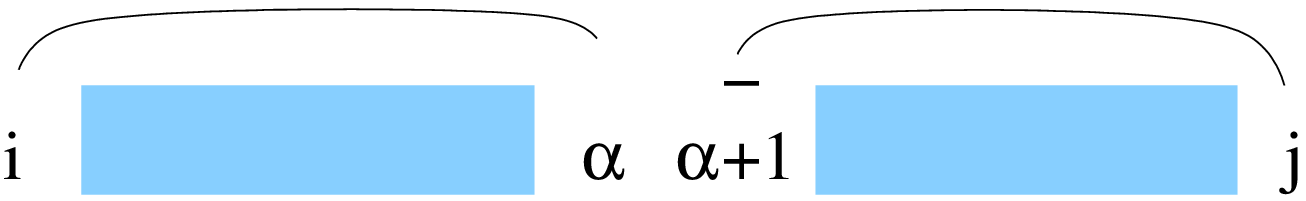}
\]
By construction 
there are no 11 pairs in the $i,\ai$ or $\aii,j$ intervals. 
The algorithm for extracting the sequences in the shaded regions will
thus operate in the same way for each sequence.
In $\lambdaa$ the algorithm generates a pair at $\ai,\aii$
as already noted,
so we may pass to an iteration where these and both shaded parts have
been dealt with.
Vertex $i$ is not involved in a pair from below (else it would be in
$\aop{\ai}\lambdaa$), and $j \in \lambdaa$, so we get a pair $\{i,j \}$ as
required. \\
Suppose $j>\aii$ and $i > j$. Then we are in a situation like
\[
\includegraphics[width=4in]{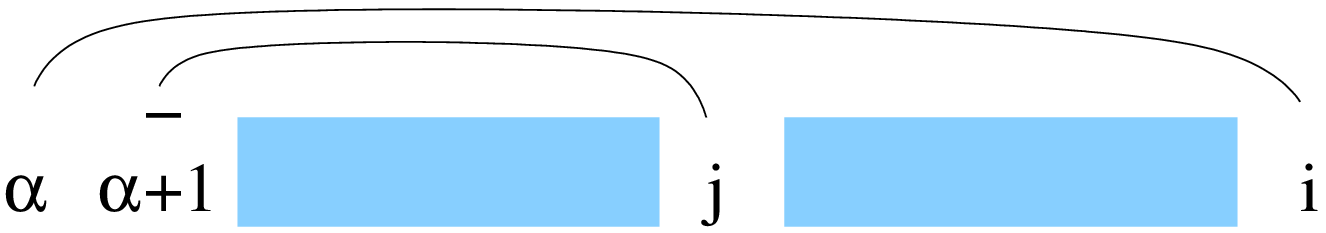}
\]
The same argument goes through until noting that 
$\ai,i \in \aop{\ai}\lambdaa$, so that there is an even number of 1s in the
remainder sequence (algorithm stage 5) left of $\ai$.
This even property
still holds for $\lambdaa$, so $j$ is not involved in a pair from below.
Again we have the required outcome.
\\
The other cases are similar.  
\Qed

\footnot{{OLD COAST STUFF:
\mdef
Suppose $S$ is a subsequence of some $I=[m,n]$ (such as $\lambda$ or 
 $\Gamma_{\lambda}$, both subsequences of $I=[1,x]$). 
We call a maximal consecutive sequence in $S$ 
 an $I$-ISLAND 
(or just an island, if $I$ is clear)
if the next number on either side both lie in $I$. 
We call this pair the COAST of the island. 
\\
We call a number a left-coast if it the next number in $I$ below any
maximal
consecutive sequence in $S$ (not necessarily one with a right coast). 

Let $S$ be given, and let $\lambda$ be another sequence. 
A $\lambda$-right coast is a coast that intersects 
$\lambda$ only in the right-hand element of the coast pair.
A $\lambda$-full coast is a coast that intersects 
$\lambda$ in both elements of the coast pair.

\mdef
Note that  $\Gamma_{\lambda}$ is built from the set of left-coasts $\alpha$ of
$\lambda$ by collecting together all the pairs 
$\{ \alpha, \;\; \alpha+1 \}$;
together with a slight variant for 12: if 12 appears in $\lambda$
then the edge/pair 12 appears in  $\Gamma_{\lambda}$ 
(as opposed to $\alpha=1$ with $\alpha+1=2$, as would appear 
for $\lambda=2346$, say).
\\
Note that we will write $\underline{12}$ for the former, 
where confusion between the
cases just mentioned might otherwise arise. 


{\rem{
We claim that if there is a pair of islands in  $\Gamma_{\lambda}$
separated by a single
integer then 

COME BACK AND FIX THIS!!

the lower of the two overlapping coast operators takes
$\lambda$ to a vertex not below $\lambda$.

Example: 23 and 56 are defined on  7-sequence 1367, so 
$\Gamma_{1367} = \{2,3,5,6 \}$. The coasts are 14 and 47. We have 
$(14) 1367 = 3467$ (not below 1367) and $(47) 1367  = 1346$. 
}}


\mdef
For each $\lambda$ define a set of pairs $\Gamma_{\lambda}'$ 
iteratively as follows.
\\ (1) Initialise  $\Gamma_{\lambda}' = \emptyset$ and $S=\Gamma_{\lambda}$.
\\ (2) If  $S$ has an island 
with a $\lambda$-right coast, 
add the coast for the highest
such to  $\Gamma_{\lambda}'$;
\\
else if  $S$ has an island 
with a $\lambda$-full coast, 
add the coast for the 
lowest 
such to  $\Gamma_{\lambda}'$; 
\\
else if $\lambda\setminus S$ is an adjacent pair, 
add this pair to $\Gamma_{\lambda}'$; 
\\ 
else  stop. 
Add  $\Gamma_{\lambda}'$ to $S$ and repeat (2). 

\mdef
REMARKS:
1. this terminates. \\
2. The order of adding right coasts can be changed without affecting
the outcome.  BULLSHIT!! \\
3. An example of the adjacent pair scenario is 1234 \\
4. cases 3 4 6 7 and 2 7 8 10 11 12 conflict with each other and 
the algorithm at present!!! OUCH!!!

}}

\newcommand{\hhh}[1]{h^{#1}}

\subsection{\KL\ polynomial Theorem}
\newcommand{\taua}{b}
\newcommand{\mua}{c} 
\newcommand{\nua}{c'} 


We continue to use  
labels $\lambdaa\subset\N$ for alcoves.
Thus the rows (and columns)
of the \pKLp\ array $\PH(\Dgroup/\Dgroup^+)$ may be
indexed by these labels. 
That is, there is a polynomial $\ph_a(b) = \ph_{a,b}$,
in the formal variable $v$,
 for each pair $a,b \in \Pow_{even}(\N)$.  
We write $p_a = \{ p_{a,b} \}_{b \in \Pow_{even}(\N)}$ for the
complete row of the array labelled by $a$.

Following on from
(\ref{para:decomp data}) we define polynomial $h^{a}_{b}$ by 
$h^{a}_{b} = v^i$ if $b$ appears in hypercube $h^{a}$ at depth $i$; and  
$h^{a}_{b} = 0$ if $b$ does not appear in $h^{a}$.

{\theo{ 
Let $\lambdaa , \taua \subset\N$ label  alcoves. 
The hypercube $h^{\lambdaa}$ gives the \pKLp s in the 
row $p_{\lambdaa}$  
as follows: 
\[
\ph_{\lambdaa,\taua}  
           =  h^{\lambdaa}_{\taua}
\]
}}

\newcommand{\pai}{\aop{\ai}} 

\noindent {\em Proof:}
We work by induction on the graph order. 
We can then get the polynomials for $\lambdaa$ by looking at the
polynomials for ${\pai}\lambdaa$, where $\ai $ labels one of the edges
in the `shoulder' of the hypercube 
associated to $\Gamma_{\lambdaa}$. 
Specifically, 
by the definition of $\PH$ in Section~\ref{ss:pKLp def},
Theorem~\ref{th:big edge}, and the inductive assumption
we need to determine
all the dominant $\alpha$ images of vertices in 
$\hhh{{\pai}\lambdaa}$.

For any $\ai$ and $\taua \in \Pow(\N)$ let 
$$
\pai \hhh{\taua} := 
 \{ \pai \mua \; | \; \mua\in \hhh{\taua}; \; \pai\mua\mbox{ defined}\}
$$ 
\fotnote{To make contact with the first part of the paper note that
  our $h^{\lambda}$ is isomorphic to each $h_{\QQ}(\mu)$
such that $e_{\QQ}(\mu)=\lambda$, and  
$\ai h^{\lambda}$ is the image of the set of vertices of 
$\ai h_{\QQ}(\mu)$ as defined in (\ref{de:ai act}). 
--- CHECK THIS!!!
}
For example $\pai h^{\pai\lambdaa} \ni \lambdaa$ since
$\pai\pai\lambdaa = \lambdaa$.
Similarly let $\pai^2 \hhh{\taua} =  
\{ \mua \; | \; \mua\in \hhh{\taua}; \; \pai\mua\mbox{ defined}\}$.

Note that there is a map 
$\pai \hhh{\pai\lambdaa} \rightarrow \pai^2 \hhh{\pai\lambdaa}$
given by $\mua \mapsto \pai\mua$, and that this is a bijection
between disjoint sets.

By Section~\ref{ss:pKLp def} (equation(\ref{eq:pKL recursion})) 
and Theorem~\ref{th:big edge}
an alcove label $\mua$ appears in $\ph_{\lambdaa}$ 
(i.e. polynomial $\ph_{\lambdaa,\mua} \neq 0$)
if there is a $\nua$ in  
$\ph_{\pai \lambdaa}$ that, as a vertex of $\Ge$, has an edge labelled
$\ai$ attached to it, and either $\mua = \nua$ or $\mua = \pai \nua$ 
(strictly speaking there is a subtraction to perform
after equation(\ref{eq:pKL recursion}), but we shall see
that all such are null).
The vertices of $\ph_{\lambdaa}$ will thus be  those occuring in
$\pai^2 \hhh{\pai\lambdaa} \cup \pai \hhh{\pai\lambdaa}$,
i.e. as a vertex set:
\[
\ph_{\lambdaa} \; \sim \; 
     \pai^2 \hhh{\pai\lambdaa} \cup \pai \hhh{\pai\lambdaa}
\]
Note that  by the bijection and the inductive hypothesis 
every  alcove label 
appears in at most one way, and hence that every
polynomial will be of form $v^i$. 

We need to check that this set of vertices 
$\ph_{\lambdaa}$
agrees with those of
$\hhh{\lambdaa}$, and that they aquire the right powers via this
identification.  

For any $\taua$ define
\[
\Gamma^{\taua}\setminus \ai = \Gamma^{\taua}\setminus \{ \ai,\aii \}
\]
\[
\Gamma^{\taua}(\ai) = 
\{ \{i,j\} \in \Gamma^{\taua} \; | \; \{i,j\}\cap\{ \ai,\aii \}=\emptyset\}
\]

\begin{figure}
\[
\hhh{\lambdaa}=\;\;\;\;\;
\raisebox{-1.199in}{\includegraphics[width=2in]{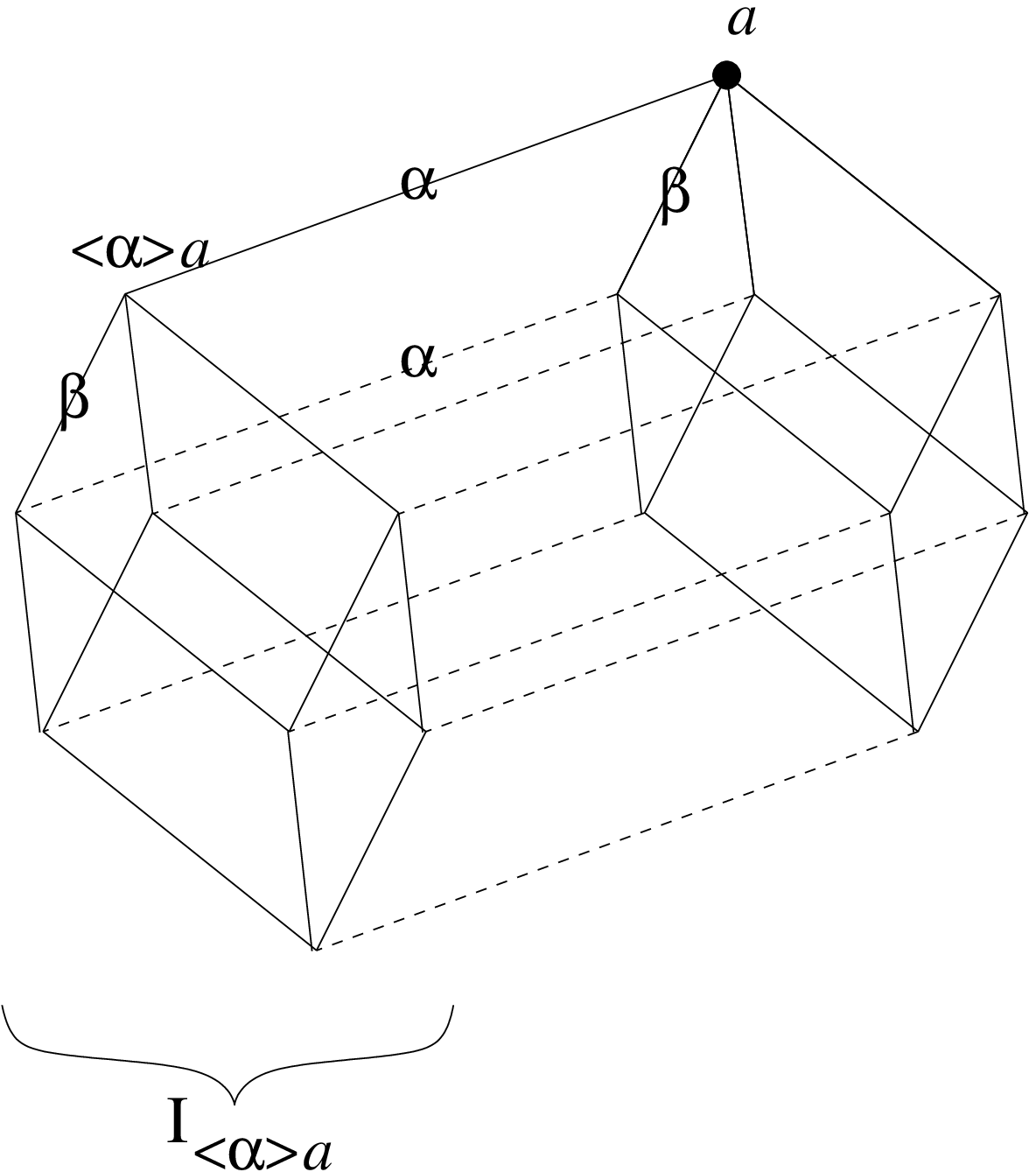}}
\] 
\caption{\label{fig:hypercub1a}}
\end{figure}

Consider the `ideal' $I_{\pai\lambdaa}$ with vertices  
$\mua \leq \pai\lambdaa$ in hypercube $\hhh{\lambdaa}$. 
Note that this sub-hypercube has shoulder
$\Gamma^{\lambdaa}\setminus \ai$; that is
\eql(use1)
I_{\pai\lambdaa} = hh(\pai\lambdaa, \Gamma^{\lambdaa}\setminus \ai)
\eq
and that the quotient of $\hhh{\lambdaa}$ 
by this ideal has the same shoulder set.
Note also that this quotient $\hhh{\lambdaa}/I_{\pai\lambdaa}$ consists 
of the images under $\ai$ of the vertices in $I_{\pai\lambdaa}$,
as exemplified in Figure~\ref{fig:hypercub1a}.

It follows from Lemma~\ref{lem:TL1}
that $\Gamma^{\lambdaa}\setminus \ai$
agrees with the set $\Gamma^{\pai\lambdaa}(\ai)$ of pairs in
$\Gamma^{\pai\lambdaa}$ that do not intersect $\ai$ or $\ai+1$,
{\em except} that if there are pairs  $\ai,i$ and $\ai\! +\! 1,j$
in $\Gamma^{\pai\lambdaa}$ then there will be a pair $i,j$ in
$\Gamma^{\lambdaa}\setminus \ai$
(that obviously does not appear in $\Gamma^{\pai\lambdaa}$):
\eql(use2)
\Gamma^{\pai\lambdaa}(\ai) \; \subseteq \; \Gamma^{\lambdaa}\setminus \ai
\eq
From (\ref{use1}) and (\ref{use2})  
we have that $hh(\pai\lambdaa, \Gamma^{\pai\lambdaa}(\ai))$
is a subgraph of $I_{\pai\lambdaa}$ and hence of $\hhh{\lambdaa}$
(albeit one layer down from the `head'), and also of 
$\hhh{\pai\lambdaa}$. 


As noted, all the vertices in the subgraph
$hh(\pai\lambdaa, \Gamma^{\pai\lambdaa}(\ai))$ of $\hhh{\pai\lambdaa}$  
have $\ai$-images (and these images are above in the graph order).
Thus all these vertices and images appear in $\ph_{\lambdaa}$
(by the inductive assumption $\ph_{\pai\lambdaa} \equiv \hhh{\pai\lambdaa}$
and the constructive definition of $\ph_{\lambdaa}$ from $\ph_{\pai\lambdaa}$).
The power of $v$ for each image vertex is inherited from the original
vertex (for example 
$\ph_{\lambdaa}(\lambdaa)  = \ph_{\lambdaa}(\pai\pai\lambdaa) 
   = \ph_{\pai\lambdaa}(\pai\lambdaa)=v^0$), 
while the power of $v$ for the original vertex is raised by
1 (example: $ \ph_{\lambdaa}(\pai\lambdaa) 
=v \ph_{\pai\lambdaa}(\pai\lambdaa)=v v^0 = v^1$). 
We see, therefore, that 
{\em all these vertices have the correct exponent}.
\\
The other vertices in the shoulder of
$\hhh{\pai\lambdaa}$ 
(the ones, if any, at the end of edges of form  
$\ai,i$ and $\ai\! +\! 1,j$)
do not have $\ai$-images.
Thus we have agreement between $\hhh{\lambdaa}$ and 
$\ph_{\lambdaa} \sim \pai^2 \hhh{\pai\lambdaa} \cup \pai \hhh{\pai\lambdaa}$
except for the ideal generated by $\aop{ij}\lambdaa$ as above (if any) in
$\hhh{\lambdaa}$ on the one hand;
and the possible descendents of 
$\aop{\ai,i} \pai\lambdaa$ and $\aop{\ai\! +\! 1,j}\pai\lambdaa$
in $\hhh{\pai\lambdaa}$ that {\em do} have $\ai$-images on the other. 


It remains to show that these contributions match up
(with the correct powers). 

\begin{figure}
\includegraphics[width=5in]{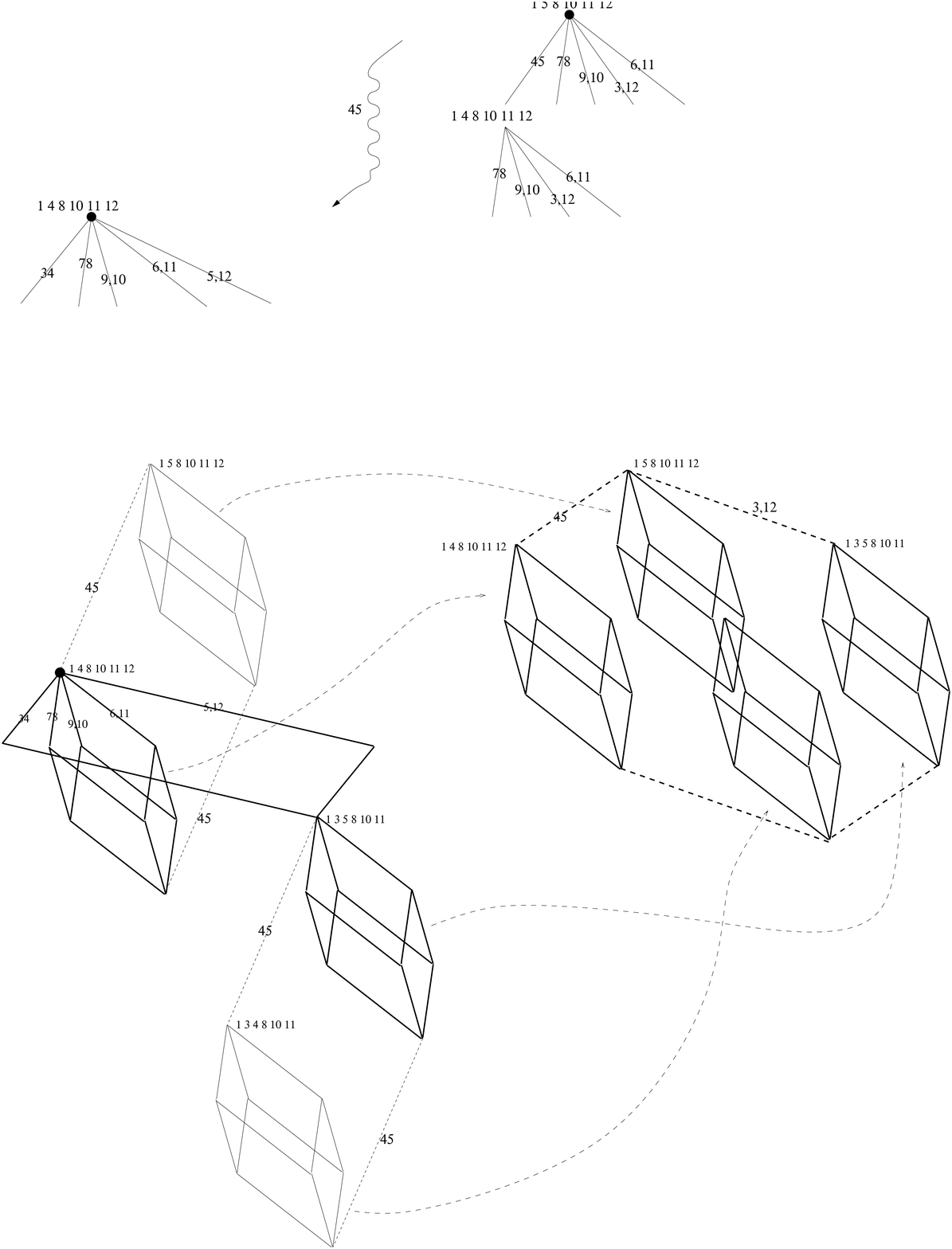}
\caption{\label{fig:eg158} A representative example. 
The top part of the figure shows the shoulder of $\hhh{\lambdaa}$ for 
$\lambdaa=\{ 1 , 5 , 8 , 10 , 11 , 12\}$; and the shoulder of the
ideal below $\aop{45}\lambdaa = \{ 1 , 4 , 8 , 10 , 11 , 12\}$ 
within $\hhh{\lambdaa}$.
Immediately below-left of this is the shoulder of $\hhh{\aop{45}\lambdaa}$
itself (showing that this hypercube is bigger than the corresponding
ideal within $\hhh{\lambdaa}$). 
The bottom-left part of the figure shows all of the vertices in 
$\hhh{\aop{45}\lambdaa}$ that have $\aop{45}$ images
(and a couple which do not, that are relevant for the construction);
together with a representation of those images.
The bottom-right part shows how all these vertices may be collected
together to constitute the vertices of  $\hhh{\lambdaa}$.
}
\end{figure}

If there is no such $\aop{ij}\lambdaa$ then one can show that there are
not descendents of 
$\aop{\ai,i} \pai\lambdaa$ and $\aop{\ai\! +\! 1,j}\pai\lambdaa$
in $\hhh{\pai\lambdaa}$ with $\ai$-images and we are done.
So let us suppose there is $\aop{ij}\lambdaa$ in $\hhh{\lambdaa}$.
Note that for our $\lambdaa$ we have
\eql(eq:cool)
\aop{ij}\lambdaa = \aop{\ai,i} \aop{\ai\! +\! 1,j}\pai\lambdaa
\eq
See Figure~\ref{fig:eg158} for a representative example of this. 
We have there 
\[
\lambdaa = \; 
1 \; 5 \; 8 \; 10 \; 11 \; 12
\stackrel{\aop{45}}{\rightarrow}
1 \; 4 \; 8 \; 10 \; 11 \; 12
\stackrel{\aop{5 \; 12}}{\rightarrow}
1 \; 4 \; 5 \; 8 \; 10 \; 11 
\stackrel{\aop{3 \; 4}}{\rightarrow}
1 \; 3 \; 5 \; 8 \; 10 \; 11 
\; = 
\aop{3 \; 12} \lambdaa
\]
A similar version works for $\aopu{ij}$ operators. 


The $\aop{ij}\lambdaa$ in $\hhh{\pai\lambdaa}$ is in level $i=2$ by 
(\ref{eq:cool}), 
and has a hypercube 
$hh(\aop{ij}\lambdaa,\Gamma^{\pai\lambdaa}(\ai))$ below it. 
All the elements of this hypercube have 
$\ai$-images, since $\pai,\aop{ij}$ commute.
Note for example that  $\aop{ij}\lambdaa$ itself
has an $\ai$-image 
(although $\aop{ij}\lambdaa$ is below $\aop{\ai\! +\! 1,j}\pai\lambdaa$,
which does not have an $\ai$-image, in the {\em graph} order), 
and that its $\ai$-image $\pai\aop{ij}\lambdaa$ is 
{\em  below} it in the graph order. 
The other labels in the ideal behave similarly.
Thus the polynomials asigned by Equation(\ref{eq:pKL recursion}) to 
the relevant part of 
$\ph_{\lambdaa} \sim \pai^2 \hhh{\pai\lambdaa} \cup \pai \hhh{\pai\lambdaa}$
are, for $v^i$ the relevant polynomial from $\ph_{\pai\lambdaa}$,
$v^i$ (for the $\ai$-image) and $v^{i-1}$ 
(the vertex `left behind') respectively.
The $-1$ compensates for the fact that the vertex appears in 
$\hhh{\pai\lambdaa}$ one layer lower than
in $\hhh{\lambdaa}$ (where it appears in the shoulder in the case of 
 $\aop{ij}\lambdaa$ itself for example), 
so subject to the 
working assumptions we verify 
$\ph_{\lambdaa} \equiv \hhh{\lambdaa}$. 

Note finally that this $-1$ increment only occurs for the vertex 
 $\aop{ij}\lambdaa$ and
those below it, and thus for polynomials $v^i$ with exponent $i\geq 2$. 
Thus we never have an increment of form 
$v^1 \rightarrow v^{1-1}=v^0$ (which would incur a subtraction in the
polynomial construction). 
The only remaining working assumption is the inductive assumption.

\Qed





\noindent
{\bf Concluding remarks}.
As already noted, a 
significant mathematical application of this work is hoped to be as
a base for corresponding calculations over fields of finite
characteristic (cf. \cite[\S 6]{\CDMii}). 
A {\em physically} motivated application is in computing eigenvectors
of the Young matrix (the adjacency matrix of the Young graph
\cite{Kerov03}),
which are involved in certain quantum spin chain computations 
(see e.g. \cite{CanduSaleur08}).
We note that formal connections between \pKLp s and Brauer algebra
decomposition matrices can be constructed in principle by 
other approaches,
such as  in \cite{OrellanaRam01}. 
However such 
formal approaches do not give 
access to the specific decomposition numbers that we compute here 
(and which are  required for the applications mentioned). 
Finally we note that \cite{EnrightShelton87} includes
formulations of Kazhdan--Lusztig polynomials  
related to the $\Dgroup/\Dgroupp$ 
case, 
considered from an entirely different perspective.



\noindent {\bf Acknowledgements}. 
I thank Anton Cox and  Maud De Visscher
for many useful discussions and suggestions;
and  
Robert Marsh for useful conversations. I thank Paula Martin for
support. 







\bibliographystyle{amsplain}
\bibliography{new31}

\end{document}